\documentclass[11pt]{article}
\usepackage{amsfonts,latexsym}
\newcommand{\cleqn}{\setcounter{equation}{0}}
\newcommand{\clth}{\setcounter{theorem}{0}}
\newcommand {\sectionnew}[1]{\section{#1}\cleqn\clth}
\newcommand{\beq}{\begin{equation}}
\newcommand{\eeq}{\end{equation}}
\newcommand{\beqa}{\begin{eqnarray}}
\newcommand{\eeqa}{\end{eqnarray}}
\newcommand{\beaa}{\begin{eqnarray*}}
\newcommand{\ben}{\begin{eqnarray*}}
\newcommand{\eaa}{\end{eqnarray*}}
\newcommand{\een}{\end{eqnarray*}}

\newcommand{\text}{\textrm}
\newcommand \nc {\newcommand}
\nc \proof {\noindent {\em{Proof.\/ }}} 
\nc \qed {$\Box$\hfill}
\newtheorem{theorem}{Theorem}[section]
\newtheorem{lemma}[theorem]{Lemma}
\newtheorem{proposition}[theorem]{Proposition}
\newtheorem{corollary}[theorem]{Corollary}
\newtheorem{definition}[theorem]{Definition}
\newtheorem{example}[theorem]{Example}
\newtheorem{remark}[theorem]{Remark}
\newtheorem{conjecture}[theorem]{Conjecture}
\newtheorem{question}[theorem]{Question}
\nc \bth[1] {\begin{theorem}\label{t#1} }
\nc \ble[1] {\begin{lemma}\label{l#1} }
\nc \bpr[1] {\begin{proposition}\label{p#1} }
\nc \bco[1] {\begin{corollary}\label{c#1} }
\nc \bde[1] {\begin{definition}\label{d#1}\rm }
\nc \bex[1] {\begin{example}\label{e#1}\rm }
\nc \bre[1] {\begin{remark}\label{r#1}\rm }
\nc \bcon[1] {\begin{conjecture}\label{con#1}\rm }
\nc \bque[1] {\begin{question}\label{que#1}\rm }
\nc {\eth} { \end{theorem} }
\nc {\ele} { \end{lemma} }
\nc {\epr}{ \end{proposition} }
\nc {\eco} { \end{corollary} }
\nc {\ede} {\end{definition} }
\nc {\eex} { \end{example} }
\nc {\ere} {\end{remark} } 
\nc {\econ} { \end{conjecture} }
\nc {\eque} {\end{question} }
\nc \eqref[1] {{\rm{(\ref{#1})}}} 
\nc \thref[1]{Theorem \ref{t#1}}
\nc \leref[1]{Lemma \ref{l#1}} 
\nc \prref[1]{Proposition
\ref{p#1}} \nc \coref[1]{Corollary \ref{c#1}} 
\nc \deref[1]{Definition \ref{d#1}} 
\nc \exref[1]{Example \ref{e#1}}
\nc \reref[1]{Remark \ref{r#1}} 
\nc \conref[1]{Conjecture\ref{con#1}}

\def\a{\alpha}
\def\b{\beta}

\def\Ga{\Gamma}
\def\e{\epsilon}

\def\l{\lambda}
\def\la{\lambda}

\def\si{\sigma}

\def \F {{\mathbb F}}

\def \Z {{\mathbb Z}}
\def \Q {{\mathbb Q}}

\def \C {{\mathbb C}}
\def \Hom{ {\mathrm{Hom}}}
\def \Aut{ {\mathrm{Aut}}}
\def \End{ {\mathrm{End}}}
\def \tr{ {\mathrm{Tr}}}
\def \coker { {\mathrm{Coker}} }

\def \mod { {\mathrm{mod}} }

\def \Im { {\mathrm{Im}} }

\renewcommand \ker { {\mathrm{Ker}} }
\def \vect {\overrightarrow }
\nc \Wr {Wr} \nc \GRN { \Gr^{(N)} }
\nc \GRA[1] { \Gr_A^{(#1)} }   
\nc \GRAN { \GRA{N} } \nc \GrA[1] { \Gr_A(#1) }\nc \GrAa {
\GrA{\alpha} }
\nc \GRB[1] { \Gr_B^{(#1)} }   
\nc \GRBN { \GRB{N} } \nc \GrB[1] { \Gr_B(#1) } \nc \GrBb {
\GrB{\beta} }
\nc \GRMB[1] { \Gr_{MB}^{(#1)} }   
\nc \GRMBN { \GRMB{N} } \nc \GrMB[1] { \Gr_{MB}(#1) } \nc \GrMBb {
\GrMB{\beta} }

\begin{document}
\title{{\LARGE\bf{Euler characteristics of arithmetic groups}}}

\author{
I. ~Horozov
\thanks{E-mail: horozov@math.brown.edu}
\\ \hfill\\ \normalsize \textit{Department of Mathematics,}\\
\normalsize \textit{ Brown University, 151 Thayer st., Providence
RI 02912 }   }
\date{}

\maketitle

\begin{abstract}
We develop a general method for computing the homological
Euler characteristic of finite index subgroups $\Ga$ of
$GL_m({\cal{O}}_K)$ where ${\cal{O}}_K$ is the ring of integers in
a number field $K$. With this method we find, that for large,
explicitly computed dimensions $m$, the homological Euler
characteristic of finite index subgroups of $GL_m({\cal{O}}_K)$
vanishes. For other cases, some of them very important for spaces
of multiple polylogarithms, we compute non-zero homological Euler
characteristic. With the same method we find all the torsion 
elements in $GL_3\Z$ up to conjugation.
Finally, our method allows us to obtain a formula
for the Dedekind zeta function at $-1$ in terms of 
the ideal class set and the multiplicative group of 
quadratic extensions of the base ring.
\end{abstract}
\setcounter{section}{-1}
\sectionnew{Introduction}

The homological Euler characteristic of a group $\Ga$ with
coefficients in a representation $V$ is defined by
$$\chi_h(\Ga,V)=\sum_i (-1)^i\dim H^i(\Ga,V)$$ (see \cite{S},
\cite{B1}). The main theoretical result of our thesis is a general
method that allows us to calculate $\chi_h(\Ga,V)$. Toward the end
of the introduction we describe briefly this method. Before that
we list the most important results of the thesis which demonstrate
the scope of the method.
Our first result is about vanishing of
homological Euler characteristics.


\bth{0.1}
 Let $\Ga$ be a finite index subgroup of $GL_m
({\cal{O}}_K)$, where ${\cal{O}}_K$ is the ring of integers in a
number field $K$. Let $V$ be a finite dimensional representation
of $\Ga$. Then\\

(a) For $K=\Q$ if $m>10$ then $\chi_h(\Ga,V)=0,$\\

(b) For $K=\Q[\sqrt{-d}]$, where $-d$ is the discriminant,\\
   if $d=4$ and $m>4$ then $\chi_h(\Ga,V)=0$;\\
   if $d=3$ and $m>6$ then $\chi_h(\Ga,V)=0$;\\
  for the other $d$'s and $m>2$ we have $\chi_h(\Ga,V)=0,$\\

(c) For the remaining number fields $K$ we always have
$\chi_h(\Ga,V)=0.$
 \eth
We obtain a similar result for $SL_m$ in place of $GL_m$.

\bth{0.2}
 Let $\Ga$ be a finite index subgroup of $SL_m
({\cal{O}}_K)$, where $m\geq 2$ and ${\cal{O}}_K$ is the ring of
integers in a number field $K$. Let $V$ be a finite dimensional
representation of $\Ga$. Then\\

(a) For $K=\Q$ if $m>10$ then $\chi_h(\Ga,V)=0,$\\

(b) For $K=\Q[\sqrt{-d}]$, where $-d$ is the discriminant,\\
   if $d=4$ and $m>4$ then $\chi_h(\Ga,V)=0$;\\
   if $d=3$ and $m>6$ then $\chi_h(\Ga,V)=0$;\\
   for the other $d$'s and $m>2$ we have $\chi_h(\Ga,V)=0,$\\

(c) For $K$ totally real field if $m>2$ then
$\chi_h(\Ga,V)=0,$\\

(d) For the remaining number fields $K$ we always have
$\chi_h(\Ga,V)=0.$
 \eth
Parts (a) and (b) follow from the previous theorem while parts
are new statements.

We also compute the homological Euler characteristic of the
arithmetic subgroups $\Ga_1(3,N)$ and $\Ga_1(4,N)$ of $GL_3(\Z)$
and $GL_4(\Z)$, respectively, where $\Ga_1(m,N)$ is the subgroup
of $GL_m(\Z)$ that fixes the vector $[0,\dots,0,1]\mbox{ } \mod
\mbox { }N$.
\bth{0.3} The homological Euler characteristic of $\Ga_1(3,N)$ and
of $\Ga_1(4,N)$ for $N$ not divisible by $2$ and $3$ is given by

$$\chi_h(\Ga_1(3,N),\Q)=-\frac{1}{12}\varphi_2(N)+\frac{1}{2}
\varphi(N),$$

$$\chi_h(\Ga_1(4,N),\Q)=\varphi(N),$$
   where $\varphi(N)$ is the
Euler $\varphi$-function, and $\varphi_2(N)$ is the multiplicative
arithmetic function generated by
$\varphi_2(p^a)=p^{2a}(1-\frac{1}{p^2})$ for $a \geq 1$.
     \eth

With the same technique we find all the torsion elements in $GL_3\Z$
up to conjugation (see proposition 3.2).
Using our method we, also, compute
$\chi_h(GL_m(\Z),S^nV_m)$ for $m=3$ and $m=4$. The computation of
$\chi_h(GL_3(\Z),S^nV_3)$ (see theorem 6.4) agrees with the
computation of the $H^i(GL_3(\Z),S^nV_3)$ in \cite{G1}, which was
used for computation of dimensions of spaces of certain multiple
polylogatithms. Also the Euler characteristic of $\Ga_1(3,N)$ with
trivial coefficients when $N$ is an odd prime, greater that 3,
agrees with the computation of $H_{inf}^i(\Ga_1(3,N),\Q)$ in the
corrected version of \cite{G1}. We compute the homological Euler
characteristic of $GL_2(\Z[i])$ and $GL_2(\Z[\xi_3])$ with
coefficients in the symmetric powers of the standard
representations. Also we compute the homological Euler
characteristic of $\Ga_1(2,\mathfrak{a})$ for an ideal
$\mathfrak{a}$ in $Z[i]$ and $\Z[\xi_3]$, respectively. These
computations can be used for finding dimensions of spaces of
certain elliptic polylogarithms, as it was explained to me by
professor A. Goncharov.
\bth{6.12}
(a) If $1+i$ does not divide ${\mathfrak{a}}$
the homological Euler characteristic of 
$\Ga_1(2,{\mathfrak{a}})\subset GL_2(\Z[i])$ is given by
$$\chi_h(\Ga_1(2,{\mathfrak{a}}),\Q)=
\frac{1}{2}\varphi_{\Z[i]}({\mathfrak{a}}),$$
where $\varphi_{\Z[i]}({\mathfrak{a}})$ is the multiplicative function defined 
on the ideals of $\Z[i]$, generated by
$$\varphi_{\Z[i]}({\mathfrak{p}}^n)=N_{\Q(i)/\Q}({\mathfrak{p}})^n(1-
\frac{1}{N_{\Q(i)/\Q}({\mathfrak{p}})}).$$

(b) If $1+\xi_6$ does not divide ${\mathfrak{a}}$
the homological Euler characteristic of 
$\Ga_1(2,{\mathfrak{a}})\subset GL_2(\Z[\xi_3])$ is given by
$$\chi_h(\Ga_1(2,{\mathfrak{a}}),\Q)=
\frac{1}{3}\varphi_{\Z[\xi_3]}({\mathfrak{a}}),$$
where $\varphi_{\Z[\xi_3]}({\mathfrak{a}})$ 
is the multiplicative function defined 
on the ideals of $\Z[\xi_3]$, generated by
$$\varphi_{\Z[\xi_3]}({\mathfrak{p}}^n)=N_{\Q(\xi_3)/\Q}({\mathfrak{p}})^n(1-
\frac{1}{N_{\Q(\xi_3)/\Q}({\mathfrak{p}})}).$$
\eth

In general the method works for any arithmetic subgroup of
$GL_m({\cal{O}}_K)$, where ${\cal{O}}_K$ is the ring of integers
in a number field $K$.

Our approach is the following: we
generalize a result of K. Brown \cite{B2} that relates the torsion
elements in the group up to conjugation to the homological Euler
characteristic of the group. Namely,
$$\chi_h(\Ga,V)=\sum_T \chi(C(T))\tr(T^{-1}|V),$$
      where the sum is over
all torsion elements $T$ of $\Ga$ up to conjugation and $C(T)$ is
the centralizer of $T$ in $\Ga$. Let us recall the definition of
orbifold Euler characteristic of $\Ga$, denoted by $\chi(\Ga)$,
which we simply call Euler characteristic. If $\Ga$ is a torsion
free group then $\chi(\Ga)=\chi_h(\Ga)$. If $\Ga$ has torsion
consider a finite index torsion free subgroup $\Ga'$. Then
$\chi(\Ga)$ is defined by $$\chi(\Ga)=[\Ga:\Ga']^{-1}\chi(\Ga').$$
Arithmetic groups do have a finite index torsion free subgroup.
So for them the Euler characteristic is defined. The main properties 
of the Euler characteristic that we are going to use are:
$$\chi(G)=\frac{1}{|G|} \mbox{ for a finite group G},$$
Given an exact sequence 
$$0\rightarrow\Ga_1 \rightarrow \Ga \rightarrow \Ga_2
\rightarrow 0,$$
we have 
$$\chi(\Ga)=\chi(\Ga_1)\chi(\Ga_2),$$
and
$$\chi(SL_m({\cal{O}}_K))=\zeta_K(-1)\dots \zeta_K(1-m).$$
The first two properties can be found in K. Brown's book \cite{B1}.
And the last one is a difficult result due to Harder \cite{H}.

In order to use the Brown's formula we need to know the torsion elements in
the group.
We develop a method for finding the torsion elements in
$GL_m({\cal{O}}_K)$, and consider another form of the above
formula that requires {\it{very few}} of the torsion elements ( in
general they can be quite a large number). 
This method involves linear algebra over number rings. More precisely
is describes a normal form of matrices over number ring. First we deal with
matrices with irreducible characteristic polynomial which lead to a 
relation to ideal classes (proposition 2.1). Then we examine matrices with
reducible characteristic polynomials which leads to resultants 
(corollary 3.5 and proposition 3.7).
Then we examine
the relation between the torsion elements in $GL_m({\cal{O}}_K)$
and the torsion elements in an arithmetic subgroup $\Ga$ of
$GL_m({\cal{O}}_K)$. That gives us the homological Euler
characteristic of $\Ga$ with coefficients in a representation.

We obtain another interesting application of our method for
computation of the values of the Dedeking zeta functions at $-1$.
Now we are going to explain how this is related to our method.
First, the Dedekind zeta function at $-1$ vanishes for number
fields which are not totally real. For totally real number fields
$K$ we consider the arithmetic group $SL_2({\cal{O}}_K)$. By a
result of Harder \cite{H} we have that the orbifold Euler
characteristic of $SL_2({\cal{O}}_K)$ is the Dedekind zeta
function at $-1$, $\zeta_K(-1)$. Also the homological Euler
characteristic is always an integer. Let us denote the 
homological Euler characteristic by an integer $N$, i.e 
$$N=\chi_h(SL_2({\cal{O}}_K)).$$

We also introduce the following notation.
Let $$Cl({\cal{O}}_K[\xi]/{\cal{O}}_K)$$
be the set of ideal classes in ${\cal{O}}_K[\xi]$
which are free as ${\cal{O}}_K$-modules.

\bth{0.4} 
Let
$K$ be a totally real number field. Let $R$ be the
ring of integers in $K$. Then
the Dedekind zeta function at $-1$ can be expressed as
$$\zeta_K(-1)=
-\frac{1}{4}\sum_{\xi}\sum_{I\in Cl({\cal{O}}_K[\xi]/{\cal{O}}_K)} 
\frac{\#{\cal{O}}_K^{\times}/
N_{K(\xi)/K}(R_I^{\times})}{\#(R_I^{\times})_{tors}}
+\frac{1}{2}N,$$
 where the the first sum is taken over all roots of $1$ 
such that $[K(\xi):K]=2$, the second sum over all ideal classes in
${\cal{O}}_K[\xi]$ that are free as ${\cal{O}}_K$-modules, the ring $R_I$
 sits between ${\cal{O}}_K[\xi]$ and its integral closure ${\cal{O}}_{K(\xi)}$.
(for precise definition of $R_I$ see theorem 8.1), and 
$N=\chi_h(SL_2({\cal{O}}_K),\Q)$.
\eth
 The organization of the thesis is the following.
In section 1 we deals with linear algebra over rings of algebraic integers.
We  give a method for classification of the matrices with integer coefficients
for the ones that are diagonalizable over the complex numbers. In section 2 
we give a method of computing the centralizer of a matrix which leads 
to explicit formulas for the homological Euler characteristics of $GL_m\Z$,
$GL_m(\Z[i])$, $GL_m(\Z[\xi_3])$ and $GL_2({\cal{O}}_{\Q(\sqrt{-d})})$,
see respectively theorems 2.10, 2.11, 2.12 and 2.13. In section 3 we
find the torsion elements in $GL_2\Z$ and $GL_3\Z$ up to conjugation.
In section 4 is very computational. We compute many 
resultants and centralizers 
needed for the homological Euler characteristics of various group. 
In section 5 we find the homological Euler characteristics of 
$GL_m({\cal{O}}_K)$ with coefficients the symmetric powers 
of the standard representation. In section 6 we find the homological 
Euler characteristics of arithmetic groups $\Ga_1(m,N)$ and 
$\Ga_1(m,{\mathfrak{a}}).$ In section 7 we examine the groups 
$SL_2({\cal{O}}_K)$ for totally real number fields $K$ which gives a relation 
to the Dedekind zeta function of the field. And in section 8 we prove the 
generalization of Brown's formula.
\begin{flushleft}
{\large\bf{Acknowledgements}}
I would like to thank most of all my advisor professor Alexander Goncharov 
for constant 
encouridgement for my work and for giving me a series of problems that lead 
to this thesis. I would like to thank professor Michael Rosen,
Bruno Harris, and Steven Lichtenbaum for the discussions we had. 
I would like to thank professor Thomas Goodwillie for refering me to
Brown's formula. 
And also I would like to thank Amir Jafary for the many useful 
conversations we had that contributed to this paper.

\end{flushleft}



\sectionnew{Conjugacy classes in $GL_m({\cal{O}}_K)$}

In this section we describe the conjugacy classes of elements in 
$GL_m({\cal{O}}_K)$,
where ${\cal{O}}_K$ is the ring of integers in a number field $K$. 
 We approach the description of conjugacy classes
in the following way. First we deal with matrices in $GL_m K$. Then we examine 
matrices in $GL_m({\cal{O}}_K)$ whose characteristic 
polynomial is irreducible over $K$.
They are described by ideal classes of a larger ring. And, finally, 
we construct an algorithm for
matrices with reducible characteristic polynomial that allows to
consider instead matrices of smaller dimension. Using this
inductive step we can describe completely conjugacy classes of
certain type, knowing the conjugacy classes of matrices of smaller
dimension. 

\subsection{Conjugacy classes in $GL_m K$}
First we recall a few statements about matrices with coefficients in an 
infinite field $K$.

    \bpr{1.1}

    Let $A$ and $B$ be matrices with coefficients in an 
infinite field $K$. If $A$
and $B$ are conjugate to each other inside $GL_m \bar{K}$,
where $\bar{K}$ is the algebraic closure of $K$, then they are
conjugate to each other inside $GL_m K.$
    \epr

    \proof Consider the vector space of matrices $M$
over the fields $\bar{K}$ and $K$ such that
$AM=MB.$ Let $$V_{\bar{K}}=\{M\in Mat_m \bar{K} : AM=MB\}$$ and $$V_{K}=\{M\in
Mat_m K : AM=MB\}.$$ Since $M$ is a solution of a linear system
over fields, we have that $$\dim_{K}V_{K}=dim_{\bar{K}}V_{\bar{K}}.$$ We
know that there is $Q\in GL_m \bar{K}$ such that $AM=MB.$ Equivalently,
there is $Q$ such that $Q\in V_{\bar{K}}$ and $\det Q \neq 0.$ Let
$$det^0_{K}=\{M\in Mat_n K: \det M=0 \}.$$ Then the set
$\det^0_{K}$ is Zariski closed subset of $V_K$ that is not the
entire space, because $K$ has infinitely may elements. 
Let $P\in V_{K}-\det^0_{K}$. Then $AP=PB$, and $P\in
GL_m K.$ Thus, $A$ and $B$ are conjugate to each other as
elements of $GL_n K.$

We are going to use a well known statement.
  \ble{1.2}
    The characteristic polynomial of
    $$\pmatrix{0 & 1  \cr \vdots & 0 & 1 \cr &&& \ddots \cr 0 &&&&1
    \cr -a_0 & -a_1 &&& -a_{m-1}}$$
     is $\la^m + a_{n-1}\la^{m-1} + \dots + a_0.$
    \ele
    For $2\times 2$-matrix can checked by direct computation
The one can use induction on the size of the matrix.\\
    Given a polynomial $f(\la)=\la^m + a_{m-1}\la^{n-1} + \dots +
a_0$, we denote the above matrix by $A_f$.


\subsection{Block-triangular form}


This subsection deals with the case when the characteristic
polynomial of a matrix is reducible over the rational numbers. The
partition of a matrix into blocks will be done in the following
way: Given an $m \times m$ matrix $A$, let $m=m_1 + \dots + m_k$
be a partition of $n.$ Then $A$ can be thought of as a $k \times
k$ block-matrix whose $A_{ij}$-entry, $i,j=1, \dots k$ is a block
(and a matrix) of size $m_i \times m_j.$ This will be the type of
block-matrices that we consider. Note that the blocks $A_{ii}$ are
square matrices.


    \bth {1.5} 
An $m\times m$ matrix $A$ over a principal ideal domain can be conjugated via
$B\in GL_m({\cal{O}}_K)$ so that
$$
BAB^{-1}=
\left[
\begin{tabular}{cccc}
$A_{11}$ & $\dots$  & $A_{1d}$\\
$0$      & $\ddots$ & $\vdots$\\
$0$      & $0$      & $A_{dd}$ 
\end{tabular}
\right],
$$
where $A_{ij}$ are matrices such that $A_{ii}$ has irreducible over $K$
characteristic polynomial and $A_{ij}=0$ for $i>j$.
    \eth

    \proof
The plan for the proof is the following: First we consider the
invariant subspaces over $\overline{K}$, where $\overline K$ is the 
algebraic closure of $K$, in order to find the
invariant spaces over $K$. Then we construct a good basis over R.\\
Let $A$ be  matrix from $GL_m({\cal{O}}_K)$. Let $\l$ be an
eigenvalue of $A$. Let $v_{1,\l}$ be one of the eigenvectors
corresponding to $\l$. Let $v_{i,\l}$, $i=2,\cdots,k$ be adjoint
vectors to $v_{1,\l}$ so that $$(A-\l I)v_{i+1,\l}=v_{i,\l}.$$ We
can choose the coefficients of $v_{1,\l}$ to be from the field
$K(\l)$, since it is determined from the linear system
$$(A-\l I)v_{1,\l}=0.$$ Then similarly, all $v_{i,\l}$'s have
coefficients in the field $K(\l)$. Let $\l_1, \cdots, \l_d$ be
the Galois conjugates to $\l$. And let $$L:=K(\l_1, \cdots,
\l_d).$$ We can consider the Galois conjugates of $v_{1,\l}$. If
$\si \e Gal(L/K)$ then $v_{1,\l}^{\si}=v_{1,\l^{\si}}$ is an
eigenvector with eigenvalue $\l^{\si}$ and adjoint vectors
$v_{i,\l}^{\si}=v_{i,\l^{\si}}$ for $i=2, \cdots, k$. We can
consider invariant subspaces with coefficients in $L$, namely the
vector space over $L$ spanned by the vectors $v_{1,\l},\cdots,
v_{i,\l}$ which we are going to denote by
$$V_{i,\l}:=L\{v_{1,\l},\cdots, v_{i,\l}\}$$ for $i=1, \cdots, k$.
Consider the space $$V_i=\oplus_{j=1}^d
V_{i,\l_j}=span\{V_{i,\l_{j_1}}, \cdots V_{i,\l_{j_d}} \}.$$ Note
that $$V_{i,\l_{j_1}} \cap V_{i,\l_{j_2}}=0.$$ Therefore the matrix
$A$ and the Galois group $Gal(L/K)$ send the vector space (over
$L$) $V_i$ to itself. Then we can take the Galois invariant
subspace of $V_i$, or equivalently the vectors with
coordinates in $K$. Let $$V_{i,K}=V_i \cap K^m.$$ Now we have a
filtration of vector spaces over $K$, namely, $$V_{1,K} \subset
V_{2,K}\subset \cdots \subset V_{k,K}.$$ Let
$$V_{i,{\cal{O}}_K}:=V_{i,K}\cap {\cal{O}}_K^m.$$ 
$V_{i,{\cal{O}}_K}$ is a free ${\cal{O}}_K$-module because it is a 
a finitely genrated, torsion free module over a principal ideal domain.
Choose a basis $$\{e_1,\cdots,
e_{i_1}\}$$ of $V_{1,{\cal{O}}_K}$. One can extend it to a basis
$$\{e_1,\cdots, e_{i_1},e_{i_1+1},\cdots, e_{i_2}\}$$ of $V_{2,{\cal{O}}_K}$.
In the same way one can extend the basis successively so that at
the end we have a basis $\{e_1,\cdots, e_{i_k}\}$ of $V_{k,{\cal{O}}_K}$
whose restriction to $V_{i,{\cal{O}}_K}$ gives again a basis. Finally, take
the basis $$\{e_1,\cdots, e_{i_k}\}$$ and extend it to a basis
$\{e_1,\cdots, e_m\}$ of ${\cal{O}}_K^m$. Let $P$ be a matrix whose column
vectors are $$\{e_1,\cdots, e_m\}$$ in the same order. Then
$P^{-1}AP$ has the following properties:(1) Consider $P^{-1}AP$ in
a block form with blocks $i_k\times i_k$, $(m-i_k)\times i_k$,
$i_k\times(m-i_k)$ and $(m-i_k)\times(m-i_k)$. Then the block
$(m-i_k)\times i_k$, under the diagonal is zero. (2) In case we
have $k>1$ we can describe the block $i_k\times i_k$. The block
$i_k\times i_k$ can be considered as a matrix consisting of $k^2$
square blocks of size $d\times d$. Moreover, under the diagonal
the block are zero.

 We can apply an induction argument on the
block $(m-i_k)\times(m-i_k)$ treating it as a matrix in
$GL_{m-i_k}({\cal{O}}_K)$. When this process is over, we have the matrix $A$
conjugated by a matrix $Q$ so that the resulting matrix $Q^{-1}AQ$
has the following property: each block on the diagonal has
irreducible characteristic polynomial, and each block under the
diagonal is zero.

If we relax the condition that ${\cal{O}}_K$ is a principal ideal domain
then instead of matrices $A_{ij}$ interpreted as homomorphism of free modules,
 we should have $A_{ij}$ to be a homomorphism of projective modules.
\bde{1.6}
Let $P$ be a finitely generated torsion free module over 
an integral Noetherian ring $R$ of dimension 1. 
Let $K$ be 
the field of fractions of $R$. Let $A$ be an endomorphism of $P$, 
$$A\in \End_R(P).$$ Then $A$ can be extended to an endomorphism of 
a finite dimensional vector space $P\otimes_R K$. Then the characteristic 
polynomial of $A$ is defined to be the characteristic polynomial of 
the induced map in $\End_K(P\otimes_R K)$.
\ede
The examples of such rings that we consider are orders in a number rings.
If the ring is also integrally closed, that is it is a Dedekind domain,
we have the following generalization of the previous theorem.

\bth{1.6}
Let $R$ be a Dedekind domain. And let $P$ be a finitely generated
projective module over $R$. Assume that the field of fractions $K$ has 
characteristic $0$. Let, also, $A$ be an endomorphism of $P$. Then 
in a suitable basis we have
$$
A=
\left[
\begin{tabular}{cccc}
$A_{11}$ & $\dots$  & $A_{1d}$\\
$0$      & $\ddots$ & $\vdots$\\
$0$      & $0$      & $A_{dd}$ 
\end{tabular}
\right],
$$
where $A_{ij}\in \Hom(P_j,P_i)$, with 
$P=P_1\oplus \dots \oplus P_d,$ so that $A_{ij}=0$ for $i>j$ and $A_{ii}$ has 
irreducible over $K$ characteristic polynomial.
\eth
\proof Consider $A$ as an endomorphism of $P\otimes_R \overline{K}$. Let $\l_1$
be an eigen values of $A$. Let $\l_1,\dots, \l_r$ be the Galois conjugates of 
$\l_1$. Let $v_1$ be an eigenvector in $P\otimes_R\overline{K}$ 
corresponding to the eigenvalue $\l_1$. Let $v_1,\dots ,v_r$ be
the Galois conjugates of $v_1$. Let $V=span(v_1,\dots,v_r)$. And let
$V_K=P\otimes_R K \cap V$ and $V_R=P\cap V$. We have that 
$$rank_R(V_R)=dim_K(V_K)=dim_{\overline{K}}(V).$$
Note that the endomorphism $A$ send $V_R$ to itself.  
We claim that $V_R$ and $P/V_R$ are projective $R$-modules and therefore 
direct sumand of $P$. 

Let $\pi:P\rightarrow P/V_R$. Suppose $P/V_R$ has torsion. Let $T$ be the 
torsion submodule of $P/V_R$. Then $V_R$ is of finite index in $\pi^{-1}(T)$.
Therefore $$\pi^{-1}(T)\subset (\pi^{-1}(T)\otimes_R K) \cap P=V_K\cap P=V_R.$$
Then $R/V_R$ is torsion free. We are working over a Dedekind domain. 
Then the finitely generated torsion free modules are projective. Thus,
$P/V_R$ is a direct summand of $P$. Set $P_d=V_R$. Note that $A$ 
induces the zero homomorphism in $\Hom(P/P_d,P_d)$. Proceed by induction 
util $P$ is factorized into projective modules 
$P=P_1 \oplus,\dots,\oplus P_d$ so that $A$ induces the trivial map in
$\Hom(P_j,P_i)$ for $i>j$.

\subsection{Irreducible blocks}
We have the following characterization of matrices with an
irreducible characteristic polynomial.
For principal ideal domains the proof is simpler.
\bpr{1.3}
Let ${\cal{O}}_K$ be the ring of integers in
a number field $K$.
Assume ${\cal{O}}_K$ is a principal ideal domain.
Let $f(t)$ be an irreducible over $K$ monic polynomial. 
Consider the set of all matrices $\{A_i\}$ which 
have characteristic polynomial $f(t)$ and which are not conjugate to each
other via $GL_m({\cal{O}}_K)$ where $m=deg(f)$. 
This set is parametrized by the ideal classes
in  ${\cal{O}}_K[t]/(f(t))$.
\epr
\proof The idea of the proof is to consider ${cal{O}}_K$-endomorphism of ideals
in ${\cal{O}}_K[\l]$, where $\l$ is a root of $f(t)$ 
in an algebraic closure of $K$. 
Let $I$ be an ideal in
$${\cal{O}}_K[\la] \cong {\cal{O}}_K[t]/(f(t)).$$
It is a free ${\cal{O}}_K$-module because ${\cal{O}}_K$ is a principal
ideal domain.
And let 
$$I={\cal{O}}_K\{\a_1, \cdots,\a_m\}$$ 
be a basis of $I$. 
Then the multiplication by $\la$ is an
endomorphism of the ideal $I$. Let $$\la \cdot \a_i=\sum_j
a_{ij}\a_j,$$ where $a_{ij} \in {\cal{O}}_K.$ If $A=(a_{ij})$ and
$$\overrightarrow{\a}=(\a_1, \dots, \a_m)^t$$ then
$\overrightarrow{\a}$ is an eigenvector of $A$ with eigenvalue
$\l$. Also the matrix $A=(a_{ij})$ has characteristic polynomial
$f$, because one of the eigenvalues of $A$ is $\la$ and because
$a_{ij}\in {\cal{O}}_K$.

 If we choose an ideal from the same ideal class
we have $I'=\a ' I$ for some $\a ' \in K[\la]$. Form the basis
$$I'={\cal{O}}_K\{\a'\a_1, \cdots, \a'\a_m\}$$ we obtain the same matrix
$A$.\\ Let $$I={\cal{O}}_K\{\b_1, \cdots,\b_m\}$$ be a different basis of
$I$. Let also $$\b_i=\sum_j b_{ij}\a_j,$$ and let $$B=(\b_{ij}).$$ Let
$$\b=(\b_1, \cdots , \b_m)$$ and  $$\a=(\a_1, \cdots, \a_m)$$ be vectors.
Then $$BAB^{-1}\cdot \b = BA\cdot\a=B\cdot \la \a= \la B\a=\la \b.$$
Then the matrix corresponding to the basis 
$$I={\cal{O}}_K\{\b_1,\cdots,\b_m\}$$ 
is $BAB^{-1}$, which is conjugate to $A$. In this
way we associate a matrix with integer coefficients to an element
of the ideal class. 

Conversely, if $A$ is a matrix with 
coefficients in ${\cal{O}}_K$ 
and with irreducible over $K$ characteristic polynomial $f$, we want to
associate an element of the ideal class of ${\cal{O}}_K[t]/(f(t)).,$ 
which is a free ${\cal{O}}_K$-module. If we
start with the ideal $(1)$ with ${\cal{O}}_K$-basis 
$${\cal{O}}_K \{1,\l, \cdots ,\l^{m-1}\}$$ 
then 
$$A_0=
\left[\begin{tabular}{cccc}
$0$      &  $1$     &           &\\ 
$\vdots$ &          & $\ddots$  &\\
$0$      &          &           & $1$\\
$-a_0$   & $\cdots$ &           & $-a_{m-1}$\\
\end{tabular} \right]
$$
is the corresponding matrix, where
$\l$ is a root of the characteristic polynomial
$$f(t)=t^m+a_{m-1}t^{m-1}+\cdots +a_0.$$ Using proposition 1.1, we
can find $P \e GL_m K$ such that $A=PA_0P^{-1}$. We can also
assume that $P$ has integer entries otherwise we can consider
$N\cdot P$ where $N$ is a suitable integer. Let 
$$\a_i=\sum_j p_{ij}\l^{j-1},$$ 
where $(p_{ij})=P$.\\ Claim: The ${\cal{O}}_K$-submodule
$$I={\cal{O}}_K\{\a_1, \cdots, \a_m\}$$ 
is an ideal in ${\cal{O}}_K[\l]$ and its
automorphism with respect to the given basis is the matrix $A$.\\
Before we prove that let 
$$\overrightarrow{\a}=(\a_1,\cdots, \a_m)^t$$ and 
$$\overrightarrow{\l}=(1,\l, \cdots ,\l^{m-1})^t$$ be vectors. 
Then $$\l\cdot \vect{\a}=\l
P\cdot \vect{\l}= P\cdot A_0 \vect{\l}=PA_0P^{-1}\cdot
P\vect{\l}=PA_0P^{-1}\cdot \vect{\a}=A\cdot \vect{\a}.$$ As a
consequence the ${\cal{O}}_K$-module $I$ is an ideal, since each of its
basis elements transforms via $$\l\a_i=\sum_j a_{ij}\a_j,$$ where
$a_{ij}\in {\cal{O}}_K$. Also, $A$ corresponds to the multiplication by
$\l$. Thus $A$ corresponds to the ideal class of $I$.
\bpr{1.3}
Let $f(t)$ be an irreducible monic polynomial. 
And let $P$ be a projective ${\cal{O}}_K$-module.
Consider the set of all endomorphisms $\{A_i\}$ of $P$ which 
have characteristic polynomial $f(t)$ and which are not conjugate to each
other via an automorphism of $P$. This set is parametrized by the ideal 
classes
in  ${\cal{O}}_K[t]/(f(t))$ which are isomorphic to $P$ as 
${\cal{O}}_K$-modules.
\epr
\proof  The idea of the proof is to consider ${\cal{O}}_K$-endomorphism 
of ideals
in ${\cal{O}}_K[\l]$, where $\l$ is a root of $f(t)$ 
in an algebraic closure of $K$. 
Let $I$ be an ideal in
$${\cal{O}}_K[\la] \cong {\cal{O}}_K[t]/(f(t)),$$
which is isomorphic to $P$ as ${\cal{O}}_K$-module.
Let $A$ be an endomorphism of $I$ induces by multiplication by $\l$.
Tensoring the endomorphism with the field $K$ we obtain that $\l$ is
 an eigenvalue. Thus, $f(t)$ is the characteristic polynomial of $A$.
Also, if $B$ is an ${\cal{O}}_K$-automorphism of $P$, then is it also an 
${\cal{O}}_K$-automorphism of $I$. Then $A$ sends $BI$ to $BI$. Therefore, 
the induced endomorphism on $I$ is $BAB^{-1}$. Thus $A$ is determined up 
to conjugation. If $J$ is in the same ideal class as $I$, then $J=\a I$, 
for $\a\in K(\la)^{\times}$. That induces an isomorphism between the 
${\cal{O}}_K$-modules $I$ and $J$. Thus, an endomorphism on one lead to 
an endomorphism of the other.
This is the one of the direction of the correspondance between 
ideal classes and endomorphisms.

Given an endomorphism of $P$ 
with characteristic polynomial $f(t)$,
we want ot associate an ideal class in ${\cal{O}}_K[\l]$.
Note that the ring ${\cal{O}}_K[A]$ is isomorphic to ${\cal{O}}_K[\l]$,
since $f(A)=0$ and $f(\l)=0$. The ring ${\cal{O}}_K[A]$ acts on $P$, because $P$ is an ${\cal{O}}_K$-module, and also $A$ acts on $P$. Therefore $P$
is an ${\cal{O}}_K[\l]$-module because ${\cal{O}}_K[\l]$ is isomprphic to
${\cal{O}}_K[A]$. It remains to show that torsion free 
${\cal{O}}_K[\l]$-modules of rank $1$ correspond to ideal classes. 
Consider the imbedding $P\rightarrow P\otimes_{{\cal{O}}_K[\l]} K(\l)$
followed by the isomorphism  
$\theta:P\otimes_{{\cal{O}}_K[\l]} K(\l) \rightarrow K(\l)$.
The image og the composition is a fractional ideal in $K(\l)$. All the
choises for the isomorphism are parametrized by $K(\l)^{\times}$.
 Therefore, the torsion free ${\cal{O}}_K[\l]$-modules of rank $1$ 
correspond to the set of ideal classes module $K(\l)^{\times}$ 
which is precisely the ideal classes in ${\cal{O}}_K[\l]$.

   \bco {1.4}
Let $f(t)\in {\cal{O}}_K[t]$ 
be an irreducible monic polynomial. Assume that ${\cal{O}}_K[t]/(f(t))$ 
is a integrally closed. Then the
conjugacy classes of matrices in $GL_m({\cal{O}}_K)$ with characteristic
polynomial $f(t)$ are in one-to-one correspondence with the
elements in 
$$\ker(K_0({\cal{O}}_K[t]/(f(t)))\rightarrow K_0({\cal{O}}_K)).$$
   \eco
\proof
It is known that for number rings ${\cal{O}}_K$ we have
$$K_0({\cal{O}}_K)=\Z \oplus Cl(K),$$ 
where $Cl(K)$ is the ideal class group of 
${\cal{O}}_K$ (see \cite{M}). 
Then the ideal classes in ${\cal{O}}_K[t]/(f(t))$ which are free 
${\cal{O}}_K$-modules are precisely 
$$\ker(K_0({\cal{O}}_K[t]/(f(t)))\rightarrow K_0({\cal{O}}_K)).$$ 
\subsection{Reducible blocks}
Using the theorem 1.5 we obtain the following: If $A$ is an $m\times m$ 
 matrix with coefficients in
${\cal{O}}_K$,
or more generally, an endomorphism of projective ${\cal{O}}_K$-module of
rank $m$, having an irreducible characteristic polynomial then $A$ is
conjugated by an element of $GL_m({\cal{O}}_K)$, or
of $Aut(P)$,
 to a $2\times 2$-block 
endomorphism 
$$
\left[
\begin{tabular}{cc}
  $A_{11}$ & $A_{12}$\\
  0      & $A_{22}$\\
\end{tabular}
\right].$$\\
Let $A_{11}$ and 
$A_{22}$ be automorphisms of projective modules $P_1$ and $P_2$, with
$$P_1\oplus P_2=P.$$ 
Now we are going to describe a method that simplifies the block
$A_{12}$ and leaves $A_{11}$ and $A_{22}$ unchanged. We can assume
that $$A=\left[
\begin{tabular}{cc}
  $A_{11}$ & $A_{12}$\\
  0      & $A_{22}$\\
\end{tabular}
\right].$$ \
Conjugate $A$ with an automorphism $$B=\left[
\begin{tabular}{cc}
  $B_{11}$ & $B_{12}$\\
  0        & $B_{22}$\\
\end{tabular}
\right]$$ of the some block form.
 that is $B_{11}\in \Aut(P_1)$,
$B_{22}\in \Aut(P_2)$ and 
$B_{12}\in \Hom(P_2,P_1)$. We want the
conjugation by $P$ to preserve $A_{11}$ and $A_{22}$ that is
     $$
\left[\begin{tabular}{cc} $B_{11}$ & $B_{12}$\\ 0 & $B_{22}$\\
\end{tabular} \right]
\cdot \left[\begin{tabular}{cc} $A_{11}$ & $A_{12}$\\ 0 &
$A_{22}$\\ \end{tabular} \right]
=
\left[\begin{tabular}{cc} $A_{11}$ & $A'_{12}$\\ 0 & $A_{22}$\\
\end{tabular} \right]
\cdot \left[\begin{tabular}{cc} $B_{11}$ & $B_{12}$\\ 0 &
$B_{22}$\\ \end{tabular} \right]. $$
     We do have that $A_{12}$ is
changed to $A'_{12}$. Then $$B_{11}A_{11}=A_{11}B_{11},$$
$$B_{22}A_{22}=A_{22}B_{22},$$ and
$$A'_{12}B_{22}-B_{11}A_{12}=B_{12}A_{22}-A_{11}B_{12}.$$ Denote by
$C(A_{ii})$ the centralizer of $A_{ii}$. Then $B_{11}\in C(A_{11})$
and $B_{22}\in C(A_{22})$. The block $B_{12}$ could be any map
in $\Hom_{{\cal{O}}_K}(P_2,P_1)$. Let 
$$P_{A_{11},A_{22}}: B_{12} \mapsto B_{12}A_{22}-A_{11}B_{12}$$ 
be a map from the space
$\Hom_{{\cal{O}}_K}(P_2,P_1)$ to itself. 
It is a linear. Then the relation
$$A'_{12}B_{22}-B_{11}A_{12}=B_{12}A_{22}-A_{11}B_{12},$$ between
$A'_{12}$ and $A_{12}$ can be written as $$A'_{12}B_{22} \equiv
B_{11}A_{12}\: \mod (\Im P_{A_{11},A_{22}}).$$


\ble{1.6}
   Let $P_{mod}=\Im (P_{A_{11},A_{22}})$ where
$P_{A_{11},A_{22}}: B_{12} \mapsto B_{12}A_{22}-A_{11}B_{12}$ Let
also $Q_{mod}=\Hom_{{\cal{O}}_K}(P_2,P_1) /P_{mod}$. 
Then $\Hom_{{\cal{O}}_K}(P_2,P_1)$,
$P_{mod}$ and $Q_{mod}$ are
 $C(A_{11})\times C(A_{22})$-modules.
   \ele
     \proof Let $B_{11}\in C(A_{11})$ and $B_{22}\in C(A_{22})$.
Then
$$
\begin{tabular}{ccl}
$(B_{11}, B_{22})\cdot P_{A_{11},A_{22}}(B_{12})$ 
& $=$ & $B_{11}\cdot P_{A_{11},A_{22}}(B_{12})\cdot B^{-1}_{22}=$\\
& $=$ & $B_{11}(B_{12}A_{22}-A_{11}B_{12})B^{-1}_{22}=$\\
& $=$ & $B_{11}B_{12}A_{22}B^{-1}_{22}-B_{11}A_{11}B_{12}B^{-1}_{22}=$\\
& $=$ & $(B_{11}B_{12}B^{-1}_{22})A_{22}-A_{11}(B_{11}B_{12}B^{-1}_{22})=$\\
& $=$ & $P_{A_{11},A_{22}}(B_{11}B_{12}B^{-1}_{22}).$\\
\end{tabular}
$$
Obviously,
$\Hom_{{\cal{O}}_K}(P_2,P_1)$ is a 
$C(A_{11})\times C(A_{22})$-module. Thus,
the quotient, as abelian group, 
$Q_{mod}=\Hom_{{\cal{O}}_K}(P_2,P_1) /P_{mod}$
has the structure of a $C(A_{11})\times C(A_{22})$-module.


\bpr{1.7}
  Let $A_{11}\in Aut(P_1)$ and $A_{22}\in Aut(P_2)$.
Suppose $A_{11}$ and $A_{22}$ have no common eigenvalues. Then the
matrices
$$\left[\begin{tabular}{cc} $A_{11}$ & $A_{12}$\\ $0$ & $A_{22}$\\
\end{tabular} \right] \mbox{ and }
\left[\begin{tabular}{cc} $A_{11}$ & $A'_{12}$\\ $0$ & $A_{22}$\\
\end{tabular} \right]$$
are conjugate to each other if and only if the projection of
$A_{12}$ and $A'_{12}$ onto the finite set $C(A_{11})\times
C(A_{22})\backslash Q_{mod}$ coincide.
   \epr
In order to prove the above proposition, we need the following
definition and lemma. Consider the linear space of $n\times
m$-matrices, or more generally of homomorphisms.
 Multiplication on the left and multiplication on the
right are linear transformations of this space. So this
transformations can be written as matrices, or as endomorphisms. 
Now we are going to
set the notation for these matrices, after tensoring with $K$.
 Let $A=(A_{ij})$ and
$B=(b_{kl})$ be two matrices. We shall write the tensor $A\otimes
B$ as a matrix in the following way: The matrix can be considered
as a block-matrix with each block being the size of $A$, and he
$(k,l)$-block being $b_{kl}A$. This tensor product has the
following properties: Let $A$ be an $m\times m$ matrix acting by
left multiplication on $M\e Mat_{m,n}$, i.e on the $m\times n$
matrices. Consider $M$ as a vector by arranging the column vectors
of $M$ one below the other. Then the left multiplication can be
expressed as a matrix. And this matrix is precisely $A\otimes
I_n$. If we multiply on the right by an $n\times n$ matrix $B$
then such a linear transform can be expressed by $I_m\otimes B^t$.
Another property that is easy to check is $(A\otimes B)(C\otimes
D)=(AC)\otimes(BD)$.


\ble{1.8}
    If $\prod_i (t-\a_i)$ and $\prod_j (t-\b_j)$ are
the characteristic polynomials of $A_{11}\in GL_{m_1}\bar{K}$ and
$A_{22}\in GL_{m_2}\bar{K}$ then the characteristic polynomial of
$$P_{A_{11},A_{22}}:X\mapsto XA_{22}-A_{11}X$$ 
is $\prod_{i,j}
(t-\b_j+\a_i)$. In particular, If $A_{11}$ and $A_{22}$ have no 
common eigenvalue then the map $P_{A_{11},A_{22}}$ is 
non-singular, and if $R$ is a number ring $Q_{mod}$ is finite, 
namely the norm of the resultant. 
\ele
  \proof
We write the matrix $X$ as a vector, arranging the columns of $X$ 
below each other. Then the linear transform $P_{A_{11},A_{22}}$ 
becomes $$I_{m_1}\otimes A_{22}^t -A_{11}\otimes I_{m_2}.$$ If 
$A_{11}$ and $A^t_{22}$ are diagonal matrices with diagonals
$$(\a_1, \cdots ,\a_{m_1})$$ and $$(\b_1,\cdots, \b_{m_2})$$ then
$P_{A_{11},A_{22}}$ would be diagonal with diagonal entries
$\b_j-\a_i$. If $A_{11}$ and $A_{22}$ are arbitrary then we can
find $B_1 \in GL_{m_1} \bar{K}$ and $B_2\in GL_{m_2}\bar{K}$ 
with such that the
matrices $B_1A_{11}B_1^{-1}$ and $B_2A_{22}^t B_2^{-1}$ are in
Jordan block form. Let $$(\a_1, \cdots ,\a_{m_1})$$ be the diagonal
entries of $B_1A_{11}B_1^{-1}$ and $$(\b_1,\cdots, \b_{m_2})$$ be the
diagonal entries of $B_2A_{22}^t B_2^{-1}$. Note that 
$$
\begin{tabular}{ccl}
$B_1\otimes B_2(P_{A_{11},A_{22}})B_1^{-1}\otimes B_2^{-1}$ 
 & $=$ & 
$B_1\otimes
B_2(I\otimes A_{22}^t -A_{11}\otimes I)B_1^{-1}\otimes B_2^{-1}$\\
& $=$ & $I\otimes B_2A_{22}^tB_2^{-1} -B_1A_{11}B_1^{-1}\otimes I$\\
\end{tabular}
$$ 
is upper triangular. Its diagonal entries are $\b_j-\a_i$. Therefore
if $\a_i \neq \b_j$ for all $i$ and $j$ then  $P_{A_{11},A_{22}}$
is non-singular. If $R$ is a number ring 
then $\Im(P_{A_{11}A_{22}})=P_{mod}$ is a
sublattice of maximal rank in $\Hom_{{\cal{O}}_K}(P_2,P_1)$ of index
the norm $N_{K/\Q}$ of the resultant. Therefore,
$Q_{mod}=\Hom(P_2,P_1)/P_{mod}$ is finite.\\
\\
\proof (of proposition 1.7) If $A_{12}$, 
$A_{12} \in \Hom_{{\cal{O}}_K}(P_2,P_1)$
are such that the the endomorphisms

       $$
 \left[\begin{tabular}{cc} $A_{11}$ & $A_{12}$\\ $0$ & $A_{22}$\\
\end{tabular} \right]
\mbox{ and } \left[\begin{tabular}{cc} $A_{11}$ & $A'_{12}$\\ $0$
& $A_{22}$\\ \end{tabular} \right]
    $$
are conjugate to each other in $\Aut_{{\cal{O}}_K}(P)$,
then there exists $B$ in
$\Aut_{{\cal{O}}_K}(P)$ of the same block form such that
$$
\left[\begin{tabular}{cc} $B_{11}$ & $B_{12}$\\ $B_{21}$ &
$B_{22}$\\ \end{tabular} \right] \cdot \left[\begin{tabular}{cc}
$A_{11}$ & $A_{12}$\\ $0$ & $A_{22}$\\ \end{tabular} \right]
=
\left[\begin{tabular}{cc} $A_{11}$ & $A'_{12}$\\ $0$ & $A_{22}$\\
\end{tabular} \right]
\cdot \left[\begin{tabular}{cc} $B_{11}$ & $B_{12}$\\ $B_{21}$ &
$B_{22}$\\ \end{tabular} \right].
$$
Then the $(2,1)$-entry of the product gives 
$$B_{21}A_{11}=A_{22}B_{21}.$$ However, the map
$$P_{A_{22}A_{11}}:B_{21}\mapsto B_{21}A_{11}-A_{22}B_{21}=0$$ is
non-singular from the previous lemma. Therefore we have
$B_{21}=0$. Now the conditions on the other blocks of $B$ become
simpler: $$B_{11}A_{11}=A_{11}B_{11},$$ $$B_{22}A_{22}=A_{22}B_{22}$$
and $$A'_{12}B_{22}-B_{11}A_{12}= B_{12}A_{22}-A_{11}B_{12}\e \Im
P_{A_{11},A_{22}}.$$ Then $$A'_{12}\equiv
B_{11}A_{12}B_{22}^{-1}\:\mod \mbox{ } \Im \:(P_{A_{11},A_{22}}).$$
Thus, the image of $A_{12}$ and $A'_{12}$ in $$C(A_{11})\times
C(A_{22})\backslash Q_{mod}$$ coincide.\\
 Conversely, if $A_{12}$
and $A'_{12}$ map to the same element in $$C(A_{11})\times
C(A_{22})\backslash Q_{mod}$$ then $$A'_{12}\equiv
B_{11}A_{12}B_{22}^{-1} \mod \mbox{ } \Im \:(P_{A_{11}A_{22}}),$$
with $B_{11}\e C(A_{11})$ and $B_{22}\e C(A_{22})$. Equivalently,
$$A'_{12}B_{22}\equiv B_{11}A_{12} \mod \mbox{ } \Im
\:(P_{A_{11}A_{22}}).$$
    Then there exists $B_{12}$ such that
$$P_{A_{11}A_{22}}(B_{12})= A'_{12}B_{22}-B_{11}A_{12}.$$
Equivalently,
$$B_{12}A_{22}-A_{11}B_{12}=A'_{12}B_{22}-B_{11}A_{12}.$$ Therefore,
   $$
    \left[\begin{tabular}{cc} $B_{11}$ & $B_{12}$\\ $0$ &
$B_{22}$\\ \end{tabular} \right] \cdot \left[\begin{tabular}{cc}
$A_{11}$ & $A_{12}$\\ $0$ & $A_{22}$\\ \end{tabular} \right]
=
\left[\begin{tabular}{cc} $A_{11}$ & $A'_{12}$\\ $0$ & $A_{22}$\\
\end{tabular} \right]
\cdot \left[\begin{tabular}{cc} $B_{11}$ & $B_{12}$\\ $0$ &
$B_{22}$\\ \end{tabular} \right].
    $$


\bco{1.9}
      If $f(t)$ is monic polynomial of degree $m$ with integer
coefficient, and $f(0)=\pm 1$ without repeated roots then there
are finitely many matrices $A\in GL_m ({\cal{O}}_K)$ with characteristic
polynomial $f(t)$, where $K$ is a number field.\\
\eco
\proof
We will proceed by induction on the number if irreducible factor of 
the characteristic polynomial. If the characteristic polynomial is 
irreducible we use the finiteness of class numbers and use corollary 1.4.
Suppose we have proven the statement when the characteristic polynomial 
factors into $n$ irreducible polynomials. Consider a matrix $A$ whose 
characteristic polynomial factors into $n+1$ irreducible polynomials.
We can assume that $A$ is of the form
$$A= \left[\begin{tabular}{cc} $A_{11}$ & $A_{12}$\\ $0$ & $A_{22}$\\
\end{tabular} \right],$$
where $A_{11}$ has irreducible characteristic polynomial and 
the characteristic polynomial of $A_{22}$ factors into $n$
irreducible polynomials. There are finitely many matrices having 
the same characteristic polynomial as $A_{11}$ and being non-conjugate to 
$A_{11}$.
Also, by assumption we can replace $A_{22}$
only by finitely many matrices in order to obtain  new matrices 
non-conjugate to the initial one and with the same characteristic 
polynomial. To prove finiteness of matrices that are non-conjugate to $A$
but have the same characteristic polynomial,
 we use that for fixed
$A_{11}$ and $A_{22}$ there are only finitely many options for $A_{12}$
which follow from proposition 1.7.


\sectionnew{Centralizers and homological Euler characteristics}

\ble{2.1}
   Let $A_{11}$ and $A_{22}$ be square matrices with coefficients in 
${\cal{O}}_K$. Suppose that they have no common eigenvalues.
And let
 $$C=\left[
\begin{tabular}{cc} $C_{11}$ & $C_{12}$\\ $C_{21}$ & $C_{22}$\\
\end{tabular}
\right] \mbox{ commutes with }
 A=\left[
\begin{tabular}{cc} $A_{11}$ & $A_{12}$\\ $0$ & $A_{22}$\\
\end{tabular}
\right]$$
    Then the admissible matrices $C$ are determined by
the following properties:
 $$C_{21}=0,$$ $$C_{11}\in C(A_{11}),$$ $$C_{22}\in C(A_{22}),$$
$$C_{11}A_{12}C_{22}^{-1}\equiv A_{12}\mod \mbox{ }\Im
P_{A_{11}A_{22}}.$$
   Also, the matrix $C_{12}$ is uniquely
determined by $C_{11}$ and $C_{22}$, and it is given by
 $$C_{12}=P_{A_{11}A_{22}}^{-1}(A_{12}C_{22}-C_{11}A_{12}).$$
In particular,
$$C(\left[
\begin{tabular}{cc} $A_{11}$ & $0$\\ $0$ & $A_{22}$\\
\end{tabular}
\right])=C(A_{11})\times C(A_{22}),$$
         where $C(A)$ denotes the
centralizer of $A$.\\
\\
\ele
  \proof We have
$$
   \left[\begin{tabular}{cc}
$C_{11}$ & $C_{12}$\\ $C_{21}$ & $C_{22}$\\
\end{tabular} \right]
\cdot \left[\begin{tabular}{cc} $A_{11}$ & $A_{12}$\\ $0$ &
$A_{22}$\\
\end{tabular} \right]
=
\left[\begin{tabular}{cc} $A_{11}$ & $A_{12}$\\ $0$ & $A_{22}$\\
\end{tabular} \right]
\cdot \left[\begin{tabular}{cc} $C_{11}$ & $C_{12}$\\ $C_{21}$ &
$C_{22}$\\
\end{tabular} \right].
  $$
The $(2,1)$-entry of the products give
$$C_{21}A_{11}=A_{22}C_{21}$$ However the map
$$P_{A_{22}A_{11}}:C_{21}\mapsto C_{21}A_{11}-A_{22}C_{21}$$ is
non-singular by lemma 1.11. Therefore $C_{21}=0$. Considering the
$(1,1)$-entry and the $(2,2)$-entry of the above matrix product we
obtain that $C_{11}\in C(A_{11})$ and $C_{22}\in C_(A_{22})$.\\
  Fix
$C_{11}$ and $C_{22}$ to be in the centralizers $C(A_{11})$ and
$C(A_{22})$, respectively. Then use that $A=CAC^{-1}$ in order to
determine the block $C_{12}$.
  $C_{11}\in C(A_{11})$ and $C_{22}\in
C(A_{22})$. Then
$$
     \begin{array}{rl} A=CAC^{-1}
= & \left[\begin{tabular}{cc} $C_{11}$ & $C_{12}$\\ $0$      &
$C_{22}$\\
\end{tabular} \right]
\cdot \left[\begin{tabular}{cc} $A_{11}$ & $A_{12}$\\ $0$      &
$A_{22}$\\
\end{tabular} \right]
\cdot \left[\begin{tabular}{cc} $C_{11}^{-1}$ &
$-C_{11}^{-1}C_{12}C_{22}^{-1}$\\ $0$            & $C_{22}^{-1}$\\
\end{tabular} \right]=\\
= & \left[\begin{tabular}{cc} $C_{11}A_{11}$ &
$C_{11}A_{12}+C_{12}A_{22}$\\ $0$            & $C_{22}A_{22}$\\
\end{tabular} \right]
\cdot \left[\begin{tabular}{cc} $C^{-1}_{11}$ &
$-C_{11}^{-1}C_{12}C_{22}^{-1}$\\ $0$      & $C^{-1}_{22}$\\
\end{tabular} \right]=\\
= & \left[\begin{tabular}{cc} $A_{11}$ & $-A_{11}C_{12}C_{22}^{-1}
            +C_{11}A_{12}C_{22}^{-1}+C_{12}A_{22}C_{22}^{-1}$\\
$0$      & $A_{22}$\\
\end{tabular} \right].
\end{array}
    $$
We have
$$
\begin{tabular}{ccl}
$A_{12}$ & $=$ & 
$C_{11}A_{12}C_{22}^{-1}+(C_{12}A_{22}-A_{11}C_{12})C_{22}^{-1}=$\\ 
& $=$ & $C_{11}A_{12}C_{22}^{-1} +P_{A_{11}A_{22}}(C_{12})C_{22}^{-1}.$ 
\end{tabular}
$$
Therefore $A_{12}$ and
$C_{11}A_{12}C_{22}^{-1}$ coincide as elements in
$$Q_{mod}=\Hom_{{\cal{O}}_K}(P_2,P_1)/\Im P_{A_{11}A_{22}}.$$ Therefore
$C_{11}\times C_{22}$ stabilizes $A_{12}$ as element of $Q_{mod}$.
Conversely, if $C_{11}\times C_{22}$ stabilizes $A_{12}$ as an
element of $Q_{mod}$ we can find a unique $C_{12}$ such that

$$
   C=\left[\begin{tabular}{cc}
$C_{11}$ & $C_{12}$\\ $0$      & $C_{22}$\\
\end{tabular}
     \right]
\mbox{ and } A=
   \left[\begin{tabular}{cc}
$A_{11}$ & $A_{12}$\\ $0$      & $A_{22}$\\
\end{tabular} \right]
    $$
   commute. Indeed, if $$A_{12}-C_{11}A_{12}C_{22}^{-1}
\mbox{ }\e \mbox{ }\Im (P_{A_{11}A_{22}})$$ then
$$A_{12}C_{22}-C_{11}A_{12} \mbox{ } \e \mbox{ }\Im
(P_{A_{11}A_{22}}).$$ And we can set
$$C_{12}=P_{A_{11}A_{22}}^{-1}(A_{12}C_{22}-C_{11}A_{12}).$$ Then
the equation $$A_{12}=C_{11}A_{12}C_{22}^{-1}
+P_{A_{11}A_{22}}(C_{12})C_{22}^{-1}$$ is satisfied.
\\

Let
  $$
  A=\left[\begin{tabular}{cc} $A_{11}$ & $A_{12}$\\ $0$      &
$A_{22}$\\
\end{tabular} \right].
    $$
We are going to define a map $i_{A_{12}}:C(A)\rightarrow
C(A_{11})\times C(A_{22})$ And let $C$ be in the centralizer of
$A$. Consider $C$ as a block matrix of the same block form as $A$.
Then by lemma 2.1 we have that $C_{21}=0$ and that $C_{11}\in
C(A_{11})$ and $C_{22}\in C(A_{22})$. Let $i_{A_{12}}$ sends $C$ to
$C_{11}\times C_{22}$. Then $i_{A_{12}}$ is a homomorphism.
\bco{2.2}
   With the above definition of $i_{A_{12}}$ we have that
$i_{A_{12}}$ is always injective and the index of
$\Im(i_{A_{12}})$ in $C_{11}\times C_{22}$ is equal to 
the number of elements in
 $C_{11}\times C_{22}$-orbit
of $A_{12}$ inside $Q_{mod}$.
    \eco
    \proof Let $C$ be in the centralizer of $A$.
Let also
     $$
     C=\left[\begin{tabular}{cc} $C_{11}$ & $C_{12}$\\ $0$ &
$C_{22}$\\
\end{tabular} \right]
     $$
be of the same block form as $A$. From lemma 2.1 we know that the
blocks of $C$ must satisfy $C_{11}\in C(A_{11})$ and $C_{22}\in
C(A_{22})$ with the relation that $A_{12}$ and
$$C_{11}A_{12}C_{22}^{-1}$$ are congruent modulo the image of
$P_{A_{11}A_{22}}$. Equivalently, $A_{12}$ and
$C_{11}A_{12}C_{22}^{-1}$ coincide in $$Q_{mod}=Mat_{m_1m_2}{\cal{O}}_K /\Im
P_{A_{11}A_{22}}.$$ Therefore the image $\Im(i_{A_{12}})$ consists
of such $C_{11}$ and $C_{22}$ that the group element $C_{11}\times
C_{22}$ fixes $A_{12}$ as an element of $Q_{mod}$. Therefore the
image $\Im(i_{A_{12}})$ coincides with the stabilizer of $A_{12}\in
Q_{mod}$, and the index of $\Im(i_{}A_{12})$ in $C(A_{11})\times
C(A_{22})$ corresponds to the $C(A_{11})\times C(A_{22})$-orbit of
$A_{12}$ in $Q_{mod}$.
\ble{2.3}
    Let $A_{11}$ and $A_{22}$ be invertible matrices
with coefficients in ${\cal{O}}_K$, or automorphisms of projective modules
having no common eigenvalues. Let $f_1$
and $f_2$ be the characteristic polynomials if $A_{11}$ and
$A_{22}$, respectively. Then
$$
   \sum \chi(C(A))=|N_{K/\Q}(R(f_1,f_2))|\chi(C(A_{11}))\chi(C(A_{22})),
   $$
where the sum is taken over all non-conjugate
torsion elements

$$
   A=\left[\begin{tabular}{cc}
$A_{11}$ & $A_{12}$\\ $0$      & $A_{22}$\\
\end{tabular} \right]
   $$
   with fixed $A_{11}$ and $A_{22}$, and $R(f_1,f_2)$ is the
resultant of the two polynomials.
   \ele
\proof By Proposition 1.10, taking the sum over all non-conjugate
matrices with fixed $A_{11}$ and $A_{22}$ is the same as varying
$A_{12}$ through representatives of $$C(A_{11})\times 
C(A_{22})\backslash Q_{mod}.$$ For a fixed $A_{12}$ by corollary 2.2,
we have that the group $C(A)$ is a finite index subgroup of
$C(A_{11})\times C(A_{22})$, and the index is equal to the number
of elements in the $C(A_{11})\times C(A_{22})$-orbit of $A_{12}$
in $Q_{mod}$. Thus, $$\chi(C(A))= \#|\mbox{orbit of }A_{12}|\cdot
\chi(C(A_{11}))\chi(C(A_{22}))$$ Summing over all orbits, we obtain
 $$\sum \chi(C(A))=
\#|Q_{mod}|\cdot \chi(C(A_{11}))\chi(C(A_{22})).$$ On the other
hand, by lemma 1.11
$$\#|Q_{mod}|=|N_{K/\Q}(det(P_{A_{11}A_{22}}))|=|N_{K/\Q}(R(f_1,f_2))|.$$

We are going to define a resultant of $k$ polynomials $f_1,\dots
,f_k$, $k \geq 2$ by $$R(f_1,\dots ,f_k)=\prod_{i<j} R(f_i,f_j).$$
\bpr{2.4}
  Let $A_{11}$, $A_{22}$ , $\cdots$, $A_{kk}$ be
invertible matrices with coefficients in ${\cal{O}}_K$ such that $A_{ii}$
and $A_{jj}$ have no common eigenvalues for $i\neq j$. Let $f_i$
be the characteristic polynomial of $A_{ii}$. Then $$\sum
\chi(C(A))= |N_{K/\Q}(R(f_1,\dots, f_k))|\cdot\chi(C(A_{11}))\dots
\chi(C(A_{kk})),$$ where the sum is taken over all non-conjugate
torsion elements $A$ is a block-diagonal form such that the blocks
on the diagonal are $A_{11}$, $\dots$, $A_{kk}$, and the blocks
under the diagonal are zero, and $R(f_1,\dots ,f_k)$ is the
resultant of $f_1,\dots, f_k$ defined above.
\epr
  Note that all $A$'s in the above sum have the same characteristic
polynomial; namely $f_1\cdot \dots \cdot f_k$.\\

\proof We are going to prove the statement by induction on $k$.
For $k=2$ it is true by the previous lemma. Assume it is true for
$k-1$. Notice that $$R(f_1,\dots ,f_k)= R(f_1,\dots ,f_{k-1})\cdot
R((f_1\cdot \dots \cdot f_{k-1}),f_k).$$ Consider the matrices $A$
in the sum as $2\times 2$-block-matrices in the following way:

$$
   A=\left[\begin{tabular}{cc} $\overline{A}_{11}$ &
$\overline{A}_{12}$\\ $0$                 & $A_{kk}$\\
\end{tabular} \right],
   $$
   with $\overline{A}_{11}$ a block-triangular matrix
with blocks on the diagonal $A_{11}, \dots A_{k-1,k-1}$, and zero
blocks under the diagonal. Let $A'$ be another matrix in the above
sum. Consider it similarly as a $2\times 2$-block matrix

$$
    A'=\left[\begin{tabular}{cc}
$\overline{A'}_{11}$ & $\overline{A'}_{12}$\\ $0$
& $A_{kk}$\\
\end{tabular} \right].
   $$
  We can assume that either $\overline{A}_{11}$ and
$\overline{A'}_{11}$ coincide, or that they are non-conjugate. We
can make that assumption for the following reason. If they are
conjugate say by a matrix $\overline{B}_{11}$ instead of the
representative $A'$
 one can take $BA'B^{-1}$ where

$$
     B=\left[\begin{tabular}{cc}
   $\overline{B}_{11}$ & $0$\\
$0$                 & $I$\\
\end{tabular} \right],
  $$
with the matrix $I$ having the size of $A_{kk}$. Then the new
matrix $A''=BAB^{-1}$ will have
$\overline{A''}_{11}=\overline{A}_{11}$. Thus, we can separate the
summation over non-conjugate matrices $A$ into two summations.
First fix $\overline{A}_{11}$ and let $\overline{A}_{12}$ very
over all representatives $$C(\overline{A}_{11})\times
C(A_{kk})\backslash Q_{mod},$$ where $Q_{mod}$ is the cokernel of
the inclusion $P_{\overline{A}_{11}A_{kk}}$. Then sum over
non-conjugate $\overline{A}_{11}$ in block-diagonal form with
$A_{11},\dots,A_{k-1,k-1}$ on the diagonal, and zero under the
diagonal. The summation over all
 The first summation leads to

$$
  \sum_{\overline{A}_{11}}
|N_{K/\Q}(R((f_1\cdot \dots \cdot f_{k-1}),f_k))|\cdot
\chi(C(\overline{A}_{11}))\chi(A_{kk}),
   $$
    by the previous lemma
because $f_1\cdot \dots \cdot f_{k-1}$ is the characteristic
polynomial of any of the $\overline{A}_{11}$ matrices. And by
induction assumption
    $$
    \sum_{\overline{A}_{11}}
\chi(C(\overline{A}_{11}))= |N_{K/\Q}(R(f_1,\dots, f_{k-1}))|\cdot
\chi(C(A_{11}))\dots \chi(C(A_{k-1,k-1})).
    $$

We need e few more lemmas on the size of the centralizer.
\ble{2.5}
      Let $A$ and $B$ be matrices in $GL_m {\cal{O}}_K$, or $\Aut(P)$
which are conjugate as elements of $GL_m K$ the $C(A)$
and $C(B)$ are comensurable.
 \ele
\proof The group $C(A)$ is an arithmetic subgroup of
$C_{GL_m K}(A)$, which is conjugate to $C_{GL_m K}(B)$.
Therefore $C_{GL_m {\cal{O}}_K}(A)$ and $C_{GL_m {\cal{O}}_K}(B)$
are comeasurable.

\ble{2.6} 
   Let $T_n \in GL_m ({\cal{O}}_K)$ be an  $n$-torsion matrix with
with irreducible characteristic polynomial, or an $n$-torsion 
automorphisms. Then the
centralizer $C(T_n)$  contains ${\cal{O}}_{K}[\xi_n]^{\times}$, 
and is comensurable to it
where $\xi_n=e^{2\pi i/n}$.\\
   \ele
\proof The matrices in $Mat_{m,m} K$ commuting with $T_n$
are precisely $$K[T_n]\cong K(\xi_n)$$ because this is the maximal
abelian sub-Lie algebra commuting with $T_n$; it is of dimension
$m$. Let the intersection of $K[T_n]$ with $Mat_{m,m} {\cal{O}}_K$ be 
$$R\subset {\cal{O}}_{K(\xi_n)},$$
where $R$ and ${\cal{O}}_{K(\xi_n)}$ are comeasurable by the previous lemma.
In this intersection $R$ the invertible
elements are 
$$C(T_n)\cong R^{\times}\subset {\cal{O}}_{K[\xi_n]}^{\times}.$$
\ble{2.6}
     Let $T_n \in GL_m\Z$ be an $n$-torsion matrix with
irreducible characteristic polynomial. Let $T$ be a
$k\times k$-block matrix with blocks on the diagonal $T_n$ and the
rest of the blocks being zero. Let $R$ be the ring of endomorphisms 
that commute with $T_n$.
Then 
$C(T) \cong GL_k R.$
         \ele
  \proof Then the matrices commuting with
$T$ are $$\Q[T_n]\otimes
Mat_{k,k} K \cong Mat_{k,k}(K(\xi_n)).$$ Among them the ones with
integer coefficients are 
$$C_{Mat_{km}\Z}(T)=C_{Mat_{km}\Z}(T_n\otimes I_k)=C_{Mat_{k}\Z}(T_n)\otimes Mat_{k,k}\Z\cong Mat_{k,k}R,$$ 
where $R$ is a order in ${\cal{O}}_{K(\xi_n)}$ isomorphic to
the ring of matrices (not necessarily with unit determinant)
with coefficients in ${\cal{O}}_K$
commuting with $T_n$
And the invertible ones with integer
coefficients are $$C(T)\cong GL_k R \subset GL_k({\cal{O}}_{K(\xi_n)}).$$

The following two propositions give bases for the proof of
the vanishing results, namely theorem 0.1 and theorem 0.2.
\bpr{2.8}
   Let $A$ be a torsion element of $GL_m {\cal{O}}_K$.
Then $\chi(C(A))\neq 0$ if and only if the set of eigenvalues of
$A$ is inside the set\\

(a) $\{1,-1,i,-i,\xi_3,\bar\xi_3,\xi_6,\bar\xi_6\}$, and the
multiplicity of $1$ and $-1$ is at most $2$ and the multiplicity
of the rest of the roots of unity is at most $1$ if $K=\Q$;\\

(b) $\{1,-1,i,-i\}$, and the
multiplicities are at most $1$ if $K=\Q(i)$;\\

(c)  $\{1,-1,\xi_3,\bar\xi_3,\xi_6,\bar\xi_6\}$  and the
multiplicities are at most $1$ if $K=\Q(\xi_3)$;\\

(d) $\{1,-1\}$, and the
multiplicities are at most $1$ if $K=\Q(\sqrt{-d})$, $d\neq 3,4$.\\

(e)  $\chi(C(A))= 0$ always when $K\neq \Q\mbox{, }\Q(\sqrt{-d})$.
    \epr
  \proof From theorem 1.5 we can assume that $A$ is
in block-triangular form with zero under the block-diagonal. We
can also assume (by grouping similar blocks together) that the
diagonal blocks have no common eigenvalues. Construct a new matrix
$A'$ having the same block diagonal as $A$ and zeroes both above
and below the diagonal. Then $A$ and $A'$ are conjugate in $GL_m K$.
By lemma 2.5 the Euler characteristic of $A$ is a non-zero
rational multiple of the Euler characteristic of $A'$. On the
other hand, $\chi(A')$ can be expressed as the product of the
Euler characteristics of the diagonal blocks. Thus, it is enough
to examine all possible blocks that cannot be decomposed further.
Such block $B$ has eigenvalues $n$-th roots of unity repeated $k$
times. Also, the centralizer of such a block is comeasurable with
the centralizer of a matrix consisting of an $n$-torsion matrix
$T_n$ sitting in each diagonal block where $T_n$ is an $n$-torsion
matrix in $GL_l{\cal{O}}_K$. By lemma 2.7
$$\chi(C(B))=\chi(GL_l\Z[\xi_n])
  =\chi(\Z[\xi_n]^{\times})\chi(SL_l\Z[\xi_n]).$$
  The number $\chi({\cal{O}}_K[\xi_n]^{\times})$ will be zero
if the units are infinitely many, because $\chi(\Z)=0$.

For $K=\Q$ we can have only $n$-th roots of unity for 
$n=1,2,3,4,6$ such that  $$\chi({\cal{O}}_K[\xi_n]^{\times})\neq 0.$$ 
Also, we have that $$\chi(SL_l\Z)=0 \mbox{ for } k\geq 3$$ and 
$$\chi(SL_l\Z[\xi_n])=\zeta_{\Q(\xi_n)}(-1)
\dots \zeta_{\Q(\xi_n)}(-l+1)=0 \mbox{ for } 
n=3,4,6, \mbox{ and }l\geq 2.$$ 
Thus, we are left with the roots
$$\{1,-1,i,-i,\xi_3,\bar\xi_3,\xi_6,\bar\xi_6\}$$ with
multiplicities of $+1$ and $-1$ at most $2$, and multiplicity of
$\pm i,\mbox{ }\pm\xi_3$ and $\pm \xi_6$ at most $1$.

For $K=\Q(i)$, we can have only $n$-th roots of unity for 
$n=1,2,4$ such that  $$\chi(\Z[i,\xi_n]^{\times})\neq 0.$$
Also, we have that $$\zeta_{\Q(i)}(-1)=0.$$ Thus,
$$\chi(SL_l\Z[i,\xi_n])
=\zeta_{\Q(i,\xi_n)}(-1)\dots \zeta_{\Q(i,\xi_n)}(-l+1)
=0 \mbox{ for } l\geq 2.$$
Thus, we are left with the roots
$$\{1,-1,i,-i\}$$ with
multiplicities at most $1$.

For $K=\Q(\xi_3)$, we can have only $n$-th roots of unity for 
$n=1,2,3,6$ such that  $$\chi(\Z[\xi_3,\xi_n]^{\times})\neq 0.$$
Also, we have that $$\zeta_{\Q(\xi_3)}(-1)=0.$$ Thus,
$$\chi(SL_l\Z[\xi_3,\xi_n])
=\zeta_{\Q(\xi_3,\xi_n)}(-1)\dots \zeta_{\Q(\xi_3,\xi_n)}(-l+1)
=0 \mbox{ for } l\geq 2.$$
Thus, we are left with the roots
$$\{1,-1,\xi_3,\bar\xi_3,\xi_6,\bar\xi_6\}$$ with
multiplicities at most $1$.

For $K=\Q(\sqrt{-d})$, we can have only $n$-th roots of unity for 
$n=1,2$ such that  $$\chi({\cal{O}}_{\Q(\sqrt{-d},\xi_n)}^{\times})\neq 0,$$
because otherwise the units in the extension will be infinitely 
many.
Also, we have that $$\zeta_{\Q(\sqrt{-d})}(-1)=0.$$ Thus,
$$\chi(SL_k({\cal{O}}_{\Q(\sqrt{-d},\xi_n)}))
=\zeta_{\Q(\sqrt{-d},\xi_n)}(-1)\dots \zeta_{\Q(\sqrt{-d},\xi_n)}(1-l)
=0 \mbox{ for } k\geq 2.$$ 
Thus, we are left with the roots
$$\{1,-1\}$$ with
multiplicities at most $1$.

If $K$ is not $\Q$, not an imaginary quadratic extension of $\Q$
then $$\chi({\cal{O}}_{K(\xi_n)}^{\times})= 0$$
for any $n$. Thus $$\chi(C(A))=0$$ for such fields $K$.

For $SL_m {\cal{O}}_K$ we have a similar statement.
\bpr{2.9}
   Let $A$ be a torsion element of $SL_m {\cal{O}}_K$.
Then $\chi(C(A))\neq 0$ if and only if the set of eigenvalues of
$A$ is inside the set\\

(a) $\{1,-1,i,-i,\xi_3,\bar\xi_3,\xi_6,\bar\xi_6\}$, and the
multiplicity of $1$ is at most $2$, the multiplicity of $-1$ is $0$ or $2$
 and the multiplicity
of the rest of the roots of unity is at most $1$ if $K=\Q$;\\

(b) $\{1,-1,i,-i\}$, and the
multiplicities are at most $1$ if $K=\Q(i)$;\\

(c)  $\{1,-1,\xi_3,\bar\xi_3,\xi_6,\bar\xi_6\}$  and the
multiplicities are at most $1$ if $K=\Q(\xi_3)$;\\

(d) $\{\xi, \xi^{-1}\}$ where $\xi$ is a root of $1$ and the dimension of 
the matrix $A$ is at most $2$, if $K$ is totally real field different
from $\Q$.\\

(e)  $\chi(C(A))= 0$ always when $K$ is not totally real
and different from $\Q(i)$ and $\Q(\xi_3)$.
\epr
\proof Then proof is similar to the one of the previous proposition.
Given a torsion matrix $A$ in $SL_m {\cal{O}}_K$ we can conjugate it
with a matrix from $GL_m K$ to a matrix $A'$. Then $C(A)$ and $C(A')$
will be comeasurable. Consider $A$ as a matrix in $GL_m {\cal{O}}_K$.
We can assume that $A$ is of block-triangular form, using theorem 1.5.
Then we can take $A'$ to be of block-diagonal form with the same blocks 
on the diagonal as $A$. With another conjugation by a matrix from 
$GL_m K$, we can modify $A'$ so that blocks with identical eigenvalues
will be represented by the same matrices. And finally, we combine the
similar matrices into bigger block so that: (1) each new block is zero 
if it is not on the diagonal; (2) a block on the diagonal consists
of smaller blocks such that on the diagonal we have smaller 
blocks repeated as many times as  needed, and off the diagonal the 
smaller blocks are zero. The centralizer of such a matrix 
can be computed using lemma 2.6 and 2.7. We need a minor modification
of lemma 2.6: We have that $C_{GL_m {\cal{O}}_K}(T_n)$ is equal 
to the group of units in ${\cal{O}}_{K(\xi_n)}$ while 
 $C_{SL_m {\cal{O}}_K}(T_n)$ is equal 
to the group $\ker(N_{K(\xi_n)/K})$ inside the
units in ${\cal{O}}_{K(\xi_n)}$. For the 
vanishing of the Euler characteristic this does not have an effect
since both groups have the same rank. 

For $K=\Q$ we can have only $n$-th roots of unity for 
$n=1,2,3,4,6$ such that  $$\chi(\Z[\xi_n]^{\times})\neq 0.$$ 
Also, we have that $$\chi(SL_l\Z)=0 \mbox{ for } k\geq 3$$ and 
$$\chi(SL_l\Z[\xi_n])=\zeta_{\Q(\xi_n)}(-1)
\dots \zeta_{\Q(\xi_n)}(-l+1)=0 \mbox{ for } 
n=3,4,6, \mbox{ and }l\geq 2.$$ 
Thus, we are left with the roots
$$\{1,-1,i,-i,\xi_3,\xi_3^{-1},\xi_6,\xi_6^{-1}\}$$ with
multiplicities of $+1$ and $-1$ at most $2$, and multiplicity of
$\pm i,\mbox{ }\pm\xi_3$ and $\pm \xi_6$ at most $1$.

For $K=\Q(i)$, we can have only $n$-th roots of unity for 
$n=1,2,4$ such that  $$\chi(\Z[i,\xi_n]^{\times})\neq 0.$$
Also, we have that $$\zeta_{\Q(i)}(-1)=0.$$ Thus,
$$\chi(SL_l\Z[i,\xi_n])
=\zeta_{\Q(i,\xi_n)}(-1)\dots \zeta_{\Q(i,\xi_n)}(1-l)
=0 \mbox{ for } l\geq 2.$$
Thus, we are left with the roots
$$\{1,-1,i,-i\}$$ with
multiplicities at most $1$.

For $K=\Q(\xi_3)$, we can have only $n$-th roots of unity for 
$n=1,2,3,6$ such that  $$\chi(\Z[\xi_3,\xi_n]^{\times})\neq 0.$$
Also, we have that $$\zeta_{\Q(\xi_3)}(-1)=0.$$ Thus,
$$\chi(SL_l\Z[\xi_3,\xi_n])
=\zeta_{\Q(\xi_3,\xi_n)}(-1)\dots \zeta_{\Q(\xi_3,\xi_n)}(1-l)
=0 \mbox{ for } l\geq 2.$$
Thus, we are left with the roots
$$\{1,-1,\xi_3,\xi_3^{-1},\xi_6,\xi_6^{-1}\}$$ with
multiplicities at most $1$.
 
For $K=\Q(\sqrt{-d})$, we can have only $n$-th roots of unity for 
$n=1,2$ such that  $$\chi({\cal{O}}_{\Q(\sqrt{-d},\xi_n)}^{\times})\neq 0,$$
because otherwise the units in the extension will be infinitely 
many.
Also, we have that $$\zeta_{\Q(\sqrt{-d})}(-1)=0.$$ Thus,
$$\chi(SL_k{\cal{O}}_{\Q(\sqrt{-d},\xi_n)})
=\zeta_{\Q(\sqrt{-d},\xi_n)}(-1)\dots \zeta_{\Q(\sqrt{-d},\xi_n)}(1-l)
=0 \mbox{ for } k\geq 2.$$ 
We are left with the roots
$$\{1,-1\}$$ with
multiplicities at most $1$. Since $-1$ will lead to a determinant $-1$,
we have that the only with the eigenvalues $1$ with multiplicity $1$,
which gives that $A$ is in $SL_1$ and for higher $SL_m$ the Euler 
characteristic will vanish.

For $K$ a totally real field different from $\Q$, we have that
there are infinitely many units and that  
$$\zeta_{K}(-1)\neq 0,$$ which follows from the functional 
equation for the Dedekind zeta function. Also, $$\zeta_{K(\xi_n)}(-1)=0$$
for $n>2$. Thus the eigenvalues $1$ or $-1$ might occur with 
multiplisity $2$, and another root of $1$ might occur with multiplicity
at most $1$. 

Suppose $T_n$ occurs as a block in $A'$, $n>2$. Recall $T_n$ is an 
$n$-torsion block that has irreducible over $K$ characteristic polynomial.
By the modification of lemma 2.6 in the beginning of this proof
we have that the centralizer of $A'$ will contain a copy of the group of
units $\a$ in ${\cal{O}}_K$ such that $$N_{K(\xi_n)/K}(\a)=1.$$ In order
to have non-vanishing Euler characteristic the rank of this group must 
be zero. This can happen only when the ranks of the groups of units in $K$
and in $K(\xi_n)$ coincide. This is the case when  $K(\xi_n)$ is an
extension of degree $2$ over $K$. Thus, the block $T_n$ is of size
$2\times 2$. 

We are going to show that $A'$ consists of at most block $T_n$.
Suppose that the diagonal of $A'$ consists of  
$T_{n_0}$, $T_{n_1}$, $\dots$, $T_{n_l}$ on 
the diagonal, where each $T_{n_i}$ for $i=1\dots l$
is $n_i$-torsion block
with irreducible characteristic polynomial of degree $2$,
and $T_{n_0}$ is either empty, or$(\pm 1)$, or $\pm I_2$.  
Let $K_i=K(\xi_{n_i})$, and $K_0=K$. 

If $T_0$ is empty that is it is not present in $A'$ then  
the centralizer of $A'$ is isomorphic to the group
$$\{(\a_1,\dots,\a_n)|N_{K_1/K}(\a_1)\dots N_{K_l/K}(\a_l)=1\}.$$
The rank of the group is $(l-1)$ times the rank of units in $K$. In order 
to have non-vanishing Euler characteristic we must start with a zero 
rank group. Thus, $l=1$ and $A'=T_{n_i}$ for some $i$.

If $T_0=(\pm 1)$ then  
the centralizer of $A'$ is isomorphic to the group
$$\{(\a_0,\dots,\a_n)|\a_0. N_{K_1/K}(\a_1)\dots N_{K_l/K}(\a_l)=1\}.$$
The rank of the group is $l$ times the rank of units in $K$. In order 
to have non-vanishing Euler characteristic we must start with a zero 
rank group. Thus, $l=0$ and $A'=(\pm 1)$.

If $T_0=\pm I_2$ then  
the Euler characteristic of the centralizer of $A'$ is 
isomorphic to the Euler characteristic of
$$SL_2({\cal{O}}_K)\times 
\{(\a_0,\dots,\a_n)|\a_0 N_{K_1/K}(\a_1)\dots N_{K_l/K}(\a_l)= 1\}.$$
The rank of the abelian group is $l$ times the rank of units in $K$. In order 
to have non-vanishing Euler characteristic we must start with a zero 
rank group. Thus, $l=0$ and $A'=\pm I_2$. This completes the case when
$K$ is totally real.

In the rest of the cases we have that the group of units of $K$ is
infinite and that $\zeta_K(-1)=0$. Then the Euler characteristic of
the centralizer of $A'$ always vanishes.

The following theorems are very useful for computational purposes.
It expresses the homological Euler characteristic as a sum of very
few terms. And each of the terms can be easily computed. We use the 
following notation: $$A=[A_{11},\dots A_{ll}]$$
means that the square block $A_{11}$ through $A_{ll}$ are placed
on the block-diagonal of $A$ and the blocks of $A$ outside the 
block-diagonal are zero blocks. Also, let
$$R(A)=R(f_1,\dots,f_l),$$
where  $$A=[A_{11},\dots A_{ll}],$$ and $f_i$ is the characteristic
polynomial of $A_{ii}$, and $f_i$ is a power of an irreducible polynomial.
As a consequence of proposition 2.4 and 2.8 we obtain
theorems 2.10, 2.11, 2.12 and 2.13.
\bth{2.10}
   Let $V$ be a finite dimensional representation
of $GL_m \Q$. Then the homological Euler characteristic of
$GL_n\Z$ with coefficients in $V$ is given by
 $$
   \chi_h(GL_m\Z,V)=
\sum_{A} |R(A)|\chi(C(A))\tr(A^{-1}|V),
    $$
  where the sum is taken over torsion matrices $A$ consisting
of blocks $A_{11},\dots A_{mm}$ on the block-diagonal and zero
blocks off the diagonal. Also the matrices $A_{ii}$ are in the set
$\{+1,+I_2,-1,-I_2,T_3,T_4,T_6\},$ 
where 
$$T_3=
\left[
\begin{tabular}{rr}
$0$ & $1$\\
$-1$ & $-1$
\end{tabular}
\right],
\mbox{ }
T_4=
\left[
\begin{tabular}{rr}
$0$ & $1$\\
$-1$ & $0$
\end{tabular}
\right],
\mbox{ }
T_6=
\left[
\begin{tabular}{rr}
$0$ & $-1$\\
$1$ & $1$
\end{tabular}
\right].
$$
and the characteristic
polynomial $f_i$ of $A_{ii}$ is a power of an irreducilbe
polynomial, and $f_i$ and $f_j$ are relatively prime.
    \eth
\bth{2.11}
   Let $V$ be a finite dimensional representation
of $GL_m \Q(i)$. Then the homological Euler characteristic of
$GL_n(\Z[i])$ with coefficients in $V$ is given by
 $$
   \chi_h(GL_m(\Z[i]),V)=
\sum_{A} |N_{\Q(i)/\Q}(R(A))|\chi(C(A))\tr(A^{-1}|V),
    $$
  where the sum is taken over torsion matrices $A$ consisting
of blocks $A_{11},\dots A_{mm}$ on the block-diagonal and zero
blocks off the diagonal. Also the matrices $A_{ii}$ are in the set
$\{+1,-1,i,-i\}$ and the characteristic
polynomials $f_i$ of $A_{ii}$ are relatively prime if $i\neq j$.
    \eth
\bth{2.12}
   Let $V$ be a finite dimensional representation
of $GL_m \Q(\xi_3)$. Then the homological Euler characteristic of
$GL_n(\Z[\xi_3])$ with coefficients in $V$ is given by
 $$
   \chi_h(GL_m(\Z[\xi_3]),V)=
\sum_{A} |N_{\Q(\xi_3)/\Q}(R(A))|\chi(C(A))\tr(A^{-1}|V),
    $$
  where the sum is taken over torsion matrices $A$ consisting
of blocks $A_{11},\dots A_{mm}$ on the block-diagonal and zero
blocks off the diagonal. Also the matrices $A_{ii}$ are in the set
$\{+1,-1,\xi_3,\bar{\xi_3},\xi_6,\bar{xi_6}\}$ and the characteristic
polynomials $f_i$ of $A_{ii}$ are relatively prime if $i\neq j$.
    \eth
\bth{2.13}
Let $P$ be a projective module of rank $2$ over the ring of integers 
${\cal{O}}_K$ in a number field $K=\Q(\sqrt{-d})$ for $d\neq 3, 4$.
   Let $V$ be a finite dimensional representation
of $GL_2(Q(\sqrt{-d}))$ 
Then the homological Euler characteristic of $\Aut(P)$
with coefficients in $V$ is given by
$$ \chi_h(\Aut(P),V)=n\tr([1,-1]|V),$$
where $n$ is the number of ways that $P$ can be 
written as a direct sum of two projective modules, $P=P_1\oplus P_2$ 
counting also the order of the summands. In particular if $P$ 
is free we obtain the homological Euler characteristic of 
$GL_2({\cal{O}}_{\Q(\sqrt{-d})})$ with coefficients in any representation.
\eth
\proof (of theorems 2.10-2.13) Using proposition 2.4, we combine the torsion 
elements that have a common characteristic polynomial. Using 
proposition 2.8, we can take the sum only over those torsion 
elements that lead to a non-zero Euler characteristic of their 
centralizer. That leads to theorems 2.10-2.12.

For theorem 2.13 we need a little bit more. By proposition 2.8 part (d)
we know that the only block-diagonal endomorphisms that will give contribution 
to the homological Euler characteristic is $$A=[1,-1],$$ 
acting as $1$ on $P_1$ and as $-1$ on $P_2$. We have to take into account 
all decompositions of $P$ into a direct sum of rank one 
projective modules.
Note that $$R([1,-1])=2$$ and
$$N_{\Q(\sqrt{-d})/\Q}(R([1,-1]))=4.$$
Also, $$C(A)=C(1)\times C(-1)=\Z/2\times \Z/2.$$
Thus, the formula becomes
$$
\begin{tabular}{ccl}
$\chi_h(\Aut(P),V)$ & $=$ 
  & $N_{\Q(\sqrt{-d})/\Q}(R([1,-1]))n \chi(C([1,-1]))\tr([1,-1]|V)=$\\
& $=$ & $4\cdot \frac{1}{4}n\tr([1,-1]|V).$\\
\end{tabular}
$$

In the rest of the section we will deal with general statements
about the homological Euler characteristics of
groups comensurable to $GL_m({\cal{O}}_K)$ or to $SL_m({\cal{O}}_K)$
where ${\cal{O}}_K$ is a number ring.

\proof (of theorem 0.1) Let $\Ga$ be a finite subgroup of 
$GL_m ({\cal{O}}_K)$. Let $A$ be a torsion element of $\Ga$. 
Then then the centralizer of $A$ inside $\Ga$, $$C_{\Ga}(A)$$
is a finite index subgroup of the centralizer of $A$ inside
$GL_m ({\cal{O}}_K)$, $$C_{GL_m {\cal{O}}_K}(A).$$ Thus,
$$\chi(C_{\Ga}(A))=q\cdot \chi(C_{GL_m {\cal{O}}_K}(A)),$$
where $q$ is a positive integer. Thus, vanishing of 
$\chi(C_{GL_m {\cal{O}}_K}(A))$ implies vanishing of
$\chi(C_{\Ga}(A))$. We use a generalization of
Kenneth Brown's formula (see section 10,and \cite{B2}):
$$\chi(\Ga)=\sum_{A:\mbox{ torsion}}\tr(A^{-1}|V)\cdot\chi(C_{\Ga}(A)),$$
where the summation is over all torsion elements in $\Ga$
up to conjugation, and $V$ is a representation of $\Ga$

Let $K=\Q$. Let also $$A\in GL_m \Z.$$ Suppose that 
$$\chi(C(A))\neq 0.$$ Then by proposition 2.9 part (a) 
we obtain that $A$ multiplicity of the eigenvalues
$1$, $-1$ at most $2$, and multiplicity of the eigenvalues
$i$, $-i$, $\xi_3$, $\bar{\xi_3}$, $\xi_6$, $\bar{\xi_6}$
at most $1$. Thus, the dimension of $A$ is at most $10$.
Therefore if the dimension of $A$ is greater than $10$ then
the Euler characteristic of the centralizers of $A$ will vanish 
and so will the homological Euler characteristic  of a 
finite index subgroup $\Ga$ of $GL_m \Z$ with coefficients in $V$, 
$$\chi_h(\Ga,V)=0.$$

Let $K=\Q(i)$. Let also $$A\in GL_m (\Z[i]).$$ Suppose that 
$$\chi(C(A))\neq 0.$$ Then by proposition 2.9 part (b) 
we obtain that $A$ multiplicity of the eigenvalues
$1$, $-1$, $i$, $-i$ 
is at most $1$. Thus, the dimension of $A$ is at most $4$.
Therefore if the dimension of $A$ is greater than $4$ then
the Euler characteristic of the centralizers of $A$ will vanish 
and so will the homological Euler characteristic  of a 
finite index subgroup $\Ga$ of $GL_m (\Z[i])$, $$\chi_h(\Ga)=0.$$

Let $K=\Q(\xi_3)$. Let also $$A\in GL_m (\Z[\xi_3]).$$ Suppose that 
$$\chi(C(A))\neq 0.$$ Then by proposition 2.9 part (c) 
we obtain that $A$ multiplicity of the eigenvalues
$1$, $-1$, $\xi_3$, $\bar{\xi_3}$, $\xi_6$, $\bar{\xi_6}$
at most $1$. Thus, the dimension of $A$ is at most $6$.
Therefore if the dimension of $A$ is greater than $6$ then
the Euler characteristic of the centralizers of $A$ will vanish 
and so will the homological Euler characteristic of a 
finite index subgroup $\Ga$ of $GL_m (\Z[\xi_3])$ with coefficients in $V$, 
$$\chi_h(\Ga,V)=0.$$

Let $K=\Q(\sqrt{-d})$, $d\neq 3,4$.
Let also $$A\in GL_m ({\cal{O}}_{\Q(\sqrt{-d})}).$$ 
Suppose that 
$$\chi(C(A))\neq 0.$$ Then by proposition 2.9 part (d) 
we obtain that $A$ multiplicity of the eigenvalues
$1$, $-1$, 
is at most $1$. Thus, the dimension of $A$ is at most $2$.
Therefore if the dimension of $A$ is greater than $2$ then
the Euler characteristic of the centralizers of $A$ will vanish 
and so will the homological Euler characteristic of a 
finite index subgroup $\Ga$ of $GL_m ({\cal{O}}_{\Q(\sqrt{-d})})$
with coefficients in $V$, 
$$\chi_h(\Ga,V)=0.$$

And finally, for all other fields $K$ by proposition 2.9 part (e),
we have that always $$\chi(C(A))=0.$$ Therefore, the homological 
Euler characteristic  of a 
finite index subgroup $\Ga$ of $GL_m ({\cal{O}}_K)$ with coefficients in $V$, 
$$\chi_h(\Ga,V)=0.$$

Similarly we obtain the vanishing result for
finite index subgroups of $SL_m ({\cal{O}}_K)$.

\proof (of theorem 0.2)  Let $\Ga$ be a finite subgroup of 
$SL_m ({\cal{O}}_K)$. Let $A$ be a torsion element of $\Ga$. 
Then then the centralizer of $A$ inside $\Ga$, $$C_{\Ga}(A)$$
is a finite index subgroup of the centralizer of $A$ inside
$SL_m ({\cal{O}}_K)$, $$C_{SL_m {\cal{O}}_K}(A).$$ Thus,
$$\chi(C_{\Ga}(A))=q\cdot \chi(C_{SL_m {\cal{O}}_K}(A)),$$
where $q$ is a positive integer. Thus, vanishing of 
$\chi(C_{SL_m {\cal{O}}_K}(A))$ implies vanishing of
$\chi(C_{\Ga}(A))$. We use a generalization of
Kenneth Brown's formula (see section 10,and \cite{B2}):
$$\chi(\Ga)=\sum_{A:\mbox{ torsion}}\tr(A^{-1}|V)\cdot\chi(C_{\Ga}(A)),$$
where the summation is over all torsion elements in $\Ga$
up to conjugation, and $V$ is a representation of $\Ga$

Let $K=\Q$. Let also $$A\in SL_m \Z.$$ Suppose that 
$$\chi(C(A))\neq 0.$$ Then by proposition 2.9 part (a) 
we obtain that $A$ multiplicity of the eigenvalues
$1$, $-1$ at most $2$, and multiplicity of the eigenvalues
$i$, $-i$, $\xi_3$, $\bar{\xi_3}$, $\xi_6$, $\bar{\xi_6}$
at most $1$. Thus, the dimension of $A$ is at most $10$.
Therefore if the dimension of $A$ is greater than $10$ then
the Euler characteristic of the centralizers of $A$ will vanish 
and so will the homological Euler characteristic  of a 
finite index subgroup $\Ga$ of $SL_m \Z$ with coefficients in $V$, 
$$\chi_h(\Ga,V)=0.$$

Let $K=\Q(i)$. Let also $$A\in SL_m (\Z[i]).$$ Suppose that 
$$\chi(C(A))\neq 0.$$ Then by proposition 2.9 part (b) 
we obtain that $A$ multiplicity of the eigenvalues
$1$, $-1$, $i$, $-i$ 
is at most $1$. Thus, the dimension of $A$ is at most $4$.
Therefore if the dimension of $A$ is greater than $4$ then
the Euler characteristic of the centralizers of $A$ will vanish 
and so will the homological Euler characteristic  of a 
finite index subgroup $\Ga$ of $SL_m (\Z[i])$ with coefficients in $V$, 
$$\chi_h(\Ga,V)=0.$$

Let $K=\Q(\xi_3)$. Let also $$A\in SL_m (\Z[\xi_3]).$$ Suppose that 
$$\chi(C(A))\neq 0.$$ Then by proposition 2.9 part (c) 
we obtain that $A$ multiplicity of the eigenvalues
$1$, $-1$, $\xi_3$, $\bar{\xi_3}$, $\xi_6$, $\bar{\xi_6}$
at most $1$. Thus, the dimension of $A$ is at most $6$.
Therefore if the dimension of $A$ is greater than $6$ then
the Euler characteristic of the centralizers of $A$ will vanish 
and so will the homological Euler characteristic of a 
finite index subgroup $\Ga$ of $SL_m (\Z[\xi_3])$ with coefficients in $V$, 
$$\chi_h(\Ga,V)=0.$$

Let $K$ be totally real field..
Let also $$A\in GL_m ({\cal{O}}_K).$$ 
Suppose that 
$$\chi(C(A))\neq 0.$$ Then by proposition 2.9 part (d) 
we obtain that the dimension of $A$ is at most $2$. 
Therefore if the dimension of $A$ is greater than $2$ then
the Euler characteristic of the centralizers of $A$ will vanish 
and so will the homological Euler characteristic of a 
finite index subgroup $\Ga$ of $GL_m ({\cal{O}}_{\Q(\sqrt{-d})})$
with coefficients in $V$, 
$$\chi_h(\Ga,V)=0.$$

And finally, for all other fields $K$ by proposition 2.9 part (e),
we have that always $$\chi(C(A))=0.$$ Therefore, the homological 
Euler characteristic  of a 
finite index subgroup $\Ga$ of $SL_m ({\cal{O}}_K)$ with coefficients in $V$, 
$$\chi_h(\Ga,V)=0.$$

\sectionnew{Torsion elements of $GL_2\Z$ and $GL_3\Z$}
This section deals with the computational task of finding
all the torsion element up to conjugation in $GL_2 \Z$ and $GL_3\Z$.
\bpr{4.1}
    All torsion elements in $GL_2\Z$
up to conjugation are listed in the table below together with
their centralizers and Euler characteristic of the centralizers.\\
   $
\begin{tabular}{|l|c|c|}
\hline $A$                          & $C(A)$   & $\chi(C(A))$\\
\hline (a) $\left[\begin{tabular}{rr} $1$ &     \\
    & $1$ \\
\end{tabular}\right]$        & $GL_2\Z$ & $-\frac{1}{24}$\\ \hline
(b) $\left[\begin{tabular}{rr} $-1$ &      \\
     & $-1$ \\
\end{tabular}\right]$        & $GL_2\Z$ & $-\frac{1}{24}$\\ \hline
(c1) $\left[\begin{tabular}{rr} $1$ &     \\
    & $-1$ \\
\end{tabular}\right]$        & $C_2\times C_2$ & $-\frac{1}{4}$\\ \hline
(c2) $\left[\begin{tabular}{rr} $1$ &  $1$ \\
    & $-1$ \\
\end{tabular}\right]$        & $C_2\times C_2$ & $-\frac{1}{4}$\\ \hline
(d) $\left[\begin{tabular}{rr} $0$ & $1$    \\ $-1$ & $-1$  \\
\end{tabular}\right]$        & $C_6$   & $\frac{1}{6}$\\ \hline
(e) $\left[\begin{tabular}{rr} $0$ & $-1$ \\ $1$ & $1$  \\
\end{tabular}\right]$        & $C_6$ & $\frac{1}{6}$\\ \hline
(f) $\left[\begin{tabular}{rr} $0$  & $1$ \\ $-1$ & $0$ \\
\end{tabular}\right]$        & $C_4$ & $\frac{1}{4}$\\ \hline
\end{tabular}
$
    \epr
   \proof Let $A$ be a torsion element in $GL_2\Z$. If the
characteristic polynomial of $A$ is irreducible over the rational
numbers then its root belong to a quadratic extension of $\Q$. The
only such options occur if the roots are $3$-rd, $4$-th or $6$-th
root of unity. The corresponding rings of integers are unique
factorization domains. Therefore by proposition 1.3 there is only
one matrix $T_n$ up to conjugation having eigenvalues $n$-th roots
of unity for $n=3$, $4$, or $6$. This gives the torsion
elements in part (d), (e) and (f).If the characteristic polynomial
of $A$ is reducible then its roots are $1$ or $-1$. If both roots
are $1$ or both are $-1$ then we have only one representative for
each case which is given in part (a) and (b). It remains to
consider the case when the eigenvalues of $A$ are $1$ and $-1$. By
theorem 1.5 we can assume that $A$ is upper triangular. What
remains to be done is to determine
 which elements in the upper right corner of $A$ give conjugate
elements. Consider $P_{A_{11}A_{22}}$ with $A_{11}=1$ and 
$A_{22}=-1$. Then the matrix representing $P_{A_{11}A_{11}}$ is
$$I\otimes A_{22} -A_{11}\otimes I=1\otimes (-1) -1\otimes 1=-2.$$
Thus,
$$\left[\begin{tabular}{cr} $1$  & $a$ \\ $0$ & $-1$ \\
\end{tabular}\right]
\mbox{ and } \left[\begin{tabular}{cr} $1$  & $b$ \\ $0$ & $-1$ \\
\end{tabular}\right]$$
   are conjugate if and only if $a\equiv b\mbox{ }\mod\mbox{ }2$. Thus
we obtain the two torsion elements in (c1) and (c2) corresponding
to even $a$ and odd $a$.

For $GL_3 \Z$ we find all the torsion elements up
to conjugation.


\bpr{4.2}
    All torsion elements in $GL_3\Z$
up to conjugation are listed in the table below together with
their centralizers and Euler characteristic of the centralizers.\\
$
\begin{tabular}{|c|c|c|c|}
\hline $A$                          & $-A$   & $C(A)$   &
$\chi(C(A))$\\ \hline (a) $\left[\begin{tabular}{ccc} $1$ &     &
\\
    & $1$ & \\
    &     & $1$\\
\end{tabular}\right]$        &
  (b)
$\left[\begin{tabular}{ccc} $-1$ &      &\\
     & $-1$ &\\
     &      & $-1$\\
\end{tabular}\right]$        & $GL_3\Z$ & $0$\\ \hline
(c1) $\left[\begin{tabular}{cc|c} $1$ & $0$ & \\ $0$ & $1$ & \\
\hline
    &     & $-1$\\
\end{tabular}\right]$        &

(d1) $\left[\begin{tabular}{cc|c} $-1$ & $0$ & \\ $0$ & $-1$ & \\
\hline
    &      & $1$\\
\end{tabular}\right]$        & $GL_2\Z\times GL_1\Z$ & $-\frac{1}{48}$\\ \hline

(c2) $\left[\begin{tabular}{cc|c} $1$ & $0$ & $1$\\ $0$ & $1$ &
$0$\\ \hline
    &     & $-1$\\
\end{tabular}\right]$        &

(d2) $\left[\begin{tabular}{cc|c} $-1$ & $0$ & $-1$\\ $0$ & $-1$ &
$0$\\ \hline
    &      & $1$\\
\end{tabular}\right]$
  & $\Gamma_1(2,2)\times GL_1\Z$ & $-\frac{1}{16}$\\ \hline

(e1)
   $\left[\begin{tabular}{cc|c} $0$ & $1$   & \\ $-1$ & $-1$ &
\\ \hline
    &     & $1$\\
\end{tabular}\right]$        &

(f1)
  $\left[\begin{tabular}{cc|c}
$0$ & $-1$ & \\ $1$ & $1$  & \\ \hline
    &      & $-1$\\
\end{tabular}\right]$        & $C_6\times C_2$ & $\frac{1}{12}$\\ \hline

(e2)
  $\left[\begin{tabular}{cc|c}
$0$ & $1$   & $1$\\ $-1$ & $-1$ & \\ \hline
    &     & $1$\\
\end{tabular}\right]$        &

(f2)
   $\left[\begin{tabular}{cc|c}
   $0$ & $-1$ & $-1$ \\
$1$ & $1$  & \\ \hline
    &      & $-1$\\
\end{tabular}\right]$        & $C_3\times C_2$ & $\frac{1}{6}$\\ \hline

(g)
   $\left[\begin{tabular}{cc|c}
$0$  & $1$  & \\ $-1$ & $-1$ & \\ \hline
     &      & $-1$\\
\end{tabular}\right]$        &

(h)
    $\left[\begin{tabular}{cc|c}
    $0$ & $-1$ & \\
  $1$ & $1$  & \\ \hline
    &      & $1$\\
\end{tabular}\right]$        & $C_6\times C_2$ & $\frac{1}{12}$\\ \hline

(i1)
   $\left[\begin{tabular}{cc|c}
$0$  & $1$ & \\ $-1$ & $0$ & \\ \hline
     &     & $1$\\
\end{tabular}\right]$        &

(j1)
    $\left[\begin{tabular}{cc|c}
$0$  & $-1$ & \\ $1$ & $0$ & \\ \hline
     &     & $-1$\\
\end{tabular}\right]$        & $C_4\times C_2$ & $\frac{1}{8}$\\ \hline

(i2)
   $\left[\begin{tabular}{cc|c}
$0$  & $1$ & $1$\\ $-1$ & $0$ & \\ \hline
     &     & $1$\\
\end{tabular}\right]$        &

(j2)
   $\left[\begin{tabular}{cc|c}
$0$  & $-1$ & $-1$\\
$1$ & $0$ & \\ \hline
     &     & $-1$\\
\end{tabular}\right]$        & $C_4\times C_2$ & $\frac{1}{8}$\\ \hline
\end{tabular}
   $

\epr
    \proof If $A$ is a torsion element in $GL_3\Z$
then its eigenvalues are roots of $1$. If $\l$ is an eigenvalue
then all of its Galois conjugates are eigenvalues. If $\l$ is an
$n$-th root of $1$ then all its  Galois conjugates are all the
primitive $n$-th roots of $1$. Their number is $\varphi(n)$ and we
must have at must $3$ of them. The inequality $\varphi(n)\leq 3$
has solutions $n=1,2,3,4,6$. In all of these cases we have
$\varphi(n)\leq 2$.That is, $\Q(\l)$ is at most quadratic
extension of $\Q$ Thus, the remaining eigenvalue must be rational.
Therefore the matrix $A$ must have an eigenvalue $+1$ or $-1$.

If all the eigenvalues are $1$ then the matrix is either the
identity or it contains a non-trivial Jordan block, which cannot
be of finite order. That gives case (a).Minus that matrix gives
case (b). If $A$ the characteristic polynomial of $A$ has a root
at $-1$ and a double root at $1$ then we need to use the material
that we developed so far. Using theorem 1.6 we know that the
matrix $A$ can be conjugated to a block-triangular matrix with a
zero block under the diagonal. The diagonal block must be
$A_{11}=I_2$ and $A_{22}=-1$. Then $$P_{A_{11},A_{22}}=I_2\otimes
A_{22}- A_{11}\otimes 1=-2I_2$$ acts on $Mat_{2,1}\Z$. Therefore
$$P_{mod}=2Mat_{2,1}\Z$$ and $$Q_{mod}=Mat_{2,1}(\Z/2).$$ The
centralizer of $A_{11}=I_2$, $C(A_{11})=GL_2\Z$ acts on $Q_{mod}$
and has two orbits: the zero vector and the non-zero vectors.
Therefore we have exactly two non-conjugate matrices up to
conjugation with diagonal blocks $I_2$ and $-1$. We assumed that
below the diagonal the block is zero. Above the diagonal there are
two cases leading to (c1) and (c2). In one of them we have zero
above the diagonal, and in the other a non-zero representative of
$Q_{mod},$ for example $[1, \: 0]^t$.

   By lemma 2.1 we obtain the
centralizer in the case of (c1). For the case (c2), we have that
the centralizer of $A_{22}=-1$ is $$C(A_{22})=\pm1,$$ and $-1$ acts
trivially on $$[1, \: 0]^t \e Q_{mod}.$$ On the other hand
$C(A_{11})=GL_2\Z$ acts on $Q_{mod}$ as $GL_2(\Z/2)$. And the
stabilizer of $[1, \: 0]^t$ is $\Gamma_1(2,2)$. We have that
$$\chi(\Gamma_1(2,2))=[\Gamma_1(2,2):GL_2\Z]\cdot \chi(GL_2\Z)=
3\cdot(-\frac{1}{24})=-\frac{1}{8}.$$ Considering $-A$, we obtain
the cases (d1) and (d2).

     In case not all the eigenvalues of $A$
are $\pm1$, we must have either third, fourth or sixth root of
$1$. Together with them the third eigenvalue is either $1$ or
$-1$. That exhausts all cases for the eigenvalues of a torsion
element in $GL_3\Z$.

  Suppose the eigenvalues of $A$ are
primitive third roots of $1$ and $1$. We can assume that the
matrix $A$ is in block-triangular form. Let the blocks be

$$A_{11}= \left[
\begin{tabular}{cc} $0$ & $1$\\ $-1$ & $-1$\\
\end{tabular}
\right]
 \: \mbox{and } \: A_{22}=1.$$
For $A_{11}$ we could have picked any other matrix with
characteristic polynomial $$t^2+t+1.$$ In $\Z[\xi_3]$ all ideals are
principal. Now using corollary 1.4, we obtain that any two
matrices with characteristic polynomial $t^2+t+1$ are conjugate to
each other inside $GL_2\Z$. We assumed that $A_{21}=0$. For the
simplification of $A_{12}$ we use $$P_{A_{11}A_{22}}=I_2\otimes
A_{22}-A_{11}\otimes 1.$$ Therefore,

 $$
   P_{A_{11}A_{22}}=
\left[
\begin{tabular}{cc} $1$ & $-1$\\ $1$ & $2$\\
\end{tabular}
\right].
   $$
Then $$\det P_{A_{11}A_{22}}=3$$ which means that $$Q_{mod}\cong C_3.$$
Then $C(A_{22})=\pm1$ acts on $Q_{mod}$ by exchanging the two
nonzero elements. Thus there are two orbits on $Q_{mod}$. A
representative of the zero orbit is $[0,\: 0]^t$ and for the
non-zero orbit $[1,\: 0]^t$. Thus, we obtain the cases (e1) and
(e2). By lemma 2.1 we obtain the centralizer of (e1). For the
centralizer of (e2), we observe that $-1$ from $C(A_{22})$ and
$-I_2$ from $C(A_{11})$ do change independently the element $$[1,\:
0]^t\e Q_{mod}$$ to $[-1,\: 0]^t$. However, if they act
simultaneously, they keep the element $[1,\: 0]^t$ fixed. This is
an element of order $2$ that fixes $[1,\: 0]^t$. An element of
order $3$ from $C(A_{11})$ must fix $[1,\: 0]^t$.  Thus, the
centralizer is $C_3\times C_2$. Considering $-A$ we obtain the six
torsion from (f1) and (f2).

     If the eigenvalues of $A$ are
primitive third roots of $1$ and $-1$ then we can again assume
that $A$ is in block-diagonal form with block

     $$
     A_{11}= \left[
\begin{tabular}{cc} $0$ & $1$\\ $-1$ & $-1$\\
\end{tabular}
\right]
 \: \mbox{and } \: A_{22}=-1.
   $$
    Then $P_{A_{11}A_{22}}=I_2\otimes A_{22}-A_{11}\otimes 1$.
Therefore,

$$
   P_{A_{11}A_{22}}= \left[
\begin{tabular}{cc} $-1$ & $-1$\\ $1$ & $0$\\
\end{tabular}
\right].
   $$
Then $$\det P_{A_{11}A_{22}}=1$$ which means that $Q_{mod}=0$.
Therefore, there is only one conjugacy class with such $A_{11}$
and  $A_{22}$, and we can take for $A_{12}$ any vector. So we take
$[0,\: 0]^t$. The centralizer of $A$ is a product of two
centralizers by lemma 2.1. That concludes case (g). Case (h) is
taking $-A$ in the previous case. Thus, we have exhausted
all $3$-torsions and all $6$-torsions in $GL_3\Z$.

Suppose we
have a matrix $A$ with eigenvalues $i$, $-i$ and $1$. (Taking the
minus sign we obtain a matrix with eigenvalues $i$, $-i$ and
$-1$.) Using theorem 1.3, we can assume that $A$ is in
block-triangular form with zero below the diagonal block, and
$$
    A_{11}= \left[
\begin{tabular}{cc} $0$ & $1$\\ $-1$ & $0$\\
\end{tabular}
\right]
 \: \mbox{and } \: A_{22}=1.$$
   Then $$P_{A_{11}A_{22}}=I_2\otimes A_{22}-A_{11}\otimes 1.$$
Therefore,
$$
  P_{A_{11}A_{22}}=
\left[
\begin{tabular}{cc} $1$ & $-1$\\ $1$ & $1$\\
\end{tabular}
\right].
    $$
Then $$\det P_{A_{11}A_{22}}=2$$ which means that $Q_{mod}=\Z/2$.
The module $Q_{mod}$ has no automorphisms therefore each of its
elements correspond to a conjugacy class. Thus there are two
classes: (i1) and (i2). The centralizer of (i1) can be computed
via Lemma 2.1. For (i2) note that every element of $C(A_{11})$ or
$C(A_{22})$ fixes the non-zero element of $Q_{mod}$. Therefore the
centralizer in the case of (i2) is isomorphic as an abstract group
to the centralizer of (i1). Taking the minus sign we obtain (j1)
and (j2).



\sectionnew{Resultants}


This section is computational. It deals with resultants and Euler
characteristics of centralizers needed for computation of 
homological Euler characteristics. Let us recall the notation
that we are using. By $A=[A_{11},\dots A_{kk}]$ we mean a matrix $A$,
whose diagonal blocks are $A_{ii}$ for $i=1,\dots k.$ We also assume that
the characteristic polynomial $f_i$ of $A_{ii}$ is a power of an
irreducible polynomial. Also,
$$R(A)=R(f_1,\dots , f_k)=\prod_{i<j}R(f_i,f_j),$$
where $R(f_i,f_j)$ is the resultant of $f_i$ and $f_j$.

In the following statements we compute $$|R(A)|\cdot \chi(A)$$ for 
various matrices $A$ for the need of computation of the 
homological Euler characteristic of an arithmetic group over $\Z$, 
$\Z[i]$ and $\Z[\xi_3]$
(see theorem 2.10).

\ble{5.1}
    The resultants needed for the homological Euler
characteristic of arithmetic subgroups of $GL_2\Z$ are given by:
 $$|R([1,-1])|\cdot \chi(C([1,-1]))=\frac{1}{2}.$$
\ele
\proof We have $$R([1,-1])=R(t-1,t+1)=2$$ and 
$$C([1,-1])=C([1])\times C([-1])=C_2\times C_2.$$ 
Then,
$$\chi(C([1,-1]))=\frac{1}{4}.$$ 

\ble{5.2}
    The resultants needed for the homological Euler
characteristic of arithmetic subgroups of $GL_3\Z$ are given by:\\

(c) $|R([I_2,-1])|\cdot \chi(C([I_2,-1]))=-\frac{1}{12},$\\

(d) $|R([-I_2,1])|\cdot \chi(C([-I_2,1]))=\frac{1}{12},$\\

(e) $|R([T_3,1])|\cdot \chi(C([T_3,1]))=\frac{1}{4},$\\

(f) $|R([T_6,1])|\cdot \chi(C([T_6,1]))=\frac{1}{12},$\\

(g) $|R([T_3,-1])|\cdot \chi(C([T_3,-1]))=\frac{1}{12},$\\

(h) $|R([T_6,-1])|\cdot \chi(C([T_6,-1]))=\frac{1}{4},$\\

(i) $|R([T_4,1])|\cdot \chi(C([T_4,1]))=\frac{1}{4},$\\

(j) $|R([T_4,-1])|\cdot \chi(C([T_4,-1]))=\frac{1}{4}.$\\
The enumeration follows the one of the table in proposition 3.2.
    \ele
\proof (c) We have $$R([I_2,-1])=R((t-1)^2,t+1)=4$$ and 
$$C([I_2,-1])=C([I_2])\times C([-1])=GL_2\Z\times C_2.$$ 
Then,
$$\chi(C([I_2,-1]))=-\frac{1}{48}.$$ 
(d) We have $$R([-I_2,1])=R((t+1)^2,t-1)=4$$ and 
$$C([-I_2,1])=C([-I_2])\times C([1])=GL_2\Z\times C_2.$$ 
Then,
$$\chi(C([-I_2,1]))=-\frac{1}{48}.$$ 
(e) We have $$R([T_3,1])=R((t-\xi_3)(t-\xi_3^{-1}),t-1)=3$$ and 
$$C([T_3,1])=C([T_3])\times C([1])=C_6\times C_2.$$ 
Then,
$$\chi(C([T_3,1]))=-\frac{1}{12}.$$
(f) We have $$R([T_6,1])=R((t-\xi_6)(t-\xi_6^{-1}),t-1)=1$$ and 
$$C([T_6,1])=C([T_6])\times C([1])=C_6\times C_2.$$ 
Then,
$$\chi(C([T_6,1]))=-\frac{1}{12}.$$
(g) We have $$R([T_3,-1])=R((t-\xi_3)(t-\xi_3^{-1}),t+1)=1$$ and 
$$C([T_3,-1])=C([T_3])\times C([-1])=C_6\times C_2.$$ 
Then,
$$\chi(C([T_3,-1]))=-\frac{1}{12}.$$
(h) We have $$R([T_6,-1])=R((t-\xi_6)(t-\xi_6^{-1}),t+1)=3$$ and 
$$C([T_6,-1])=C([T_6])\times C([1])=C_6\times C_2.$$ 
Then,
$$\chi(C([T_6,1]))=-\frac{1}{12}.$$
(i) We have $$R([T_6,1])=R((t-i)(t+i),t-1)=2$$ and 
$$C([T_4,1])=C([T_4])\times C([1])=C_4\times C_2.$$ 
Then,
$$\chi(C([T_4,1]))=-\frac{1}{8}.$$
(j) We have $$R([T_4,-1])=R((t-i)(t+i),t+1)=2$$ and 
$$C([T_4,-1])=C([T_4])\times C([-1])=C_4\times C_2.$$ 
Then,
$$\chi(C([T_4,-1]))=-\frac{1}{8}.$$

\ble{5.3}
    The resultants needed for the homological Euler
characteristic of arithmetic subgroups of $GL_4\Z$ are given by:\\

(a) $|R([I_2,-I_2])|\cdot \chi(C([I_2,-I_2]))=\frac{1}{36},$\\

(b) $|R([I_2,T_3])|\cdot \chi(C([I_2,T_3]))=-\frac{1}{16},$\\

(c) $|R([I_2,T_6])|\cdot \chi(C([I_2,T_6]))=-\frac{1}{144},$\\

(d) $|R([I_2,T_4])|\cdot \chi(C([I_2,T_4]))=-\frac{1}{24},$\\

(e) $|R([-I_2,T_3])|\cdot \chi(C([-I_2,T_3]))=-\frac{1}{144},$\\

(f) $|R([-I_2,T_6])|\cdot \chi(C([-I_2,T_6]))=-\frac{1}{16},$\\

(g) $|R([-I_2,T_4])|\cdot \chi(C([-I_2,T_4]))=-\frac{1}{24},$\\

(h) $|R([1,-1,T_3])|\cdot \chi(C([1,-1,T_3]))=\frac{1}{4},$\\

(i) $|R([1,-1,T_6])|\cdot \chi(C([1,-1,T_6]))=\frac{1}{4},$\\

(j) $|R([1,-1,T_4])|\cdot \chi(C([1,-1,T_4]))=\frac{1}{2},$\\

(k) $|R([T_3,T_6])|\cdot \chi(C([T_3,T_6]))=\frac{1}{9},$\\

(l) $|R([T_3,T_4])|\cdot \chi(C([T_3,T_4]))=\frac{1}{24},$\\

(m) $|R([T_6,T_4])|\cdot \chi(C([T_6,T_4]))=\frac{1}{24}.$
       \ele
\proof (a) We have 
$$R([I_2,-I_2])=R((t-1)^2,(t+1)^2)=2^4.$$
Also,
$$C(I_2,-I_2)\cong GL_2\Z \times GL_2\Z.$$
Then
$$\chi(C([I_2,-I_2]))=(-\frac{1}{24})^2.$$
We obtain
$$|R([I_2,-I_2])|\chi(C([I_2,-I_2]))=2^4(-\frac{1}{24})^2=\frac{1}{36}.$$
(b) We have 
$$R([I_2,T_3])=R((t-1)^2,(t-\xi_3)(t-\xi_3^{-1}))=3^2.$$
Also,
$$C(I_2,T_3)\cong GL_2\Z \times C_6.$$
Then
$$\chi(C([I_2,T_3]))=-\frac{1}{24}\cdot\frac{1}{6}.$$
We obtain
$$|R([I_2,T_3])|\chi(C([I_2,T_3]))=3^2(-\frac{1}{24})\frac{1}{6}
=-\frac{1}{16}.$$
(c) We have 
$$R([I_2,T_6])=R((t-1)^2,(t-\xi_6)(t-\xi_6^{-1}))=1.$$
Also,
$$C(I_2,T_6)\cong GL_2\Z \times C_6.$$
Then
$$\chi(C([I_2,T_6]))=-\frac{1}{24}\cdot\frac{1}{6}.$$
We obtain
$$|R([I_2,T_6])|\chi(C([I_2,T_6]))=(-\frac{1}{24})\frac{1}{6}
=-\frac{1}{144}.$$
(d) We have 
$$R([I_2,T_4])=R((t-1)^2,(t-i)(t+i))=2^2.$$
Also,
$$C(I_2,T_4)\cong GL_2\Z \times C_4.$$
Then
$$\chi(C([I_2,T_4]))=-\frac{1}{24}\cdot\frac{1}{4}.$$
We obtain
$$|R([I_2,T_4])|\chi(C([I_2,T_4]))=2^2(-\frac{1}{24})\frac{1}{4}
=-\frac{1}{24}.$$
(e) We have 
$$R([-I_2,T_3])=R((t+1)^2,(t-\xi_3)(t-\xi_3^{-1}))=1.$$
Also,
$$C(-I_2,T_3)\cong GL_2\Z \times C_6.$$
Then
$$\chi(C([-I_2,T_3]))=-\frac{1}{24}\cdot\frac{1}{6}.$$
We obtain
$$|R([-I_2,T_3])|\chi(C([-I_2,T_3]))=-\frac{1}{24}\cdot\frac{1}{6}
=-\frac{1}{144}.$$
(f) We have 
$$R([-I_2,T_6])=R((t+1)^2,(t-\xi_6)(t-\xi_6^{-1}))=3^2.$$
Also,
$$C(-I_2,T_6)\cong GL_2\Z \times C_6.$$
Then
$$\chi(C([-I_2,T_6]))=-\frac{1}{24}\cdot\frac{1}{6}.$$
We obtain
$$|R([-I_2,T_6])|\chi(C([-I_2,T_6]))=3^2(-\frac{1}{24})\frac{1}{6}
=-\frac{1}{16}.$$
(g) We have 
$$R([-I_2,T_4])=R((t+1)^2,(t-i)(t+i))=2^2.$$
Also,
$$C(-I_2,T_4)\cong GL_2\Z \times C_4.$$
Then
$$\chi(C([-I_2,T_4]))=-\frac{1}{24}\cdot\frac{1}{4}.$$
We obtain
$$|R([-I_2,T_4])|\chi(C([-I_2,T_4]))=2^2(-\frac{1}{24})\frac{1}{4}
=-\frac{1}{24}.$$
(h) We have 
$$R([1,-1,T_3])=R(t-1,t+1,(t-\xi_3)(t-\xi_3^{-1}))=2\cdot 3\cdot1.$$
Also,
$$C(1,-1,T_3)\cong C_2\times C_2 \times C_6.$$
Then
$$\chi(C([1,-1,T_3]))=\frac{1}{2}\cdot\frac{1}{2}\cdot\frac{1}{6}.$$
We obtain
$$|R([-I_2,T_3])|\chi(C([-I_2,T_3]))=2\cdot 3\frac{1}{2}\cdot\frac{1}{2}
\cdot\frac{1}{6}=\frac{1}{4}.$$
(i) We have 
$$R([1,-1,T_6])=R(t-1,t+1,(t-\xi_6)(t-\xi_6^{-1}))=2\cdot 1\cdot 3.$$
Also,
$$C(1,-1,T_6)\cong C_2\times C_2 \times C_6.$$
Then
$$\chi(C([1,-1,T_6]))=\frac{1}{2}\cdot\frac{1}{2}\cdot\frac{1}{6}.$$
We obtain
$$|R([1,-1,T_6])|\chi(C([1,-1,T_6]))=2\cdot 3\frac{1}{2}\cdot\frac{1}{2}
\cdot\frac{1}{6}=\frac{1}{4}.$$
(j) We have 
$$R([1,-1,T_4])=R(t-1,t+1,(t-i)(t+i))=2^3.$$
Also,
$$C(1,-1,T_4)\cong C_2\times C_2 \times C_4.$$
Then
$$\chi(C([1,-1,T_4]))=\frac{1}{2}\cdot\frac{1}{2}\cdot\frac{1}{4}.$$
We obtain
$$|R([1,-1,T_4])|\chi(C([1,-1,T_4]))=2^3\frac{1}{2}\cdot\frac{1}{2}
\cdot\frac{1}{4}=\frac{1}{2}.$$
(k) We have 
$$R([T_3,T_6])=R((t-\xi_3)(t-\xi_3^{-1}),(t-\xi_6)(t-\xi_6^{-1}))=
2^2.$$
Also,
$$C(T_3,T_6)\cong C_6 \times C_6.$$
Then
$$\chi(C([T_3,T_6]))=\cdot\frac{1}{6}\cdot\frac{1}{6}.$$
We obtain
$$|R([T_3,T_6])|\chi(C([T_3,T_6]))=2^2\frac{1}{6}
\cdot\frac{1}{6}=\frac{1}{9}.$$
(l) We have 
$$R([T_3,T_4])=R((t-\xi_3)(t-\xi_3^{-1}),(t-i)(t+i))=1.$$
Also,
$$C(T_3,T_4)\cong C_6 \times C_4.$$
Then
$$\chi(C([T_3,T_4]))=\cdot\frac{1}{6}\cdot\frac{1}{4}.$$
We obtain
$$|R([T_3,T_4])|\chi(C([T_3,T_4]))=\frac{1}{6}\cdot\frac{1}{4}
=\frac{1}{24}.$$
(m) We have 
$$R([T_6,T_4])=R((t-\xi_6)(t-\xi_6^{-1}),(t-i)(t+i))=1.$$
Also,
$$C(T_6,T_4)\cong C_6 \times C_4.$$
Then
$$\chi(C([T_6,T_4]))=\cdot\frac{1}{6}\cdot\frac{1}{4}.$$
We obtain
$$|R([T_6,T_4])|\chi(C([T_6,T_4]))=\frac{1}{6}\cdot\frac{1}{4}
=\frac{1}{24}.$$


From theorem 2.8 part (b) we have that
the only torsion element in $GL_2\Z[i]$ are these element
whose Euler characteristic of their centralizer is not zero are
the ones whose eigenvalues are different and belong to the set
$$\{\pm1.\pm i\}.$$


\ble{8.3} This is the computation of $|N_{\Q(i)/\Q}(R(T))|\cdot \chi(C(T))$ for
all torsion elements up to conjugation $T$ in $GL_2\Z[i]$ such
that $\chi(C(T))\neq 0$.\\

(a) $R([i^k,i^{k+1}])\chi(C([i^k,i^{k+1}]))=\frac{1}{8},$ for $k=0,1,2,3,$\\

(b) $R([i^k,i^{k+2}])\chi(C([i^k,i^{k+2}]))=\frac{1}{4},$ for
$k=0,1.$
  \ele
  \proof If $k$ is not congruent to $l$ modulo $4$ then
$$\chi(C([i^k,i^l]))=\chi(C([i^k]))\cdot
\chi(C([i^l]))=\frac{1}{16}.$$ It remain to compute $R(A)$ in the
two cases.
$$R([i^k,i^{k+1}])=N_{\Q(i)/\Q}(\det(P_{[i^k],[i^{k+1}]}))=
N_{\Q(i)/\Q}(i^{k+1}-i^k)=2.$$ And
$$R([i^k,i^{k+2}])=N_{\Q(i)/\Q}(\det(P_{[i^k],[i^{k+2}]}))=
N_{\Q(i)/\Q}(i^{k+2}-i^k)=4.$$

From theorem 2.8 part (b) we have that
 the only torsion element in $GL_2 (\Z[\xi_3])$ are these
elements whose Euler characteristic of their centralizer is not
zero are the ones whose eigenvalues are different and belong to
the set $$\{\pm1.\xi_3^{\pm 1}, \xi_6^{\pm 1}\}.$$

\ble{8.3} This is the computation of $|N_{\Q(i)/\Q}(R(T))|\cdot \chi(C(T))$ for
all torsion elements up to conjugation $T$ in $GL_2\Z[\xi_3]$ such
that $\chi(C(T))\neq 0$. Let $\xi_6=e^{\frac{2\pi i}{6}}$.\\

(a) $R([\xi_6^k,\xi_6^{k+1}])\chi(C([\xi_6^k,\xi_6^{k+1}]))=\frac{1}{36},$
for $k=0,1,\dots, 5,$\\

(b) $R([\xi_6^k,\xi_6^{k+2}])\chi(C([\xi_6^k,\xi_6^{k+2}]))=\frac{1}{12},$
for $k=0,1,\dots, 5,$\\

(c) $R([\xi_6^k,\xi_6^{k+3}])\chi(C([\xi_6^k,\xi_6^{k+3}]))=\frac{1}{9},$ for
$k=0,1,2.$
  \ele
\proof  If $k$ is not congruent to $l$ modulo $6$ then
$$\chi(C([\xi_6^k,\xi_6^l]))=\chi(C([\xi_6^k]))\cdot
\chi(C([\xi_6^l]))=\frac{1}{36}.$$ It remain to compute $R(A)$ in the
two cases.
$$R([\xi_6^k,\xi_6^{k+1}])=N_{\Q(\xi_6)/\Q}(\det(P_{[\xi_6^k],[\xi_6^{k+1}]}))=
N_{\Q(\xi_6)/\Q}(\xi_6^{k+1}-\xi_6^k)=1.$$ 
And 
$$R([\xi_6^k,\xi_6^{k+2}]) =
N_{\Q(\xi_6)/\Q}(\det(P_{[\xi_6^k],[\xi_6^{k+2}]}))=
N_{\Q(\xi_6)/\Q}(\xi_6^{k+2}-\xi_6^k)=3.$$ 
And
$$R([\xi_6^k,\xi_6^{k+3}])=N_{\Q(\xi_6)/\Q}(\det(P_{[\xi_6^k],[\xi_6^{k+3}]}))=
N_{\Q(\xi_6)/\Q}(\xi_6^{k+3}-\xi_6^k)=4.$$

\sectionnew{Homological Euler characteristic of $GL_m({\cal{O}}_K)$ 
with coefficients symmetric powers}
In this section we compute the homological Euler characteristic of 
$GL_n({\cal{O}}_K)$ for $K=\Q$, $Q(i)$ and $\Q(\xi_3)$. For other
imaginary quadratic fields the formula is given in theorem 2.13. 
And for all other fields the Euler characteristic vanishes 
(see theorem 0.1 part (c)).
\subsection{Homological Euler characteristic of $GL_m\Z$.}
\ble{5.4}
     For the symmetric powers of the standard representation $V_2$
of $GL_2\Z$, we have:\\
(a)   $\tr(I_2|S^nV_2)=n+1$ \\
\\
(b)   $\tr(-I_2|S^nV_2)=(-1)^{n+1}(n+1)$ \\
\\
(c)  
$\tr(\left[\begin{tabular}{rr} $1$ & $0$ \\
$0$ & $-1$ \\
\end{tabular}\right]|S^{2n+k}V_2)
=
\left\{\begin{tabular}{cc} $1$ & $k=0$ \\ $0$ & $k=1$ \\
\end{tabular}\right.
$\\
\\
\\
(d)  
$\tr(T_3|S^{3n+k}V_2)= \left\{\begin{tabular}{rr} $1$  & $k=0$ \\
$-1$ & $k=1$ \\ $0$  & $k=2$ \\
\end{tabular}\right.
$\\
\\
\\
(e)   $\tr(T_6|S^{6n+k}V_2)= \left\{\begin{tabular}{rr}

$1$  & $k=0$ \\

$1$  & $k=1$ \\

$0$  & $k=2$ \\

$-1$ & $k=3$ \\

$-1$ & $k=4$ \\

$0$  & $k=5$ \\

\end{tabular}\right.$\\
\\
\\
(f)  
  $\tr(T_4|S^{4n+k}V_2)=
\left\{\begin{tabular}{rr}

$1$  & $k=0$ \\

$0$  & $k=1$ \\

$-1$ & $k=2$ \\

$0$  & $k=3$ \\
\end{tabular}\right.
    $\\
\ele
  \proof If $\l_1$ and $\l_2$ are the two eigenvalues of $A$
acting on $V_2$ then $$\tr(A|S^nV_2)=\sum_{i=0}^n \l_1^i
\l_2^{n-i}.$$ Note that $\tr(A|S^nV_2)$ depends on the on the
eigenvalues of $A$ on $V_2$ not the conjugacy class that $A$
belongs to. Also, $$\tr(-A|S^nV_2)=(-1)^n tr(A|S^nV_2).$$ From
$$\tr(I_2|S^nV_2)=dim(S^nV_2)=n+1,$$ we obtain part (a) and (b). If
$\l_1=1$ and $\l_2=-1$ then $$\tr(A|S^nV_2)=1-1+1-\dots.$$ and we
have $n+1$ summands. Thus $\tr(A|S^nV_2)$ is $1$ if $n$ is even and
$0$ if $n$ is odd. This proves part (c). For part (d) the
eigenvalues of $T_3$ on $V_2$ are $\xi_3$ and $\xi_3^{-1}$. Thus,
$$\tr(T_3|S^nV_2)=\xi_3^n+\xi_3^{n-2}+\dots +\xi_3^{-n}.$$ Also the
sum of three successive summands is zero, and the total number of
summands is $n+1$. Therefore, if $n\equiv 0\mod 3$ then
$$\tr(T_3|S^nV_2)=\xi_3^n=1.$$ If $n\equiv 1\mod 3$ then
$$\tr(T_3|S^nV_2)=\xi_3^n+\xi_3^{n-2}=\xi_3+\xi_3^{-1}=-1.$$ And if
$n\equiv 2\mod 3$ then $$\tr(T_3|S^nV_2)=0.$$ This proves part (d).
It also proves part (e) because $-T_3$ is conjugate to $T_6$.
Similarly, for the $4$-torsion $T_4$ we have
$$\tr(T_4|S^nV_2)=i^n+i^{n-2}+\dots +i^{-n}.$$ Also, every two
consecutive summands add up to zero. Therefore if $n$ is odd the
trace is zero. If $n$ is even we have $$\tr(T_4|S^nV_2)=i^n,$$ which
is $1$ if $n\equiv\mbox{ } 0\mbox{ } \mod 4$ and $-1$ if 
$n\equiv 2 \mbox{ }\mod \mbox{ }4$. This
proves part (f).


\bth{5.5}
   Let $S^nV_2$ be the $n$-th symmetric power of the
standard representation of $GL_2$. Then\\

$\chi_h(GL_2\Z, S^{12n+k}V_2)= \left\{\begin{tabular}{cl} 
$-n+1$ & $k=0$ \\ 
$-n$    & $k=2$ \\ 
$-n$    & $k=4$ \\ 
$-n$    & $k=6$ \\ 
$-n$    & $k=8$ \\
$-n-1$  & $k=10$ \\
$0$     & $k=odd,$ \\
\end{tabular}\right.
     $\\
            and\\

$\chi_h(GL_2\Z, S^{12n+k}V_2\otimes \det)=
\left\{\begin{tabular}{cl}
$-n$   & $k=0$ \\
$-n-1$ & $k=2$ \\ 
$-n-1$ & $k=4$ \\ 
$-n-1$ & $k=6$ \\ 
$-n-1$ & $k=8$ \\ 
$-n-2$ & $k=10$ \\ 
$0$    & $k=odd.$ \\
   \end{tabular}\right.
    $
        \\
     \eth
\proof If $k$ is odd then $-I_2$ acts non-trivially on the
symmetric power. Therefore the cohomology groups of $GL_2Z$ with
coefficient that symmetric power vanishes. Consequently, the
homological Euler characteristic vanishes. If $k=0$, using lemma 5.1,
we obtain
$$\tr(I_2|S^{12n}V_2)=12n+1,$$ 
$$\tr(-I_2|S^{12n}V_2)=12n+1,$$
$$\tr([1,-1]|S^{12n}V_2)=1,$$ 
$$\tr(T_3|S^{12n}V_2)=1,$$
$$\tr(T_6|S^{12n}V_2)=1,$$  
$$\tr(T_4|S^{12n}V_2)=1.$$
Using the Euler characteristic of the centralizers of the torsion
elements in $GL_2\Z$ from proposition 3.1, we obtain\\
$$
\begin{tabular}{ccl}
$\chi_h(GL_2\Z,S^{12n}V_2)$ & $=$ 
 &$\chi(GL_2\Z)(12n+1)+\chi(GL_2\Z)(12n+1)+$\\
\\
 && $+2\cdot\chi(C([1,-1]))\cdot 1 +\chi(C(T_3))\cdot 1+$\\
\\
 && $+ \chi(C(T_6))\cdot 1 + \chi(C(T_4))\cdot 1=$\\
\\
 & $=$ & $-\frac{1}{24}(12n+1)-\frac{1}{24}(12n+1)+$\\ 
\\
 && $+2\cdot \frac{1}{4}\cdot 1+\frac{1}{6}\cdot 1+ 
       \frac{1}{6}\cdot 1+\frac{1}{4}\cdot 1=$\\
\\
 & $=$ & $-n+1.$
\end{tabular}
$$

For $k=2$ from lemma 5.1 we have 
$$\tr(I_2|S^{12n+2}V_2)=12n+3,$$
$$\tr(-I_2|S^{12n+2}V_2)=12n+3,$$ 
$$\tr([1,-1]|S^{12n+2}V_2)=1,$$
$$\tr(T_3|S^{12n+2}V_2)=0,$$
$$\tr(T_6|S^{12n+2}V_2)=0,$$  
$$\tr(T_4|S^{12n+2}V_2)=-1.$$
Using proposition 3.1 we obtain
$$
\begin{tabular}{ccl}
$\chi_h(GL_2\Z,S^{12n+2}V_2)$ & $=$ 
 & $-\frac{1}{24}(12n+3)-\frac{1}{24}(12n+3)+$\\
\\
 && $+2\cdot\frac{1}{4}\cdot 1+\frac{1}{6}\cdot 0+ \frac{1}{6}\cdot
0+\frac{1}{4}\cdot (-1)=$\\
\\
 & $=$ & $-n.$\\ 
\end{tabular}
$$

For $k=4$ from lemma 5.1 we have
$$\tr(I_2|S^{12n+4}V_2)=12n+5,$$ 
$$\tr(-I_2|S^{12n+4}V_2)=12n+5,$$
$$\tr([1,-1]|S^{12n+4}V_2)=1,$$ 
$$\tr(T_3|S^{12n+4}V_2)=-1,$$
$$\tr(T_6|S^{12n+4}V_2)=-1,$$  
$$\tr(T_4|S^{12n+4}V_2)=1,$$
Using proposition 3.1 we obtain
$$
\begin{tabular}{ccl}
$\chi_h(GL_2\Z,S^{12n+4}V_2)$ & $=$
 & $-\frac{1}{24}(12n+5)-\frac{1}{24}(12n+5)+$\\
\\
 && $+2\cdot\frac{1}{4}\cdot 1+\frac{1}{6}\cdot (-1)+ \frac{1}{6}\cdot
(-1)+\frac{1}{4}\cdot 1=$\\
\\
 & $=$ & $-n.$
\end{tabular}
$$

For $k=6$ from lemma 5.1 we have
$$\tr(I_2|S^{12n+6}V_2)=12n+7,$$ 
$$\tr(-I_2|S^{12n+6}V_2)=12n+7,$$
$$\tr([1,-1]|S^{12n+6}V_2)=1,$$
$$\tr(T_3|S^{12n+6}V_2)=1,$$
$$\tr(T_6|S^{12n+6}V_2)=1,$$
$$\tr(T_4|S^{12n+6}V_2)=-1.$$
Using proposition 3.1 we obtain
$$
\begin{tabular}{ccl}
$\chi_h(GL_2\Z,S^{12n+6}V_2)$ & $=$ 
 & $-\frac{1}{24}(12n+7)-\frac{1}{24}(12n+7)+$\\
\\
 && $+2\cdot\frac{1}{4}\cdot 1+\frac{1}{6}\cdot 1+ \frac{1}{6}\cdot
1+\frac{1}{4}\cdot (-1)=$\\
\\
 & $=$ & $-n.$\\
\end{tabular}
$$

For $k=8$ from lemma 5.1 we have
$$\tr(I_2|S^{12n+8}V_2)=12n+9,$$ 
$$\tr(-I_2|S^{12n+8}V_2)=12n+9,$$
$$\tr([1,-1]|S^{12n+8}V_2)=1,$$ 
$$\tr(T_3|S^{12n+8}V_2)=0,$$
$$\tr(T_6|S^{12n+8}V_2)=0,$$  
$$\tr(T_4|S^{12n+8}V_2)=1.$$
Using proposition 3.1 we obtain
$$
\begin{tabular}{ccl}
$\chi_h(GL_2\Z,S^{12n+8}V_2)$ & $=$
 & $-\frac{1}{24}(12n+9)-\frac{1}{24}(12n+9)+$\\
\\
 && $+2\cdot\frac{1}{4}\cdot 1+\frac{1}{6}\cdot 0+ \frac{1}{6}\cdot
0+\frac{1}{4}\cdot 1=$\\
\\
 &$=$ & $-n.$\\
\end{tabular}
$$

For $k=10$ from lemma 5.1 we have
$$\tr(I_2|S^{12n+10}V_2)=12n+11,$$ 
$$\tr(-I_2|S^{12n+10}V_2)=12n+11,$$
$$\tr([1,-1]|S^{12n+10}V_2)=1,$$ 
$$\tr(T_3|S^{12n+10}V_2)=-1,$$
$$\tr(T_6|S^{12n+2}V_2)=-1,$$  
$$\tr(T_4|S^{12n+2}V_2)=-1.$$
Using proposition 4.1 we obtain
$$
\begin{tabular}{ccl}
$\chi_h(GL_2\Z,S^{12n+10}V_2)$ & $=$
 & $-\frac{1}{24}(12n+11)-\frac{1}{24}(12n+11)+$\\
\\
 && $+2\cdot\frac{1}{4}\cdot 1+\frac{1}{6}\cdot (-1)+ \frac{1}{6}\cdot
(-1)+\frac{1}{4}\cdot (-1)=$\\
\\
 & $=$ & $-n-1.$
\end{tabular}
$$


Note that $$tr(-A|S^nV_3)=(-1)^n tr(A|S^nV_3).$$ Therefore we need
only formulas for half of the torsion elements. One other remark:
we need only to examine the trace for torsion elements whose Euler
characteristic of the centralizer is not zero. For $GL_3\Z$ these
elements are $I_3$ and $-I_3$.

 \ble{5.2}
    The trace of the torsion
elements in $GL_3\Z$ acting on the symmetric power of the standard
representation are given by:\\

(c) \mbox{ } $\tr([I_2,-1]|S^{2n+k}V_3)= \left\{\begin{tabular}{cc}
$n+1$ & $k=0$ \\ $n+1$ & $k=1$ \\
\end{tabular}\right.
$ \\

(d) \mbox{ } $\tr([-I_2,1]|S^{2n+k}V_3)= \left\{\begin{tabular}{rc}
$n+1$  & $k=0$ \\ $-n-1$ & $k=1$ \\
\end{tabular}\right.
$ \\

(e) \mbox{ } $\tr([T_3,1]|S^{3n+k}V_3)=
\left\{\begin{tabular}{cc} $1$ & $k=0$ \\ $0$ & $k=1$ \\ $0$ &
$k=2$ \\
\end{tabular}\right.
$\\

(f) \mbox{ } $\tr([T_6,1]|S^{6n+k}V_3)
=
\left\{\begin{tabular}{rc} $1$ & $k=0$ \\ $2$ & $k=1$ \\ $2$ &
$k=2$ \\ $1$ & $k=3$ \\ $0$ & $k=4$ \\ $0$ & $k=5$ \\
\end{tabular}\right.
$\\
\\

(g) \mbox{ } $\tr([T_3,-1]|S^{6n+k}V_3)
=
\left\{\begin{tabular}{rc} $1$  & $k=0$ \\ $-2$ & $k=1$ \\ $2$  &
$k=2$ \\ $-1$ & $k=3$ \\ $0$  & $k=4$ \\ $0$  & $k=5$ \\
\end{tabular}\right.
$\\
\\

(h) \mbox{ } $\tr([T_6,-1]|S^{6n+k}V_3)
=
\left\{\begin{tabular}{rc} $1$  & $k=0$ \\ $0$  & $k=1$ \\ $0$  &
$k=2$ \\ $-1$ & $k=3$ \\ $0$  & $k=4$ \\ $0$  & $k=5$ \\
\end{tabular}\right.
$\\
\\

(i)  \mbox{ } $\tr([T_4,1]|S^{4n+k}V_3)= \left\{\begin{tabular}{cc}
$1$  & $k=0$ \\ $1$  & $k=1$ \\ $0$  & $k=2$ \\ $0$  & $k=3$ \\
\end{tabular}\right.
$\\
\\

(j) \mbox{ } $\tr([T_4,-1]|S^{4n+k}V_3)= \left\{\begin{tabular}{cc}
$1$  & $k=0$ \\ $-1$ & $k=1$ \\ $0$  & $k=2$ \\ $0$  & $k=3$ \\
\end{tabular}\right.
$\\
   \ele
   We have omitted part (a) and (b) corresponding to $I_3$ and $-I_3$
because they give no contribution to the homological Euler
characteristic since $$\chi(C(\pm I_3))=\chi(GL_3\Z)=0.$$ 
\proof
The above formulas can be derived from the following fact: let
$A\e GL_k\Z$ and $B\e GL_l\Z$. Set $$f(n)=\tr([A,B]|S^n(V_k\oplus
V_l),$$ $$g(n)=\tr(A|S^nV_k)$$ and $$h(n)=\tr(B|S^nV_k).$$ Then
$$f(n)=\sum_{i=0}^n g(i)h(n-i)=(g*h)(n).$$ Then
$$
\begin{tabular}{ccl}
$\tr([I_2,-1|S^{2n}V_3)$ & $=$
 & $\tr(I_2|S^{2n}V_2)-\tr(I_2|S^{2n-1}V_2)+\dots=$\\
\\
 & $=$ & $2n-(2n-1)+(2n-2)-\dots=$\\
\\
 & $=$ & $n+1.$ 
\end{tabular}
$$
And
$$
\begin{tabular}{ccl}
$\tr([I_2,-1|S^{2n+1}V_3)$ & $=$
 & $\tr(I_2|S^{2n+1}V_2)-\tr(I_2|S^{2n}V_2)+\dots=$\\
\\
 & $=$ & $(2n+1)-2n+(2n-1)-\dots=$\\
\\
 & $=$ & $n+1.$ 
\end{tabular}
$$
This proves part (c) and (d).

For part (e) we use lemma 5.1 part(d) and\\
$$\tr([T_3,1]|S^nV_3)=\sum_{i=0}^n \tr([T_3]|S^iV_2)\tr([1]|S^{n-i}V_1).$$ 
Then
$$\tr([T_3,1]|S^{3n}V_3)=1-1+0+\dots 1 -1+0+1=1.$$ 
Similarly,
$$\tr([T_3,1]|S^{3n+1}V_3)=1-1+0+\dots 1 -1+0+1-1=0.$$
And
$$\tr([T_3,1]|S^{3n+2}V_3)=1-1+0+\dots 1 -1+0=0.$$ 
This proves part (e) and (h).

  For part (f) we use lemma 5.1 part (e) and
$$\tr([T_6,1]|S^nV_3)=\sum_{i=0}^n
\tr([T_6]|S^iV_2)\tr([1]|S^{n-i}V_1).$$
Then
$$\tr([T_6,1]|S^{6n}V_3)=1+1+0-1-1+0\dots 1+1+0-1-1+1=1.$$
    Similarly,
$$
\begin{tabular}{l}
   $\tr([T_6,1]|S^{6n+1}V_3)=1+1+0-1-1+0\dots
1+1+0-1-1+1+1=2,$\\
\\
 $\tr([T_6,1]|S^{6n+2}V_3)=1+1+0-1-1+0\dots
1+1+0-1-1+1+1+0=2,$\\
\\
    $\tr([T_6,1]|S^{6n+3}V_3)=1+1+0-1-1+0\dots
1+1+0-1-1+1+1+0-1=1,$\\
\\  
  $\tr([T_6,1]|S^{6n+4}V_3)=1+1+0-1-1+0\dots
1+1+0-1-1+1+1+0-1-1=0,$\\
\\  
   $\tr([T_6,1]|S^{6n+4}V_3)=1+1+0-1-1+0\dots
1+1+0-1-1+0=0.$
\end{tabular}
$$
This proves part (f) and part (g).

For part (e) we use lemma 5.1 part (f) and
$$\tr([T_4,1]|S^nV_3)=\sum_{i=0}^n
\tr([T_4]|S^iV_2)\tr([1]|S^{n-i}V_1).$$
     Then
$$\tr([T_4,1]|S^{4n}V_3)=1+0-1+0+\dots 1+0 -1+0+1=1.$$
    Similarly,
$$
\begin{tabular}{l}
$\tr([T_4,1]|S^{4n+1}V_3)=1+0-1+0+\dots 1+0 -1+0+1+0=1,$\\
\\
$\tr([T_4,1]|S^{4n+1}V_3)=1+0-1+0+\dots 1+0 -1+0+1+0-1=0,$\\
\\
$\tr([T_4,1]|S^{4n+3}V_3)=1+0-1+0+\dots 1+0 -1+0=0.$
\end{tabular}
$$
     This proves part (i) and (j).

\bth{5.5}
    $\chi_h(GL_3\Z,S^nV_3)=\chi_h(GL_2\Z,S^nV_2),$
   \eth
      \proof We are going to use theorem 2.10, lemma 4.2 and lemma 5.3.
Also, for $$\chi_h(GL_3\Z,S^nV_3)$$ we only need to consider even
$n$ because if $n$ is odd then $-I_3$ acts non-trivially on the
representation. Another observation that we need to make is that
if $A$ is a torsion element in $GL_3\Z$ then the torsion element
$-A$ is non-conjugate to $A$. However
$$\tr(-A|S^{2n}V_3)=\tr(A|S^{2n}V_3).$$ Thus we only need to consider
half of the torsion elements, and just multiply by $2$ in order to
incorporate the other half of the torsion elements. For the
representation $S^{12n}V_3$ we have
 $$\chi_h(GL_3\Z,S^nV_3)=2\sum
R(A)\chi(C(A))\tr(A^{-1}|S^nV_3),$$ where the sum is taken over $A$
varying in the set
   $$\{[I_2,-1],[T_3,1],[T_6,1],[T_4,1]\}.$$ Then
$$
\begin{tabular}{lcl}
$\chi_h(GL_3\Z,S^{12n}V_3)$ & $=$
 & $2(4(-\frac{1}{48})(6n+1)+ 3\cdot
\frac{1}{12}\cdot 1 +1\cdot\frac{1}{12}\cdot 1+ 2\cdot
\frac{1}{8}\cdot 1)=$\\
\\
 & $=$ & $-n+1,$
\end{tabular}
$$
$$
\begin{tabular}{lcl}
\\
$\chi_h(GL_3\Z,S^{12n+2}V_3)$ & $=$
 & $2(4(-\frac{1}{48})(6n+2)+ 3\cdot
\frac{1}{12}\cdot 0 +1\cdot\frac{1}{12}\cdot 2+ 2\cdot
\frac{1}{8}\cdot 0)=$\\
\\
 & $=$ & $-n,$
\end{tabular}
$$
$$
\begin{tabular}{lcl}
$\chi_h(GL_3\Z,S^{12n+4}V_3)$ & $=$
 & $2(4(-\frac{1}{48})(6n+3)+ 3\cdot
\frac{1}{12}\cdot 0 +1\cdot\frac{1}{12}\cdot 0+ 2\cdot
\frac{1}{8}\cdot 1)=$\\
\\
 & $=$ & $-n,$
\end{tabular}
$$
$$
\begin{tabular}{lcl}
$\chi_h(GL_3\Z,S^{12n+6}V_3)$ & $=$
 & $2(4(-\frac{1}{48})(6n+4)+ 3\cdot
\frac{1}{12}\cdot 1 +1\cdot\frac{1}{12}\cdot 1+ 2\cdot
\frac{1}{8}\cdot 0)=$\\
\\
 & $=$ & $-n,$
\end{tabular}
$$
$$
\begin{tabular}{lcl}
$\chi_h(GL_3\Z,S^{12n+8}V_3)$ & $=$
 & $2(4(-\frac{1}{48})(6n+5)+ 3\cdot
\frac{1}{12}\cdot 0 +1\cdot\frac{1}{12}\cdot 2+ 2\cdot
\frac{1}{8}\cdot 1)=$\\
\\
 & $=$ & $-n,$
\end{tabular}
$$
$$
\begin{tabular}{lcl}
$\chi_h(GL_3\Z,S^{12n+10}V_3)$ & $=$
 & $2(4(-\frac{1}{48})(6n+6)+ 3\cdot
\frac{1}{12}\cdot 0 +1\cdot\frac{1}{12}\cdot 0+ 2\cdot
\frac{1}{8}\cdot 0)=$\\
\\
 & $=$ & $-n-1.$
\end{tabular}
$$
     Thus,
$$\chi_h(GL_3\Z,S^nV_3)=\chi_h(GL_2\Z,S^nV_2).$$


\ble{5.7}
    The trace of the torsion elements of $GL_4\Z$
acting on the symmetric powers of the standard representation are
given by\\

(a) \mbox{ } $\tr([I_2,-I_2]|S^{2n+k}V_4)=
\left\{\begin{tabular}{cc} 
$n+1$ & $k=0$ \\ 
$0$   & $k=1$ \\
\end{tabular}\right.
$ \\
\\

(b) \mbox{ }
 $\tr([I_2,T_3]|S^{3n+k}V_4)=
\left\{\begin{tabular}{cc} 
$n+1$ & $k=0$ \\ 
$n+1$ & $k=1$ \\ 
$n+1$ & $k=2$ \\
\end{tabular}\right.
$\\
\\

(c) \mbox{ } $\tr([I_2,T_6]|S^{6n+k}V_4)
=
\left\{\begin{tabular}{rc} 
$6n+1$ & $k=0$ \\ 
$6n+3$ & $k=1$ \\
$6n+5$ & $k=2$ \\ 
$6n+6$ & $k=3$ \\ 
$6n+6$ & $k=4$ \\ 
$6n+6$ & $k=5$ \\
\end{tabular}\right.
$\\
\\

(d)  \mbox{ } $\tr([I_2,T_4]|S^{4n+k}V_4)=
\left\{\begin{tabular}{cl} 
$2n+1$  & $k=0$ \\ 
$2n+2$  & $k=1,2,3$\\
\end{tabular}\right.
$\\
\\

(e) \mbox{ } $\tr([-I_2,T_3]|S^{6n+k}V_4)
=
\left\{\begin{tabular}{rc}
  $6n+1$  & $k=0$ \\
$-(6n+3)$ & $k=1$ \\
  $6n+5$  & $k=2$ \\
$-(6n+6)$ & $k=3$ \\
  $6n+6$  & $k=4$ \\
$-(6n+6)$ & $k=5$ \\
\end{tabular}\right.
$\\
\\
(f) \mbox{ } $\tr([-I_2,T_6]|S^{6n+k}V_4)
=
\left\{\begin{tabular}{rc} 
$2n+1$    & $k=0$ \\ 
$-(2n+1)$ & $k=1$ \\
$2n+1$    & $k=2$ \\ 
$-(2n+2)$ & $k=3$ \\ 
$2n+2$    & $k=4$ \\
$-(2n+2)$ & $k=5$ \\
\end{tabular}\right.
$\\
\\

(g)  \mbox{ } $\tr([-I_2,T_4]|S^{4n+k}V_4)=
\left\{\begin{tabular}{cl}
  $2n+1$   & $k=0$ \\
$-(2n+2)$  & $k=1$ \\
  $2n+2$   & $k=2$ \\
$-(2n+2)$  & $k=3$ \\
\end{tabular}\right.
$\\
\\

(h) \mbox{ }
 $\tr([1,-1,T_3]|S^{6n+k}V_4)=
\left\{\begin{tabular}{cc} 
$1$  & $k=0$ \\ 
$-1$ & $k=1$ \\ 
$1$  & $k=2$ \\ 
$0$  & $k=3$ \\ 
$0$  & $k=4$ \\ 
$0$  & $k=5$ \\
\end{tabular}\right.
$\\
\\

(i) \mbox{ } $\tr([1,-1,T_6]|S^{6n+k}V_4)
=
\left\{\begin{tabular}{rc} $1$ & $k=0$ \\
   $1$ & $k=1$ \\
   $1$ & $k=2$ \\
   $0$ & $k=3$ \\
   $0$ & $k=4$ \\
   $0$ & $k=5$ \\
\end{tabular}\right.
$\\
\\

(j)  \mbox{ }
   $\tr([1,-1,T_4]|S^{4n+k}V_4)=
\left\{\begin{tabular}{cl}
   $1$ & $k=0$ \\
   $0$ & $k=1$ \\
   $0$ & $k=2$ \\
   $0$ & $k=3$ \\
\end{tabular}\right.
$\\
\\

(k) \mbox{ }
$\tr([T_3,T_6]|S^{6n+k}V_4)
=
\left\{\begin{tabular}{rc}
   $1$ & $k=0$ \\
   $0$ & $k=1$ \\
   $-1$ & $k=2$ \\
   $0$ & $k=3$ \\
   $0$ & $k=4$ \\
   $0$ & $k=5$ \\
\end{tabular}\right.
$\\
\\

(l) \mbox{ }
    $\tr([T_3,T_4]|S^{12n+k}V_4)=
   \left\{\begin{tabular}{rlrl}
     $1$  & $k=0,$ & $1$  & $k=6$ \\
     $-1$ & $k=1,$ & $1$  & $k=7$ \\
     $-1$ & $k=2,$ & $-1$ & $k=8$ \\
     $2$  & $k=3,$ & $0$  & $k=9$ \\
     $0$  & $k=4,$ & $0$  & $k=10$ \\
     $-2$ & $k=5,$ & $0$  & $k=11$ \\
    \end{tabular}\right.
 $\\
\\

(m) \mbox{ }
  $\tr([T_6,T_4]|S^{12n+k}V_4)=
\left\{\begin{tabular}{rlrl}
   $1$  & $k=0,$ & $1$  & $k=6$ \\
   $1$  & $k=1,$ & $-1$ & $k=7$ \\
   $-1$ & $k=2,$ & $-1$ & $k=8$ \\
   $-2$ & $k=3,$ & $0$  & $k=9$ \\
   $0$  & $k=4,$ & $0$  & $k=10$ \\
   $2$  & $k=5,$ & $0$  & $k=11$ \\
 \end{tabular}\right.
$\\
  \ele
{\bf Remark:} In the above lemma we have listed only
 the torsion elements whose Euler characteristic of their
centralizer is non zero. That is only these torsion elements will
have a contribution toward the homological Euler characteristic of
$GL_4\Z$ or of arithmetic subgroup.\\ 
\proof We are going to use
that $$\tr(-A|S^nV_4)=(-1)^n\tr(A|S^nV_4).$$ For part (a) we use lemma
6.3 part(c) and
  $$\tr([I_2,-I_2]|S^nV_4)=\sum_{i=0}^n
\tr([-I_2,1]|S^iV_3)\tr([1]|S^{n-i}V_1).$$
   Then
$$
\begin{tabular}{lcl}
$\tr([I_2,-I_2]|S^{2n}V_4)$ & $=$ & $n+1-n+n+\dots -1+1=$\\
\\
 & $=$ & $n+1.$\\
\\
$\tr([I_2,-I_2]|S^{2n+1}V_4)$ & $=$ & $-(n+1)+(n+1)-n+n+\dots -1+1=$\\
\\
 & $=$ & $0.$\\
\end{tabular}
$$
    For
part (b) we use lemma 5.3 part(e) and
$$\tr([I_2,T_3]|S^nV_4)=\sum_{i=0}^n
\tr([T_3,1]|S^iV_3)\tr([1]|S^{n-i}V_1).$$
   Then
$$
\begin{tabular}{lcl}
$\tr([I_2,T_3]|S^{3n}V_4)$ & $=$ & $1+0+0\dots 1+0+0+1=$\\
\\
 & $=$ & $n+1$
\end{tabular}
$$
$$
\begin{tabular}{lcl}
$\tr([I_2,T_3]|S^{3n+1}V_4)$ & $=$ & $1+0+0\dots 1+0+0+1+0=$\\
\\
 & $=$ & $n+1$
\end{tabular}
$$
$$
\begin{tabular}{lcl}
$\tr([I_2,T_3]|S^{3n+2}V_4)$ & $=$ & $1+0+0\dots 1+0+0+1+0+0=$\\
\\
 & $=$ & $n+1$
\end{tabular}
$$
 This proves part (b)and part (f).

For part (d) we use lemma 5.3
part(i) and\\
    $$\tr([I_2,T_4]|S^nV_4)=\sum_{i=0}^n
\tr([T_4,1]|S^iV_3)\tr([1]|S^{n-i}V_1).$$
   Then
$$
\begin{tabular}{lcl}
$\tr([I_2,T_3]|S^{4n}V_4)$ & $=$ & $1+1+0+0\dots 1+1+0+0+1=$\\
\\
 & $=$ & $2n+1,$
\end{tabular}
$$
$$
\begin{tabular}{lcl}
$\tr([I_2,T_3]|S^{4n+1}V_4)$ & $=$ & $1+1+0+0\dots 1+1+0+0+1+1=$\\
\\
 & $=$ & $2n+2,$
\end{tabular}
$$
$$
\begin{tabular}{lcl}
$\tr([I_2,T_3]|S^{4n+2}V_4)$ & $=$ & $1+1+0+0\dots 1+1+0+0+1+1+0=$\\
\\
 & $=$ & $2n+2,$
\end{tabular}
$$
$$
\begin{tabular}{lcl}
$\tr([I_2,T_3]|S^{4n+3}V_4)$&$=$&$1+1+0+0\dots 1+1+0+0+1+1+0+1=$\\
\\
 & $=$ & $2n+2.$
\end{tabular}
$$
This proves part (d) and (g).

For part (f) we use lemma 5.3
part(f) and\\
    $$\tr([I_2,T_6]|S^nV_4)=\sum_{i=0}^n
\tr([T_6,1]|S^iV_3)\tr([1]|S^{n-i}V_1).$$
   Then
$$
\begin{tabular}{lcl}
$\tr([I_2,T_6]|S^{6n}V_4)$ & $=$ & $1+2+2+1+0+0\dots +1+2+2+1+0+0+1=$\\
\\
 & $=$ & $6n+1,$
\end{tabular}
$$
$$
\begin{tabular}{lcl}
$\tr([I_2,T_6]|S^{6n+1}V_4)$ & $=$ & $1+2+2+1+0+0\dots+1+2+2+1+0+0+1+2=$\\
\\
 & $=$ & $6n+3,$
\end{tabular}
$$
$$
\begin{tabular}{lcl}
$\tr([I_2,T_6]|S^{6n+2}V_4)$ & $=$ & $1+2+2+1+0+0\dots +1+2+2=$\\
\\
 & $=$ & $6n+5,$
\end{tabular}
$$
$$
\begin{tabular}{lcl}
$\tr([I_2,T_6]|S^{6n+3}V_4)$ & $=$ & $1+2+2+1+0+0\dots +1+2+2+1=$\\
\\
 & $=$ & $6n+6,$
\end{tabular}
$$
$$
\begin{tabular}{lcl}
$\tr([I_2,T_6]|S^{6n+4}V_4)$ & $=$ & $1+2+2+1+0+0\dots +1+2+2+1+0=$\\
\\
 & $=$ & $6n+6,$
\end{tabular}
$$
$$
\begin{tabular}{lcl}
$\tr([I_2,T_6]|S^{6n+5}V_4)$ & $=$ & $1+2+2+1+0+0\dots +1+2+2+1+0+0=$\\
\\
 & $=$ & $6n+6$
\end{tabular}
$$
This proves parts (c) and (e).

For part (h) we use lemma 5.3
part(g) and
 $$\tr([1,-1,T_3]|S^nV_4)=\sum_{i=0}^n
\tr([T_3,-1]|S^iV_3)\tr([1]|S^{n-i}V_1).$$
   Then
$$
\begin{tabular}{lcl}
$\tr([1,-1,T_3]|S^{6n}V_4)$ & $=$ & $1-2+2-1+0+0\dots 1-2+2-1+0+0+1=$\\
\\
 & $=$ & $1,$
\end{tabular}
$$
$$
\begin{tabular}{lcl}
$\tr([1,-1,T_3]|S^{6n+1}V_4)$ & $=$ & $1-2+2-1+0+0\dots 1-2=$\\
\\
 & $=$ & $-1,$
\end{tabular}
$$
$$
\begin{tabular}{lcl}
$\tr([1,-1,T_3]|S^{6n+2}V_4)$ & $=$ & $1-2+2-1+0+0\dots 1-2+2=$\\
\\
 & $=$ & $1,$
\end{tabular}
$$
$$
\begin{tabular}{lcl}
$\tr([1,-1,T_3]|S^{6n+3}V_4)$ & $=$ & $1-2+2-1+0+0\dots 1-2+2-1=$\\
\\
 & $=$ & $0,$
\end{tabular}
$$
$$
\begin{tabular}{lcl}
$\tr([1,-1,T_3]|S^{6n+4}V_4)$ & $=$ & $1-2+2-1+0+0\dots 1-2+2-1+0=$\\
\\
 & $=$ & $0,$
\end{tabular}
$$
$$
\begin{tabular}{lcl}
$\tr([1,-1,T_3]|S^{6n+5}V_4)$ & $=$ & $1-2+2-1+0+0\dots 1-2+2-1+0+0=$\\
\\
 & $=$ & $0.$
\end{tabular}
$$
This proves parts (h) and (i).

For part (j) we use lemma 5.3 part(j)
and
   $$\tr([1,-1,T_4]|S^nV_4)=\sum_{i=0}^n
\tr([T_4,-1]|S^iV_3)\tr([1]|S^{n-i}V_1).$$
    Then
$$
\begin{tabular}{lcl}
$\tr([1,-1,T_4]|S^{4n}V_4)$ & $=$ & $1-1+0+0\dots 1-1+0+0+1=$\\
\\
 & $=$ & $1,$
\end{tabular}
$$
$$
\begin{tabular}{lcl}
$\tr([1,-1,T_4]|S^{4n+1}V_4)$ & $=$ & $1-1+0+0\dots 1-1+0+0+1-1=$\\
\\
 & $=$ & $0,$
\end{tabular}
$$
$$
\begin{tabular}{lcl}
$\tr([1,-1,T_4]|S^{4n+2}V_4)$ & $=$ & $1-1+0+0\dots 1-1+0+0+1-1+0=$\\
\\
 & $=$ & $0,$
\end{tabular}
$$
$$
\begin{tabular}{lcl}
$\tr([1,-1,T_4]|S^{4n+3}V_4)$ & $=$ & $1-1+0+0\dots 1-1+0+0+1-1+0+0=$\\
\\
 & $=$ & $0.$
\end{tabular}
$$

For part (k) we use lemma 5.1 parts (d) and (e), and
   $$\tr([T_3,T_6]|S^nV_4)=\sum_{i=0}^n
\tr([T_3]|S^iV_2)\tr([T_6]|S^{n-i}V_2).$$
    Then
$$
\begin{tabular}{lcl}
$\tr([T_3,T_6]|S^{6n}V_4)$ & $=$ 
& $n(1\cdot1 -1\cdot0 +0\cdot(-1)+1\cdot(-1) -1\cdot0+0\cdot 1)+$\\
\\
&&$+1\cdot1=$\\
\\
 & $=$ & $1,$
\end{tabular}
$$
$$
\begin{tabular}{lcl}
$\tr([T_3,T_6]|S^{6n+1}V_4)$ & $=$ 
& $n(1\cdot1 -1\cdot1 +0\cdot0+1\cdot(-1) -1\cdot(-1)+0\cdot 0)+$\\
\\
&&$+1\cdot1-1\cdot1=$\\
\\
 & $=$ & $0,$
\end{tabular}
$$
$$
\begin{tabular}{lcl}
$\tr([T_3,T_6]|S^{6n+2}V_4)$ & $=$ 
& $n(1\cdot0 -1\cdot1 +0\cdot1+1\cdot0 -1\cdot(-1)+0\cdot(-1))+$\\
\\
&&$+1\cdot0-1\cdot1+0\cdot1=$\\
\\
 & $=$ & $-1,$
\end{tabular}
$$
$$
\begin{tabular}{lcl}
$\tr([T_3,T_6]|S^{6n+3}V_4)$ & $=$ 
& $n(1\cdot(-1) -1\cdot0 +0\cdot1+1\cdot1 -1\cdot0+0\cdot(-1))+$\\
\\
&&$+1\cdot(-1)-1\cdot0+0\cdot1+1\cdot1=$\\
\\
 & $=$ & $0,$
\end{tabular}
$$
$$
\begin{tabular}{lcl}
$\tr([T_3,T_6]|S^{6n+4}V_4)$ & $=$ 
& $n(1\cdot(-1) -1\cdot(-1) +0\cdot0+1\cdot1 -1\cdot1+0\cdot0)+$\\
\\
&&$+1\cdot(-1)-1\cdot(-1)+0\cdot0+1\cdot1-1\cdot1=$\\
\\
 & $=$ & $0,$
\end{tabular}
$$
$$
\begin{tabular}{lcl}
$\tr([T_3,T_6]|S^{6n+5}V_4)$ & $=$ 
& $(n+1)(1\cdot0 -1\cdot(-1) +0\cdot(-1)+1\cdot0 -1\cdot1+0\cdot 1)+$\\
\\
 & $=$ & $0.$
\end{tabular}
$$

For parts (l) and (m) we use lemma 5.1 parts (d) and (f), and
   $$\tr([T_3,T_4]|S^nV_4)=\sum_{i=0}^n
\tr([T_3]|S^iV_2)\tr([T_4]|S^{n-i}V_2).$$
    Then
$$
\begin{tabular}{lcl}
$\tr([T_3,T_4]|S^{12n}V_4)$ & $=$ 
& $n(1\cdot1 -1\cdot0 +0\cdot(-1)+1\cdot0 -1\cdot1+0\cdot0+$\\
\\
&& $+1\cdot(-1) -1\cdot0 +0\cdot1+1\cdot0 -1\cdot(-1)+0\cdot0 )+$\\
\\

&&$+1\cdot1=$\\
\\
 & $=$ & $1,$
\end{tabular}
$$
$$
\begin{tabular}{lcl}
$\tr([T_3,T_4]|S^{12n+1}V_4)$ & $=$ 
& $n(1\cdot0 -1\cdot1 +0\cdot0+1\cdot(-1) -1\cdot0+0\cdot1+$\\
\\
&& $+1\cdot0 -1\cdot(-1) +0\cdot0+1\cdot1 -1\cdot0+0\cdot(-1))+$\\
\\

&&$+1\cdot0-1\cdot1=$\\
\\
 & $=$ & $-1,$
\end{tabular}
$$
$$
\begin{tabular}{lcl}
$\tr([T_3,T_4]|S^{12n+2}V_4)$ & $=$ 
& $n(1\cdot(-1) -1\cdot0 +0\cdot1+1\cdot0 -1\cdot(-1)+0\cdot0+$\\
\\
&& $+1\cdot1 -1\cdot0 +0\cdot(-1)+1\cdot0 -1\cdot1+0\cdot0 )+$\\
\\

&&$+1\cdot(-1)-1\cdot0+0\cdot1=$\\
\\
 & $=$ & $-1,$
\end{tabular}
$$
$$
\begin{tabular}{lcl}
$\tr([T_3,T_4]|S^{12n+3}V_4)$ & $=$ 
& $n(1\cdot0 -1\cdot(-1) +0\cdot0+1\cdot1 -1\cdot0+0\cdot(-1)+$\\
\\
&& $+1\cdot0 -1\cdot1 +0\cdot0+1\cdot(-1) -1\cdot0+0\cdot1)+$\\
\\

&&$+1\cdot0-1\cdot(-1)+0\cdot0+1\cdot1=$\\
\\
 & $=$ & $2,$
\end{tabular}
$$
$$
\begin{tabular}{lcl}
$\tr([T_3,T_4]|S^{12n+4}V_4)$ & $=$ 
& $n(1\cdot1 -1\cdot0 +0\cdot(-1)+1\cdot0 -1\cdot1+0\cdot0+$\\
\\
&& $+1\cdot(-1) -1\cdot0 +0\cdot1+1\cdot0 -1\cdot(-1)+0\cdot0 )+$\\
\\

&&$+1\cdot1-1\cdot0 +0\cdot(-1)+1\cdot0 -1\cdot1=$\\
\\
 & $=$ & $0,$
\end{tabular}
$$
$$
\begin{tabular}{lcl}
$\tr([T_3,T_4]|S^{12n+5}V_4)$ & $=$ 
& $n(1\cdot0 -1\cdot1 +0\cdot0+1\cdot(-1) -1\cdot0+0\cdot1+$\\
\\
&& $+1\cdot0 -1\cdot(-1) +0\cdot0+1\cdot1 -1\cdot0+0\cdot(-1))+$\\
\\

&&$+1\cdot0-1\cdot1 +0\cdot0+1\cdot(-1) -1\cdot0+0\cdot1 =$\\
\\
 & $=$ & $-2,$
\end{tabular}
$$
$$
\begin{tabular}{lcl}
$\tr([T_3,T_4]|S^{12n+6}V_4)$ & $=$ 
& $n(1\cdot(-1) -1\cdot0 +0\cdot1+1\cdot0 -1\cdot(-1)+0\cdot0+$\\
\\
&& $+1\cdot1 -1\cdot0 +0\cdot(-1)+1\cdot0 -1\cdot1+0\cdot0 )+$\\
\\

&&$+1\cdot(-1)-1\cdot0+0\cdot1+1\cdot0 -1\cdot(-1)+0\cdot0+$\\
\\
&&  $+1\cdot1 =$\\
\\
 & $=$ & $1,$
\end{tabular}
$$
$$
\begin{tabular}{lcl}
$\tr([T_3,T_4]|S^{12n+7}V_4)$ & $=$ 
& $n(1\cdot0 -1\cdot(-1) +0\cdot0+1\cdot1 -1\cdot0+0\cdot(-1)+$\\
\\
&& $+1\cdot0 -1\cdot1 +0\cdot0+1\cdot(-1) -1\cdot0+0\cdot1)+$\\
\\

&&$+1\cdot0-1\cdot(-1)+0\cdot0+1\cdot1  -1\cdot0+0\cdot(-1)+$\\
\\
&&$+1\cdot0 -1\cdot1=$\\
\\
 & $=$ & $1,$
\end{tabular}
$$
$$
\begin{tabular}{lcl}
$\tr([T_3,T_4]|S^{12n+8}V_4)$ & $=$ 
& $n(1\cdot1 -1\cdot0 +0\cdot(-1)+1\cdot0 -1\cdot1+0\cdot0+$\\
\\
&& $+1\cdot(-1) -1\cdot0 +0\cdot1+1\cdot0 -1\cdot(-1)+0\cdot0 )+$\\
\\

&&$+1\cdot1 -1\cdot0 +0\cdot(-1)+1\cdot0 -1\cdot1+0\cdot0+$\\
\\
&&$+1\cdot(-1) -1\cdot0 +0\cdot1=$\\
\\
 & $=$ & $-1,$
\end{tabular}
$$
$$
\begin{tabular}{lcl}
$\tr([T_3,T_4]|S^{12n+9}V_4)$ & $=$ 
& $n(1\cdot0 -1\cdot1 +0\cdot0+1\cdot(-1) -1\cdot0+0\cdot1+$\\
\\
&& $+1\cdot0 -1\cdot(-1) +0\cdot0+1\cdot1 -1\cdot0+0\cdot(-1))+$\\
\\
&&$+1\cdot0-1\cdot1 +0\cdot0+1\cdot(-1) -1\cdot0+0\cdot1+$\\
\\
&& $+1\cdot0 -1\cdot(-1) +0\cdot0+1\cdot1=$\\
\\
 & $=$ & $0,$
\end{tabular}
$$
$$
\begin{tabular}{lcl}
$\tr([T_3,T_4]|S^{12n+10}V_4)$ & $=$ 
& $n(1\cdot(-1) -1\cdot0 +0\cdot1+1\cdot0 -1\cdot(-1)+0\cdot0+$\\
\\
&& $+1\cdot1 -1\cdot0 +0\cdot(-1)+1\cdot0 -1\cdot1+0\cdot0 )+$\\
\\

&&$+1\cdot(-1)-1\cdot0+0\cdot1+1\cdot0 -1\cdot(-1)+0\cdot0+$\\
\\
&& $+1\cdot1 -1\cdot0 +0\cdot(-1)+1\cdot0 -1\cdot1=$\\
\\
 & $=$ & $0,$
\end{tabular}
$$
$$
\begin{tabular}{lcl}
$\tr([T_3,T_4]|S^{12n+11}V_4)$ & $=$ 
& $(n+1)(1\cdot0 -1\cdot(-1) +0\cdot0+1\cdot1 -1\cdot0+0\cdot(-1)+$\\
\\
&& $+1\cdot0 -1\cdot1 +0\cdot0+1\cdot(-1) -1\cdot0+0\cdot1)=$\\
\\
 & $=$ & $0.$
\end{tabular}
$$


\bth{6.6}
   $\chi_h(GL_4\Z,S^nV_4)=\chi_h(GL_2\Z,S^nV_2)$
    \eth
\proof If we consider odd symmetric powers then $-I_4$ acts non-trivially.
Therefore the cohomology $H^i(GL_4\Z,S^{2n+1}V_4)$ will vanish.
And so will the homological Euler characteristic. In order to 
commute the homological Euler characteristic with coefficients
in the even symmetric powers we use theorem 2.10, lemma 5.5 and lemma 4.3.
Another observation is that for torsion matrices $A$ in $GL_m \Z$
we have that $A^{-1}$ and $A$ have the same eigenvalues. Thus, we can sum 
traces of the type $\tr(A|V)$. Thus the formula from theorem 2.10
can be simplified to 
$$\chi_h(GL_4\Z,S^nV_4)=\sum R(A)\chi(C(A))\tr(A|V).$$
Then
$$
\begin{tabular}{rrl}
$\chi_h(GL_4\Z,S^{12n}V_4)$ &=&
 $\frac{1}{36}(6n+1)-\frac{1}{16}(4n+1)
      -\frac{1}{144}(12n+1)-\frac{1}{24}(6n+1) -$\\
\\
&& $-\frac{1}{144}(12n+1)-\frac{1}{16}(4n+1)-\frac{1}{24}(6n+1) +$\\
\\
&& $+\frac{1}{4}\cdot1+\frac{1}{4}\cdot1+\frac{1}{2}\cdot1+
   \frac{1}{9}\cdot1+\frac{1}{24}\cdot1+\frac{1}{24}\cdot1=$\\
\\
&=&$1-n,$
\end{tabular}
$$
$$
\begin{tabular}{rrl}
$\chi_h(GL_4\Z,S^{12n+2}V_4)$ &=&
 $\frac{1}{36}(6n+2)-\frac{1}{16}(4n+1)
      -\frac{1}{144}(12n+5)-\frac{1}{24}(6n+2) -$\\
\\
&& $-\frac{1}{144}(12n+5)-\frac{1}{16}(4n+1)-\frac{1}{24}(6n+2) +$\\
\\
&& $+\frac{1}{4}\cdot1+\frac{1}{4}\cdot1+\frac{1}{2}\cdot0+
   \frac{1}{9}\cdot(-1)+\frac{1}{24}\cdot(-1)+\frac{1}{24}\cdot(-1)=$\\
\\
&=&$-n,$
\end{tabular}
$$
$$
\begin{tabular}{rrl}
$\chi_h(GL_4\Z,S^{12n+4}V_4)$ &=&
 $\frac{1}{36}(6n+3)-\frac{1}{16}(4n+2)
      -\frac{1}{144}(12n+6)-\frac{1}{24}(6n+3) -$\\
\\
&& $-\frac{1}{144}(12n+6)-\frac{1}{16}(4n+2)-\frac{1}{24}(6n+3) +$\\
\\
&& $+\frac{1}{4}\cdot0+\frac{1}{4}\cdot0+\frac{1}{2}\cdot1+
   \frac{1}{9}\cdot0+\frac{1}{24}\cdot0+\frac{1}{24}\cdot0=$\\
\\
&=&$-n,$
\end{tabular}
$$
$$
\begin{tabular}{rrl}
$\chi_h(GL_4\Z,S^{12n+6}V_4)$ &=&
 $\frac{1}{36}(6n+4)-\frac{1}{16}(4n+3)
      -\frac{1}{144}(12n+7)-\frac{1}{24}(6n+4) -$\\
\\
&& $-\frac{1}{144}(12n+7)-\frac{1}{16}(4n+3)-\frac{1}{24}(6n+4) +$\\
\\
&& $+\frac{1}{4}\cdot1+\frac{1}{4}\cdot1+\frac{1}{2}\cdot0+
   \frac{1}{9}\cdot1+\frac{1}{24}\cdot1+\frac{1}{24}\cdot1=$\\
\\
&=&$-n,$
\end{tabular}
$$
$$
\begin{tabular}{rrl}
$\chi_h(GL_4\Z,S^{12n+8}V_4)$ &=&
 $\frac{1}{36}(6n+5)-\frac{1}{16}(4n+3)
      -\frac{1}{144}(12n+11)-\frac{1}{24}(6n+5) -$\\
\\
&& $-\frac{1}{144}(12n+11)-\frac{1}{16}(4n+3)-\frac{1}{24}(6n+5) +$\\
\\
&& $+\frac{1}{4}\cdot1+\frac{1}{4}\cdot1+\frac{1}{2}\cdot1+
   \frac{1}{9}\cdot(-1)+\frac{1}{24}\cdot(-1)+\frac{1}{24}\cdot(-1)=$\\
\\
&=&$-n,$
\end{tabular}
$$
$$
\begin{tabular}{rrl}
$\chi_h(GL_4\Z,S^{12n+10}V_4)$ &=&
 $\frac{1}{36}(6n+6)-\frac{1}{16}(4n+4)
      -\frac{1}{144}(12n+12)-\frac{1}{24}(6n+6) -$\\
\\
&& $-\frac{1}{144}(12n+12)-\frac{1}{16}(4n+4)-\frac{1}{24}(6n+6) +$\\
\\
&& $+\frac{1}{4}\cdot0+\frac{1}{4}\cdot0+\frac{1}{2}\cdot0+
   \frac{1}{9}\cdot0+\frac{1}{24}\cdot0+\frac{1}{24}\cdot0=$\\
\\
&=&$-1-n,$
\end{tabular}
$$

\subsection{Homological Euler characteristic of $GL_m(\Z[i])$}
\ble{8.5}
 For the symmetric power of the standard representation
of $GL_2(\Z[i])$ we have $$\tr([i^k,i^l]|S^{4n}V_2)=1,$$ for $k$ and
$l$ not congruent to each other modulo $4$.
 \ele
\proof First, consider the problem in dimension $1$. We have
$$\tr([i^k]|S^nV_1)=i^{kn},$$ To obtain the result in dimension $2$
we use that following observation: If $$g(n)=\tr(A|S^nV_1),$$
 $$h(n)=\tr(B|S^nV_2),$$ and $$f(n)=\tr([A,B]|S^n(V_1\oplus V_2),$$
then $$f(n)=(g*h)(n).$$ Therefore,
$$\tr([1,i]|S^{4n}V_2)=i^{4n}+i^{4n-1}+ \dots +1=1.$$ Because we take
$4$-th symmetric power we have
$$\tr([i^k,i^{k+1}]|S^{4n}V_2)=\tr([1,i]|S^{4n}V_2)=1.$$ Also,
$$\tr([1,-1]|S^{4n}V_2)=(-1)^{4n}+(-1)^{4n-1} + \dots +1=1.$$ And
similarly, $$\tr([i^k,i^{k+2}]|S^{4n}V_2)= \tr([1,-1]|S^{4n}V_2)=1.$$
Therefore,  $$\tr([i^k,i^l]|S^{4n}V_2)=1,$$ for $k$ and $l$ not
congruent to each other modulo $4$.

\bth{8.7}
  The homological Euler characteristic of $GL_2(\Z[i])$
with coefficients the symmetric powers of the standard
representation are given by\\
  $$\chi_h(GL_2(\Z[i]),S^nV_2)=
\left\{\begin{tabular}{ll} $1$ & $n\equiv 0 \mbox{ } \mod \mbox{
}4,$\\ $0$ & otherwise.\\
\end{tabular}\right.$$
  \eth
\proof If in the symmetric power $n$ is not divisible by $4$ then
the central element $[i,i]$ acts non-trivially on the cohomology
$$H^i(GL_2(\Z[i]),S^nV_2).$$ Therefore the cohomology and the
homological Euler characteristic vanish. If the symmetric power is
$4n$ then we can use lemma 4.4 and 5.7, and we obtain
$$\chi_h(GL_2(\Z[i]),S^{4n}V_2)=\sum R(A)\chi(C(A))\tr(A|S^{4n}V_2)
=4\cdot\frac{1}{8} + 2\cdot\frac{1}{4}=1.$$
\subsection{Homological Euler characteristic of $GL_m(\Z[\xi_3])$}
\ble{8.6}
    For the symmetric power of the standard representation
of $GL_2(\Z[\xi_3])$, we have
$$\tr([\xi_6^k,\xi_6^l]|S^{6n}V_2)=1,$$ for $k$ and $l$ not congruent to
each other modulo $6$.
  \ele
\proof  First, consider the problem in dimension $1$. We have
$$\tr([\xi_6^k]|S^nV_1)=\xi_6^{kn}.$$ Then
$$\tr([1,\xi_6]|S^{6n}V_2)=\xi_6^{6n}+\xi_6^{6n-1}+ \dots +1=1.$$
Since $\xi_6^6=1$ we have
$$\tr([\xi_6^{k},\xi_6^{k+1}]|S^{6n}V_2)=\tr([1,\xi_6]|S^{6n}V_2)=1.$$
Similarly, 
$$\tr([1,\xi_6^2]|S^{6n}V_2)=\xi_3^{6n}+\xi_3^{6n-1}+
\dots +1=1.$$ 
Again,
$$\tr([\xi_6^{k},\xi_6^{k+2}]|S^{6n}V_2)=\tr([1,\xi_6^2]|S^{6n}V_2)=1.$$
And $$\tr([1,-1]|S^{6n}V_2)=(-1)^{6n}+(-1)^{6n-1}+ \dots +1=1.$$
Also, 
$$\tr([\xi_6^{k},\xi_6^{k+3}]|S^{6n}V_2)=
\tr([1,-1]|S^{6n}V_2)=1.$$ 
Thus, 
$$\tr([\xi_6^k,\xi_6^l]|S^{6n}V_2)=1,$$
for $k$ and $l$ not congruent to each other modulo $6$.


\bth{8.8}
  The homological Euler characteristic of $GL_2\Z[\xi_3]$
with coefficients the symmetric powers of the standard
representation are given by\\
  $$\chi_h(GL_2(\Z[\xi_3]),S^nV_2)=
\left\{\begin{tabular}{ll} $1$ & $n\equiv 0 \mbox{ } \mod \mbox{
}6,$\\
  $0$ & otherwise.\\
\end{tabular}\right.$$
 \eth
\proof If in the symmetric power $n$ is not divisible by $6$ then
the central element $[\xi_6,\xi_6]$ acts non-trivially on the
cohomlogy $$H^i(GL_2(\Z[\xi_3]),S^nV_2).$$ Therefore the cohomology and
the homological Euler characteristic vanish. If the symmetric
power is $6n$ then we can use lemma 4.5 and 5.9, and we obtain
$$\chi_h(GL_2(\Z[\xi_3]),S^{6n}V_2)=\sum
R(A)\chi(C(A))\tr(A|S^{6n}V_2) =6\cdot\frac{1}{36}+
6\cdot\frac{1}{12} + 3\cdot\frac{1}{9}=1.$$


\sectionnew{Homological Euler characteristic of $\Gamma_1(m,N)$\\
and $\Ga_1(m,{\mathfrak{a}})$.}



We are going to compute the homological Euler characteristic of
$\Gamma_1(2,N)$, $\Gamma_1(3,N)$ and $\Ga_1(4,N)$ using 
the torsion elements in
the groups. Recall $\Gamma_1(m,N)$ is a subgroup of $GL_m\Z$
stabilizing the covector $$[0,\cdots,0,\: 1]$$ modulo $N$. The same
notation is used in \cite{G1}. Recall, also,  the
homologocal Euler characteristic of a group $G$ with coefficients
a finite dimensional representation $V$ over a field $k$ over  is

$$
   \chi_{_h}(G,V)=\sum_i (-1)^i dim_{k}(H^i(G, V)).
   $$
   We are going to use a generalization of a formula due to K. Brown
(see \cite{B2}) 
$$\chi_{_h}(G,V)=\sum_{A\e\cal{C}}\chi(C(A))\tr(A^{-1}|V),$$ 
where $\cal{C}$ consists of
representatives of the torsion elements of $G$ up to conjugation,
and $C(A)$ is the centralizer of $A$.In order to make use of this
formula, we have to find the torsion elements of $\Gamma_1(m,N)$
using the torsion elements in $GL_m\Z$. In a more general setting
let $G$ be a group and $\Gamma$ be a subgroup. It is possible to
find two elements of $\Gamma$ which are conjugate to each other in
$G$ but not in $\Gamma$. In the case of $\Ga=\Ga_1(m,N)$ and
$G=GL_m\Z$, we need to find all element of $\Ga$ which are not
conjugate to any of the others however they become conjugate when
considered as elements of the bigger group $G$. To do that we
construct a set $N^\Ga_G(A)$ with the property that the double
quotient $$\Ga\backslash N^\Ga_G(A)/C_G(A)$$ parametrizes the
elements in $\Ga$ that are non-conjugate to each other but all
become conjugate to $A$ when considered as elements of $G$.
The same set  $N^\Ga_G(A)$ is used in K. Brown paper \cite{B2}.
\ble{6.1}
   Let $G$ be a group and $\Gamma$ be a subgroup.
Given an element $A\in \Gamma$, the set of elements in which are
conjugate to A in $G$ but not in $\Gamma$ is parametrized by the
elements of the double quotient

$$
  \Ga\backslash N_G^{\Gamma}(A)/ C_G(A),
    $$
where $N_G^{\Gamma}(A)=\{X\in G: XAX^{-1}\in \Gamma\}$ and $C_G(A)$
is the centralizer of $A$ inside the group $G$.
\ele
\proof Let $A_1$ and $A_2$ be two elements of $\Gamma$
both conjugate to $A$ in $G$ and Conjugate to each other in
$\Gamma$. Then there exist $X_1,\: X_2$ in the bigger group $G$
such that $$A_1=X_1AX_1^{-1}$$ and $$A_2=X_2AX_2^{-1}.$$ Since $A_1$
and $A_2$ are conjugate in the smaller group $\Ga$, there exists
$Y\in \Ga$ such that $$YA_1Y^{-1}=A_2.$$ Then
$$YX_1AX_1^{-1}Y^{-1}=X_2AX_2^{-1}$$ and, equivalently,
$$(X_2^{-1}YX_1)A=A(X_2^{-1}YX_1).$$ We obtain that $$X_2^{-1}YX_1 \e
C_G(A).$$ Therefore $YX_1\in X_2C_G(A)$ and $X_1 \in \Ga X_2C_G(A)$
We obtain that $X_1$ and $X_2$ belong to the same double quotient
in $$\Ga \backslash N_G^{\Ga}(A)/C_G(A).$$ Conversely, suppose $X_1$
and $X_2$ belong to the same double quotient. We are going to show
that both elements $X_1AX_1^{-1}$ and $X_2AX_2^{-1}$ belong to
$\Ga$, and also that they are conjugated to each other in $\Ga$.
By definition of $N_G^{\Ga}(A)$ we have that $X_1AX_1^{-1}$ and
$X_2AX_2^{-1}$ belong to $\Ga$. Since the two elements belong to
the same double quotient, we have that $X_2=YX_1C$, where $Y\in
\Ga$ and $C\in C_G(A)$. Then
$$X_2Ax_2^{-1}=YX_1CAC^{-1}X_1^{-1}Y^{-1}= Y(X_1AX_1^{-1})Y^{-1}.$$
Thus, the elements $X_1AX_1^{-1}$ and $X_2AX_2^{-1}$ are conjugate
to each other in $\Ga$.

\subsection{Homological Euler characteristic of $\Ga_1(m,N)$}
\ble{6.2}
Let $A$ be an $k$-torsion element of $\Ga_1(m,N)$.
If $k$ and $N$ are relatively prime then $A$ has an eigenvalue
$+1$.
\ele
\proof Let $L$ be the field obtained by adjoining all eigenvalues of $A$
to $\Q$. Let $p$ be a prime dividing $N$. And let ${\mathfrak{p}}$ be
a prime ideal in $L$ sitting above $p$. Denote, also by $\mu(L)$ 
the roots of $1$ in $L$. And let 
$$\pi:L\rightarrow \F_q,$$
where
$$\F_q=L/{\mathfrak{p}}.$$ 
Then $\pi$ maps $\mu(L)$ onto $\F_q$ because there are no $p$-roots of $1$
in $L$ since $p$ does not divide $k$. So the eigenvalues on $A$ in $L$
are mapped onto $\F_q$. However, $1$ is an eigenvalue of 
$A\mbox{ }\mod \mbox{ }p$ in $\F_q.$ Therefore $1$ is an eigenvalue of
$A$ in $L$.


By proposition 2.8 part (a) we know that for a torsion element $A$, we have
$\chi(C(A))\neq 0$ only when the eigenvalues of $A$ are among
$$\{1,-1,i,-i,\xi_3,\bar\xi_3,\xi_6,\bar\xi_6\}.$$

 \ble{6.3}
Let $A\in\Ga_1(m,N)$ be a torsion element with $\chi(C(A))\neq 0$,
let $N$ be relatively prime with $2$ and $3$ and let $f$ be the
characteristic polynomial of $A$. Then +1 is a root of $f$ and it
has multiplicity $1$ or $2$. \\

   (a) If the root $+1$ has
multiplicity $1$ then the set 
$$\Ga_1(m,N)\backslash N^{\Ga_1(m,N)}_{CL_m\Z}(A)/C_{GL_m \Z}(A)$$ 
has $\frac{1}{2}\varphi(N)$ elements.\\

(b) If the root $+1$ has multiplicity $2$ then the set
$$\Ga_1(m,N)\backslash N^{\Ga_1(m,N)}_{CL_m\Z}(A)/C_{GL_m \Z}(A)$$
has one element.\\
  \ele
\proof (a) In order to save some space set $$\Ga=\Ga_1(m,N)$$ and
$$G=CL_m\Z.$$ Let $X\in  N^{\Ga}_G(A)$ and let $B=XAX^{-1}$. Then
$B\in \Ga_1(m,N)$. We can write the matrices $A$, $B$ and $X$ in
$2\times 2$ block form with diagonal blocks $A_{11}$, $B_{11}$ and
$X_{11}$ of size $(m-1)\times (m-1)$ and $A_{22}$, $B_{22}$ and
$X_{22}$ of size $1\times1$. Then

 $$
   \left[\begin{tabular}{cc}
$X_{11}$ & $X_{12}$\\ $X_{21}$ & $X_{22}$\\
\end{tabular} \right]
\cdot
  \left[\begin{tabular}{cc}
 $A_{11}$ & $A_{12}$\\
   $0$      & $1$     \\
\end{tabular} \right]
\equiv
  \left[\begin{tabular}{cc}
  $B_{11}$ & $B_{12}$\\
   $0$      & $1$     \\
  \end{tabular} \right]
\cdot
 \left[\begin{tabular}{cc}
  $X_{11}$ & $X_{12}$\\
  $X_{21}$ & $X_{22}$\\
   \end{tabular} \right].
  \mod\: N
   $$
Consider the $(2,1)$-block of the product. The left hand side is
$$X_{21}A_{11}$$ and the right hand side is $$1\cdot X_{21}\: \mod \:
N.$$ Examine the map $$P_{B_{22}A_{11}}:X_{21}\mapsto 
X_{21}A_{11}-B_{22}X_{21}$$ with $B_{22}=1$. The eigenvalues of
$A_{11}$ and $B_{22}$ are different. From lemma 1.11 follows that
the eigenvalues of $$P_{B_{22}A_{11}}$$ are $\l_i-1$ where $\l_i$
runs through the eigenvalues of $A_{11}$. From the previous lemma
we have that $\l_i - 1$ is invertible $\mod\:N$. Thus,
$$P_{B_{22}A_{11}}$$ is non-singular $\mod \: N$. Therefore,
$$P_{B_{22}A_{11}}(X_{21})= X_{21}A_{11}-B_{22}X_{21} \equiv 0
\mod\: N$$ implies that $$X_{21}\equiv 0\mod \:N.$$ Therefore,
$$N^{\Ga}_G(A) \subset \Ga_0(m,N),$$
where $\Ga_0(m,N)$ is defined as the subgroup of $GL_m \Z$
that sends $[0,\dots,0,1]$ to $[0,\dots,0,a]$ modulo $N$ for any 
$a\neq 0 \mbox{ }\mod \mbox{ }N.$

On the other hand, $\Ga_0(m,N)$
lies inside the normalizer of $\Ga_1(m,N)$. 
Also,
$$\Ga_0(m,N)\subset N^{\Ga}_G(A)$$ from the definition of
$N^{\Ga}_G(A)$. We obtain that $$N^{\Ga}_G(A)=\Ga_0(m,N).$$ Let $X$ be
an element in $\Ga_0(m,N)$ of the same block form as
before. Inside the quotient $$\Ga_1(m,N)\backslash \Ga_0(m,N)$$ the
element $X$ is determined by $$X_{22} \:\mod \: N.$$ Note that there
are $\varphi (N)$ options for $X_{22}$ because it has to be
invertible modulo $N$. From the block triangular theorem (theorem
2.3) we know that we can choose $A$ so that $A_{21}=0$. If $C$ is
in the centralizer of $A$ inside $GL_m\Z$, by lemma 3.2 we know
that $C_{21}=0$. Therefore the centralizer of $A$ modulo
$\Ga_1(m,N)$ is determined by $C_{22}$ which could be $+1$ or
$-1$. Therefore the double quotient consist of
$$\frac{1}{2}\varphi(N)$$ elements.

   (b) If $n=2$ then $A=I_2$
$$N^{\Ga}_G(A)=GL_2\Z$$ and $$C_{GL_2\Z}(A)=GL_2\Z.$$ Therefore, the
double quotient has one element. Assume that $m>2$. By the block
triangular theorem we can conjugate $A$ to a matrix of the block
form
 $$
   \left[
\begin{tabular}{ccc}
  $A_{11}$ & $A_{12}$ & $A_{13}$ \\
   $0$      & $1$      & $0$\\
   $0$      & $0$      & $1$\\
\end{tabular}
\right]$$
   Note that the block $A_{11}$ is an $(m-2)\times (m-2)$
matrix and it does not have an eigenvalue $1$. We can assume that
$A$ is of the above form. Let $X$ be an element of $N^{\Ga}_G(A)$,
and let $B=XAX^{-1}$. If we write this equation with respect to
the block form of $A$ we obtain:

$$
  \left[\begin{tabular}{ccc}
   $X_{11}$ & $X_{12}$ & $X_{13}$\\
  $X_{21}$ & $X_{22}$ & $X_{23}$\\
   $X_{31}$ & $X_{32}$ & $X_{33}$\\
   \end{tabular} \right]
\cdot \left[\begin{tabular}{ccc}
  $A_{11}$ & $A_{12}$ & $A_{13}$ \\
   $0$      & $1$      & $0$\\
   $0$      & $0$      & $1$\\
   \end{tabular} \right]
   \equiv
    $$
\\

$$
   \equiv
\left[\begin{tabular}{ccc}
   $B_{11}$ & $B_{12}$ & $B_{13}$ \\
   $B_{21}$ & $B_{22}$ & $B_{23}$\\
   $0$      & $0$      & $1$\\
\end{tabular} \right]
\cdot
   \left[\begin{tabular}{ccc}
   $X_{11}$ & $X_{12}$ & $X_{13}$\\
 $X_{21}$ & $X_{22}$ & $X_{23}$\\
 $X_{31}$ & $X_{32}$ & $X_{33}$\\
\end{tabular} \right]
\mod\: N
  $$
  Consider the $(3,1)$-block of the product. From the
left hand side we obtain $$X_{31}A_{11}$$ and from the right hand
side we obtain $1\cdot X_{31}$. We need to examine the map
$$X_{31}\mapsto X_{31}A_{11}- 1\cdot X_{31} \:\mod\:N.$$ Set
$$B_{33}=1.$$ We have to consider the map
$$P_{B_{33}A_{11}}:X_31\mapsto  X_{31}A_{11}- 1\cdot X_{31}.$$ If
$\l_k$ are the eigenvalues of $A_{11}$ then $\l_k-1$ are the
eigenvalues of $P_{B_{33}A_{11}}$ from the previous lemma it
follows that $\l_k-1$ is not divisible by $p$ for any $p$ dividing
$N$. This, the map $$P_{B_{33}A_{11}}$$ is non-singular modulo $N$.
We have $$P_{B_{33}A_{11}}(X_{31})\equiv 0\:\mod\:N.$$ Therefore
$$X_{31}\equiv 0\:\mod\:N.$$ Let $$N_{31}=\{(X_{ij})\e GL_m\Z: 1\leq
i,j\leq 3, X_{31}\equiv 0 \mbox{ }\mod\mbox{ }N\}.$$ Then
$$N^{\Ga}_G(A)\subset N_{31}.$$ We are going to show that the
inclusion is in fact an equality. We are going to prove the
following inclusions modulo $N$:

$$
   N^{\Ga}_G(A)\subset N_{31}
\subset \Ga_1(m,N)\cdot H'
  \subset \Ga_1(m,N)\cdot C_{GL_m\Z}(A)
\subset N^{\Ga}_G(A),
  $$
  where $H'$ is a subgroup of $C_{GL_m\Z}(A).$
The last inclusion always holds, not only modulo $N$. Modulo any
prime number $q$ relatively prime to $N$, we have
$$\Ga_1(m,N)=SL^{\pm}_m (\Z/q\Z)=N^{\Ga}_G(A),$$
where $SL^{\pm}_m(\Z/q\Z)$ is the group of $m\times m$ matrices 
with coefficients in $\Z/q\Z$ whose determinant is $\pm 1$.
 Therefore, the only
restriction on $N^{\Ga}_G(A)$ become apearant modulo $N$. From
these inclusions it will follow that $$N^{\Ga}_G(A)=
\Ga_1(m,N)\cdot C_{GL_m\Z}(A).$$ Therefore the double quotient
$$\Ga_1(m,N)\backslash N^{\Ga}_G(A)/C_{GL_m\Z}(A)$$ will have one
element.

  We have proved the first inclusion. The last
inclusion follows directly from the definition of $N^{\Ga}_G(A)$.
And the second to the last inclusion holds because $H'$ is a
subgroup of $C_{GL_m\Z}(A)$. Now we are going to prove the second
inclusion.\\
 Consider $A$ as a $2\times 2$-block matrix
$$
   A=
\left[\begin{tabular}{cc}
  $A_{11}$ & $\overline{A}_{12}$\\
  $0$      & $I_2$     \\
\end{tabular} \right],
   $$
with $$\overline{A}_{12}=[A_{12},\: A_{13}].$$ Let $$C\in
C_{GL_n\Z}(A)$$ be in the centralizer of $A$. Consider $C$ as a
$2\times 2$ block matrix. Namely,

$$
  C=
\left[\begin{tabular}{cc} $C_{11}$ & $C_{12}$\\
   $0$      &
$C_{22}$     \\
\end{tabular} \right],
     $$
  with $C_{22}$ a $2\times 2$ matrix.
The admissible matrices $C$ are such that
$C_{11}\overline{A}_{12}C_{22}^{-1}$ and $\overline{A}_{12}$ map
to the same element in $Q_{mod}$, where $$Q_{mod}=Mat_{n-2,2}\Z/\Im
\mbox{ }P_{A_{11},I_2},$$ the matrix $C_{12}$ is determined
uniquely from $C_{11}$ and $C_{22}$. Let $$H=\{C\in C(A):
C_{11}=I_{m-2}\}$$ be a subgroup of $C(A)$. We are going to show
that $$N_{31}\subset \Ga_1(m,N)\cdot H',$$ where $$H'\subset H
\subset C(A)$$ are subgroups. In this setting $C_{12}$ is uniquely
determined by $C_{11}$; and we are going to write it as
$$C_{12}(C_{11}).$$ Let $$M=det(P_{A_{11},I_2}).$$ Then we have that
$$M\cdot P_{A_{11},I_2}^{-1}$$ has integer entries. Therefore, $$\Im
\mbox{ }(P_{A_{11},I_2})\subset M\cdot Mat_{m-2,2}\Z.$$ If $C_{22}$
is in $\Ga(2,M)$ then $$\overline{A}_{12} C_{22}^{-1}\equiv
\overline{A}_{12}\mbox{ } \mod \mbox{ } M\cdot Mat_{m-2,2}\Z.$$
Therefore $$\overline{A}_{12} C_{22}^{-1} \equiv
\overline{A}_{12}\mbox{ } \mod \mbox{ }\Im \mbox{
}(P_{A_{11},I_2}).$$ Set

 $$
 H'=\{C\in C(A): C_{11}=I_{m-2}, C_{22}\in
\Ga(2,M),C_{12}=C_{12}(C_{22})\}.
    $$
    Then $$H'\cong \Ga(2,M).$$
However, we need $H'$ for comparison of groups modulo $N$. Note
that $M$ and $N$ are relatively prime because
$$M=det(P_{A_{11}I_2})$$ which has prime factors only $2$ and $3$.
And by assumption $N$ is relatively prime to $2$ and $3$ Then
$$\Ga(2,M)\equiv SL_2\Z \mbox{ }\mod\mbox{ } N.$$ We need to show
that $$N_{31}\subset \Ga_1(n,N)\cdot H'$$ modulo $N$. Let $X\in
N_{31}$. We have that $X_{32}$ and $X_{33}$ are relatively prime
modulo $N$. Otherwise their common factor will divide $\det(X)$,
which is $\pm 1$. Then we can find $Y$ and $Z$ such that
$$YX_{33}-ZX_{32}=1.$$ Let

   $$
   C_{22}=
\left[\begin{tabular}{cc}
 $Y$      & $Z$     \\
  $X_{32}$ & $X_{33}$\\
\end{tabular} \right].
  $$
   Then
Let $C$ be a corresponding element in $H'$ with $C_{22}$ congruent
to the above matrix modulo $M$ and $$C_{11}=I_{n-2}.$$ Then
$$XC^{-1}\in \Ga_1(m,N),$$ and we are done.\\


\bth{6.4}
    Let $V$ be a representation of $GL_m\Z$, and let
$N$ be an integer relatively prime to $2$ and $3$. Then\\

  $$
  \begin{tabular}{cl}
$\chi_h(\Ga_1(m,N),V)=$ &
 $\varphi(N)\sum_{A=[A_1,1]} R(A)\chi(C_{GL_{m-1}\Z}(A_1))\tr(A^{-1}|V)+$ \\
 & $+\varphi_2(N)\sum_{A=[A_2,I_2]}
R(A)\chi(C_{GL_{m-2}\Z}(A_2))\tr(A^{-1}|V),$\\
\end{tabular}$$
\\
  where $A_1$ and $A_2$ are block-diagonal matrices with zero blocks off the
diagonal, and with block on the diagonal
$A_{ii}$ varying through $$\{-1,-I_2,T_3,T_4,T_6\}$$ and $A_{ii}$ and
$A_{jj}$ have no common eigenvalues. We set 
$$T_3=
\left[
\begin{tabular}{rr}
$0$ & $1$\\
$-1$ & $-1$
\end{tabular}
\right],
\mbox{ }
T_4=
\left[
\begin{tabular}{rr}
$0$ & $1$\\
$-1$ & $0$
\end{tabular}
\right],
\mbox{ }
T_6=
\left[
\begin{tabular}{rr}
$0$ & $-1$\\
$1$ & $1$
\end{tabular}
\right].
$$
Also $\varphi(N)$ is the Euler
function of $N$, and $\varphi_2(N)$ is the arithmetic function
generated by $$\varphi_2(p^n)=p^{2k}(1-\frac{1}{p^2}).$$
       \eth
   \proof Let $A$ be a torsion element of $$\Ga=\Ga_1(m,N)$$
such that $$\chi(C_{\Ga}(A))\neq 0.$$ Then $$\chi(C_{GL_m\Z}(A))\neq
0.$$ By lemma 6.2 we have that $1$ is an eigenvalue of $A$. Also,
if  $$\chi(C_{\Ga}(A))\neq 0$$ then the multiplicity of the
eigenvalue $1$ is at most $2$. We are going to prove that if the
eigenvalue has multiplicity $1$ then
  $$\chi(C_{\Ga}(A))=2\cdot \chi(C_{GL_m\Z}(A)),$$ and
if the eigenvalue $1$ has multiplicity $2$ then
 $$\chi(C_{\Ga}(A))=\varphi_2(N)\cdot \chi(C_{GL_m\Z}(A)).$$
Assume that we have proven these two formulas. From the
generalization of Brown's formula, we have
 $$\chi_h(\Ga_1(m,N),V)=\sum_{A:\mbox{ torsion}}
\chi(C_{\Ga_1(m,N)}(A))\tr(A^{-1}|V).$$ In the summation it is enough to
sum over torsion elements $A$ that have eigenvalue $1$ either with
multiplicity $1$ or with multiplicity $2$. First take the sum over
elements $A$ that become conjugate to each other in $GL_m\Z$. Let
$A$ has eigenvalue $1$ with multiplicity $1$. By lemma 6.3 $$\sum
\chi(C_{\Ga_1(m,N)} (A''))\tr(A''^{-1}|V)
=\varphi(N)\chi(C_{GL_n\Z}(A'))\tr(A'^{-1}|V),$$ where the sum is taken
over all non-conjugate $A''$ in $\Ga_1(m,N)$ that become conjugate
to $A'$ in $GL_n\Z$. Now sum over the torsion elements $A'$ in
$GL_m\Z$ that become conjugate to each other in $GL_m\C$. 

Similarly, in the case
when the eigenvalue $1$ has multiplicity $2$, we have
$$\varphi_2(N)=[\Ga_1(m,N):GL_m\Z]$$ instead of $\varphi(N)$. Thus, it
remains to prove that if the eigenvalue $1$ has multiplicity $1$
then
  $$\chi(C_{\Ga}(A))=2\cdot \chi(C_{GL_m\Z}(A)),$$ and
if the eigenvalue $1$ has multiplicity $2$ then
 $$\chi(C_{\Ga}(A))=\varphi_2(N)\cdot \chi(C_{GL_m\Z}(A)).$$\\
Suppose $A$ has eigenvalue $1$ with multiplicity one. We can
assume that $A$ is in block-triangular 
$$
\left[\begin{tabular}{cc}
$A_{11}$ & $A_{22}$ \\
   0     & $1$\\
\end{tabular}\right]
$$
with $A_{22}=1$ and $A_{21}=0$. If $C$ is a matrix in $GL_m\Z$
commuting with $A$ then $C$ is if the same block form with
$C_{21}=0$ and $C_{22}=\pm 1$. Exactly one of the matrices $C$ and
$-C$ belongs to $\Ga_1(m,N)$ because of the $C_{22}$ entry. Thus
the centralizer of $A$ in $GL_m\Z$ contains twice as many element
as the centralizer inside $\Ga(m,N)$. Therefore, 
$$\chi(C_{\Ga}(A))=2\cdot \chi(C_{GL_m\Z}(A)).$$

Suppose $A$ has eigenvalue $1$ with multiplicity two. We can
assume that $A$ is of the form 
$$
\left[\begin{tabular}{cc}
  $A_{11}$
& $A_{22}$ \\
   0     & $I_2$\\
\end{tabular}\right]
$$
with $$A_{22}=I_2$$ and $$A_{21}=0.$$ If $C$ commutes with $A$ then it
is of the same block type with $$C_{21}=0.$$ The possibilities of
$C_{22}$ are determined by relations modulo the image of
$$P_{A_{11}I_2}.$$ Note that if $$\chi(C(A))\neq 0$$ then the
eigenvalues of $A$ are among $$\{\pm1,\pm i,\pm \xi_3\}.$$ Then
$$\det(P_{A_{11},I_2})$$ has only prime factors $2$ and $3$.
Therefore the conditions on $C_{22}$ imposed by $$P_{A_{11},I_2}$$
are independent from the condition $C\e \Ga_1(n,N)$, because that
last condition is a congruence modulo $N$ but $N$ is relatively
prime $2$ and $3$. Therefore
 $$
\begin{tabular}{ccl}
$\chi(C_{\Ga_1(m,N)}(A))$ & $=$
 & $[C_{\Ga_1(m,N)}(A):C_{GL_m\Z}(A)] \chi(C_{GL_m\Z}(A))=$\\
\\
 & $=$ & $\varphi_2(N)\cdot \chi(C_{GL_m\Z}(A)).$
\end{tabular}
$$


\bco{6.5}
   For $N$ not divisible by $2$ and $3$
the homological Euler characteristics of $\Ga_1(2,N)$ with
coefficients the symmetric powers of $V_2$ are given by

$$
   \chi_h(\Ga_1(2,N),S^{2n+k}V_2)=
\left\{
\begin{tabular}{ll}
$-\frac{1}{24}\varphi_2(N)(2n+1)+\frac{1}{2}\varphi(N)$ & $k=0$\\
\\
$-\frac{1}{24}\varphi_2(N)(2n+2)$ & $k=1$
\end{tabular}
\right.
    $$
where $\varphi(N)$ is the multiplicative Euler function generated
by $$\varphi(p^n)=p^n(1-\frac{1}{p})$$ and $\varphi_2(N)$ is the
multiplicative function generated by
$$\varphi_2(p^n)=p^{2n}(1-\frac{1}{p^2}).$$
   \eco
\proof From theorem 6.4, we have that
$$
\begin{tabular}{lll}
$\chi_h(\Ga_1(2,N),S^{2n}V_2)$ & = &
   $\varphi(N)R([-1,1])\chi(C([-1,1]))\tr([-1,1]|S^{2n}V_2)+$\\
\\
&  & $\varphi_2(N)\chi(\Ga_1(2,N))\tr(I_2|S^{2n}V_2)=$\\ 
\\
&= & $\varphi(N)\cdot 2\cdot \frac{1}{4}\cdot 1
       + \varphi_2(N)(-\frac{1}{24})(2n+1)=$\\
\\
&= & $-\frac{1}{24}\varphi_2(N)(2n+1)+\frac{1}{2}\varphi(N).$
\end{tabular}
$$

Similarly,
$$
\begin{tabular}{lll}
$\chi_h(\Ga_1(2,N),S^{2n+1}V_2)$ &=&
 $\varphi(N)R([-1,1])\chi(C([-1,1]))\tr([-1,1]|S^{2n+1}V_2)+$\\ 
\\
 & & $+\varphi_2(N)\chi(\Ga_1(2,N))\tr(I_2|S^{2n+1}V_2)=$\\
\\
&=& $\varphi(N)\cdot 2\cdot \frac{1}{4}\cdot 0+
\varphi_2(N)(-\frac{1}{24})(2n+2).$
\end{tabular}
$$

\bth{6.6}
  The torsion elements of $\Ga_1(3,N)$, for $N$ not
divisible by $2$ or $3$, are given in the following table together
with the Euler characteristic of their centralizers. In the table
$$\varphi_2(N)=[\Ga_1(2,N):GL_2\Z]$$ is the multiplicative function
generated by $$\varphi_2(p^n)=p^{2n}(1-\frac{1}{p^2}).$$ For
abbreviation  set $$\Ga=\Ga_1(3,N) \mbox{ and }G=GL_3\Z.$$\\
$
\begin{tabular}{|c|c|c|c|}
\hline $A$                          & $|\Ga\backslash
N_G^{\Ga}(A)/C_G(A)|$  &
 $C_{\Ga}(A)$ & $\chi(C_{\Ga}(A))$\\ \hline

(a) $\left[\begin{tabular}{ccc} $1$ &     & \\
    & $1$ & \\
    &     & $1$\\
\end{tabular}\right]$        &
$1$                  & $\Ga_1(3,N)$ & $0$\\ \hline

(c1)
 $\left[\begin{tabular}{c|cc}
  $-1$ &     &      \\ \hline
     & $1$ & $0$  \\
     & $0$ & $1$  \\
\end{tabular}\right]$        &
$1$ & $C_2\times \Ga_1(2,N)$ & $-\frac{1}{48}\varphi_2(N)$\\
\hline

(c2) $\left[\begin{tabular}{c|cc}
   $-1$ & $0$ & $1$\\ \hline
  $0$ & $1$ & $0$\\
    &     & $1$\\
\end{tabular}\right]$        &
$1$ & $C_2\times \Ga_1(2,2N)$ & $-\frac{1}{16}\varphi_2(N)$\\
\hline

(d1)
  $\left[\begin{tabular}{cc|c}
 $-1$ & $0$ & \\
 $0$ & $-1$ & \\ \hline
    &      & $1$\\
\end{tabular}\right]$
& $\frac{1}{2}\varphi(N)$ & $GL_2\Z$ & $-\frac{1}{24}$\\ \hline

(d2)
   $\left[\begin{tabular}{cc|c}
  $-1$ & $0$ & $-1$\\
  $0$ & $-1$ & $0$\\ \hline
    &      & $1$\\
\end{tabular}\right]$
& $\frac{1}{2}\varphi(N)$ & $\Gamma_1(2,2)$ & $-\frac{1}{8}$\\
\hline

(e1)
 $\left[\begin{tabular}{cc|c}
 $0$ & $1$   & \\
 $-1$ & $-1$ & \\ \hline
    &     & $1$\\
\end{tabular}\right]$        &
 $\frac{1}{2}\varphi(N)$ & $C_6$ & $\frac{1}{6}$\\ \hline

(e2)
 $\left[\begin{tabular}{cc|c}
$0$ & $1$   & $1$\\
   $-1$ & $-1$ & \\ \hline
    &     & $1$\\
\end{tabular}\right]$        &
 $\frac{1}{2}\varphi(N)$  & $C_3$ & $\frac{1}{3}$\\ \hline

(h)
 $\left[\begin{tabular}{cc|c}
 $0$ & $-1$ & \\
  $1$ & $1$  & \\ \hline
    &      & $1$\\
\end{tabular}\right]$
& $\frac{1}{2}\varphi(N)$ & $C_6$ & $\frac{1}{6}$\\ \hline

(i1)
   $\left[\begin{tabular}{cc|c}
   $0$  & $1$ & \\
   $-1$ & $0$ & \\ \hline
     &     & $1$\\
\end{tabular}\right]$        &
 $\frac{1}{2}\varphi(N)$ & $C_4$ & $\frac{1}{4}$\\ \hline

(i2)
   $
   \left[\begin{tabular}{cc|c}
   $0$  & $1$ & $1$\\
  $-1$ & $0$ & \\ \hline
     &     & $1$\\
\end{tabular}\right]$        &
 $\frac{1}{2}\varphi(N)$ & $C_4$ & $\frac{1}{4}$\\ \hline
\end{tabular}
   $
\eth

\proof This table follows from the table in proposition 3.2 and from the
previous two lemmas on the calculation of $$\Ga\backslash
N_G^{\Ga}(A)/C_G(A).$$ The only part that requires a proof given
two elements $A$ and $B$ in $\Ga_1(3,N)$ such that both are
conjugate in $GL_3\Z$ to the same element form the above list, we
have $$C_{\Ga_1(3,N)}(A)\cong C_{\Ga_1(3,N)}(B).$$ Indeed, if $A$
and $B$ are conjugate to each other in $\Ga_1(3,N)$ then we are
done. If they are not conjugate to each other in $\Ga_1(3,n)$ then
the double quotient $\Ga\backslash N_G^{\Ga}(A)/C_G(A)$ has more
than one element. From lemma 5.3 we have that
$$N_G^{\Ga}(A)=\Ga_0(3,N)$$ which is the normalizer of $\Ga_1(3,N)$
in $GL_3\Z$. Let $$P\e N_G^{\Ga}(A)$$ such that $$B=PAP^{-1}$$ then
$$C_{\Ga_1(3,N)}(B)=P\cdot C_{\Ga_1(3,N)}(A) \cdot P^{-1}$$ because
$P$ is in the normalizer of $\Ga_1(3,N)$ in $GL_3\Z.$


\bco{6.7}
    For $N$ not divisible by $2$ and $3$
the homological Euler characteristic of $\Ga_1(3,N)$ is given by

$$
  \chi_h(\Ga_1(3,N))=-\frac{1}{12}\varphi_2(N)+\frac{1}{2}\varphi(N).
  $$
where $\varphi(N)$ is the multiplicative Euler function generated
by $$\varphi(p^n)=p^n(1-\frac{1}{p})$$ and $\varphi_2(N)$ is the
multiplicative function generated by
$$\varphi_2(p^n)=p^{2n}(1-\frac{1}{p^2}).$$
  \eco
   \proof
Brown's formula gives the homological Euler characteristic as a
sum of the usual (orbifold) Euler characteristic of the
centralizers of the torsion elements. From the previous
proposition we obtain
$$
\begin{tabular}{ccl}
$\chi_h(\Ga_1(3,N))$ & $=$ &
 $-\frac{1}{48}\varphi_2(N)-\frac{1}{16}\varphi_2(N)
    -\frac{1}{48}\varphi(N)--\frac{1}{16}\varphi(N)+$\\
\\
 & &$+\frac{1}{12}\varphi(N)
+\frac{1}{6}\varphi(N)+\frac{1}{12}\varphi(N)
+\frac{1}{8}\varphi(N)+\frac{1}{8}\varphi(N)=$\\
\\
 &=& $-\frac{1}{12}\varphi_2(N)+\frac{1}{2}\varphi(N).$
\end{tabular}
$$

\bco{6.8}
   For $N$ not divisible by $2$ and $3$
the homological Euler characteristic of $\Ga_1(4,N)$ is given by

$$\chi_h(\Ga_1(4,N))=\varphi(N).$$
   \eco
   \proof
Using theorem 6.4 and lemma 4.3 we obtain
$$
\begin{tabular}{ccl}
$\chi_h(\Ga_1(4,N))$ &=& $\varphi_2(N)(\frac{1}{36}-
\frac{1}{16}-\frac{1}{144}-\frac{1}{24})+$\\
\\
 && $+\varphi(N)(\frac{1}{4}+\frac{1}{4}+\frac{1}{2})=$\\
\\
 &=& $\varphi(N).$
\end{tabular}
$$
\subsection{Homological Euler characteristic
of $\Ga_1(m,{\mathfrak{a}})\subset GL_m(\Z[i]).$}
\ble{6.13}
Let $A$ be an $k$-torsion element of $\Ga_1(m,{\mathfrak{a}})$.
If $2k$ and ${\mathfrak{a}}$ are relatively prime then $A$ has an eigenvalue
$+1$.
\ele
\proof Let $L$ be the field obtained by adjoining all eigenvalues of $A$
to $\Q(i)$. Let ${\mathfrak{p}}$ be a prime ideal dividing ${\mathfrak{a}}$. 
And let ${\mathfrak{P}}$ be
a prime ideal in $L$ sitting above ${\mathfrak{p}}$. Denote, also by $\mu(L)$ 
the roots of $1$ in $L$. And let 
$$\pi:L\rightarrow \F_q,$$
where
$$\F_q=L/{\mathfrak{P}}.$$ 
Let $p$ be the rational prime sitting below ${\mathfrak{p}}$.
Then $\pi$ maps $\mu(L)$ onto $\F_q$ because there are no $p$-roots of $1$,
in $L$ since ${\mathfrak{p}}$ does not divide $k$. 
So the eigenvalues on $A$ in $L$
are mapped onto $\F_q$. However, $1$ is an eigenvalue of 
$A\mbox{ }\mod \mbox{ }p$ in $\F_q.$ Therefore $1$ is an eigenvalue of
$A$ in $L$.


By proposition 2.8 part (b) we know that for a torsion element $A$, we have
$\chi(C(A))\neq 0$ only when the eigenvalues of $A$ are among
$$\{1,-1,i,-i\}.$$

 \ble{6.14}
Let $A\in\Ga_1(m,{\mathfrak{a}})$ be a torsion element with $\chi(C(A))\neq 0$,
let ${\mathfrak{a}}$ be relatively prime with $(1+i)$, and let $f$ be the
characteristic polynomial of $A$. Then +1 is a root of $f$ and it
has multiplicity $1$. And the set 
$$\Ga_1(m,N)\backslash N^{\Ga_1(m,{\mathfrak{a}})}_{GL_m\Z[i]}(A)/
C_{GL_m \Z[i]}(A)$$ 
has $\frac{1}{4}\varphi_{\Z[i]}({\mathfrak{a}})$ elements,
where $\varphi_{\Z[i]}({\mathfrak{a}})$ is a multiplicative function
on the ideals of $\Z[i]$ generated by
$$\varphi_{\Z[i]}((\mathfrak{p})^n)=N_{\Q(i)/\Q}({\mathfrak{p}})^n(1
-\frac{1}{N_{\Q(i)/\Q}({\mathfrak{p}})}).$$
  \ele
\proof  In order to save some space set $$\Ga=\Ga_1(m,{\mathfrak{a}})$$ and
$$G=CL_m(\Z[i]).$$ Let $X\in  N^{\Ga}_G(A)$ and let $B=XAX^{-1}$. Then
$B\in \Ga$. We can write the matrices $A$, $B$ and $X$ in
$2\times 2$ block form with diagonal blocks $A_{11}$, $B_{11}$ and
$X_{11}$ of size $(m-1)\times (m-1)$ and $A_{22}$, $B_{22}$ and
$X_{22}$ of size $1\times1$. Then

 $$
   \left[\begin{tabular}{cc}
$X_{11}$ & $X_{12}$\\ $X_{21}$ & $X_{22}$\\
\end{tabular} \right]
\cdot
  \left[\begin{tabular}{cc}
 $A_{11}$ & $A_{12}$\\
   $0$      & $1$     \\
\end{tabular} \right]
\equiv
  \left[\begin{tabular}{cc}
  $B_{11}$ & $B_{12}$\\
   $0$      & $1$     \\
  \end{tabular} \right]
\cdot
 \left[\begin{tabular}{cc}
  $X_{11}$ & $X_{12}$\\
  $X_{21}$ & $X_{22}$\\
   \end{tabular} \right].
  \mod\: {\mathfrak{a}}
   $$
Consider the $(2,1)$-block of the product. The left hand side is
$$X_{21}A_{11}$$ and the right hand side is 
$$1\cdot X_{21}\: \mod \:{\mathfrak{a}}.$$ 
Examine the map $$P_{B_{22}A_{11}}:X_{21}\mapsto 
X_{21}A_{11}-B_{22}X_{21}$$ with $B_{22}=1$. The eigenvalues of
$A_{11}$ and $B_{22}$ are different. From lemma 1.8 follows that
the eigenvalues of $$P_{B_{22}A_{11}}$$ are $\l_i-1$ where $\l_i$
runs through the eigenvalues of $A_{11}$. From the previous lemma
we have that $\l_i - 1$ is invertible $\mod\:{\mathfrak{a}}$. Thus,
$$P_{B_{22}A_{11}}$$ is non-singular $\mod \: {\mathfrak{a}}$. Therefore,
$$P_{B_{22}A_{11}}(X_{21})= X_{21}A_{11}-B_{22}X_{21} \equiv 0
\mod\: {\mathfrak{a}}$$ implies that $$X_{21}\equiv 0\mod \:{\mathfrak{a}}.$$ 
Therefore,
$$N^{\Ga}_G(A) \subset \Ga_0(m,{\mathfrak{a}}),$$
where $\Ga_0(m,{\mathfrak{a}})$ is defined as the subgroup of $GL_m(\Z[i])$
that sends $[0,\dots,0,1]$ to $[0,\dots,0,\a]$ modulo ${\mathfrak{a}}$ for any 
$\a\neq 0 \mbox{ }\mod \mbox{ }{\mathfrak{a}}.$

On the other hand, $\Ga_0(m,{\mathfrak{a}})$
lies inside the normalizer of $\Ga_1(m,{\mathfrak{a}})$. 
Also,
$$\Ga_0(m,{\mathfrak{a}})\subset N^{\Ga}_G(A)$$ from the definition of
$N^{\Ga}_G(A)$. We obtain that $$N^{\Ga}_G(A)=\Ga_0(m,{\mathfrak{a}}).$$ 
Let $X$ be
an element in $\Ga_0(m,{\mathfrak{a}})$
of the same block form as
before. Inside the quotient $$\Ga_1(m,{\mathfrak{a}})\backslash 
\Ga_0(m,{\mathfrak{a}})$$ the
element $X$ is determined by $$X_{22} \:\mod \: {\mathfrak{a}}.$$ 
Note that there
are $\varphi_{\Z[i]} ({\mathfrak{a}})$ options for $X_{22}$ 
because it has to be
invertible modulo ${\mathfrak{a}}$. From the block-triangular theorem (theorem
2.3) we know that we can choose $A$ so that $A_{21}=0$. If $C$ is
in the centralizer of $A$ inside $GL_m(\Z[i])$, by lemma 3.2 we know
that $C_{21}=0$. Therefore the centralizer of $A$ modulo
$\Ga_1(m,{\mathfrak{a}})$ is determined by $C_{22}$ which could be $+1$,
$-1$, $i$ or $-i$. Therefore the double quotient consist of
$$\frac{1}{4}\varphi_{\Z[i]}({\mathfrak{a}})$$ elements.


\bth{6.4}
    Let $V$ be a representation of $GL_m(\Z[i])$, and let
${\mathfrak{a}}$ be an ideal relatively prime to $(1+i)$. Then
$$
\chi_h(\Ga_1(m,{\mathfrak{a}}),V)=
\varphi_{\Z[i]}({\mathfrak{a}})\sum_{A=[A_0,1]} 
R(A)\chi(C_{GL_{m-1}\Z[i]}(A_0))\tr(A^{-1}|V), 
$$
  where $$A_0=[A_{11},\dots A_{kk}],$$
$A_{ii}$ vary through $$\{-1,i,-i\}$$ and $A_{ii}$ and
$A_{jj}$ are distinct,
$\varphi_{\Z[i]}({\mathfrak{a}})$ is the arithmetic function
defined on the ideals of $\Z[i]$
generated by 
$$\varphi_{\Z[i]}({\mathfrak{p}}^n)=N_{\Q(i)/\Q}({\mathfrak{p}})^{k}(1-
\frac{1}{N_{\Q(i)/\Q}({\mathfrak{p}})}).$$
       \eth
   \proof Let $A$ be a torsion element of $$\Ga=\Ga_1(m,{\mathfrak{a}})$$
such that $$\chi(C_{\Ga}(A))\neq 0.$$ Then $$\chi(C_{GL_m\Z[i]}(A))\neq
0.$$ By lemma 6.9 we have that $1$ is an eigenvalue of $A$. Also,
if  $$\chi(C_{\Ga}(A))\neq 0$$ then the multiplicity of the
eigenvalue $1$ is at most $1$. We are going to prove that 
  $$\chi(C_{\Ga}(A))=4\cdot \chi(C_{GL_m\Z[i]}(A)),$$ 
Assume that we have proven this formula. From the
generalization of Brown's formula, we have
 $$\chi_h(\Ga_1(m,{\mathfrak{a}}),V)=\sum_{A:\mbox{ torsion}}
\chi(C_{\Ga_1(m,{\mathfrak{a}})}(A))\tr(A^{-1}|V).$$ 
In the summation it is enough to
sum over torsion elements $A$ that have eigenvalue $1$ either with
multiplicity $1$. First take the sum over
elements $A$ that become conjugate to each other in $GL_m(\Z[i])$. Let
$A$ has eigenvalue $1$ with multiplicity $1$. By lemma 6.9 $$\sum
\chi(C_{\Ga_1(m,{\mathfrak{a}})} (A''))\tr(A''^{-1}|V)
=\varphi_{\Z[i]}({\mathfrak{a}})
\chi(C_{GL_n\Z[i]}(A'))\tr(A'^{-1}|V),$$ where the sum is taken
over all non-conjugate $A''$ in $\Ga_1(m,{\mathfrak{a}})$ that become conjugate
to $A'$ in $GL_n(\Z[i])$.

We can assume that $A$ is in block-triangular 
$$
\left[\begin{tabular}{cc}
$A_{11}$ & $A_{22}$ \\
   0     & $1$\\
\end{tabular}\right]
$$
with $A_{22}=1$ and $A_{21}=0$. If $C$ is a matrix in $GL_m(\Z[i])$
commuting with $A$ then $C$ is if the same block form with
$C_{21}=0$ and $C_{22}=i^k$. Exactly one of the matrices $C$, $-C$, $iC$ and
$-iC$ belongs to $\Ga_1(m,{\mathfrak{a}})$ because of the $C_{22}$ entry. Thus
the centralizer of $A$ in $GL_m(\Z[i])$ contains 4 times as many element
as the centralizer inside $\Ga(m,{\mathfrak{a}})$. Therefore, 
$$\chi(C_{\Ga}(A))=4\cdot \chi(C_{GL_m\Z[i]}(A)).$$
By the previous lemma we have that the number of non-conjugate matrices in 
$\Ga$ that become conjugate in $G$ is 
$\varphi_{\Z[i]}({\mathfrak{a}})/4.$
This proves the theorem.
\bth{6.12}
The homological Euler characteristic of 
$\Ga_1(2,{\mathfrak{a}})\subset GL_2(\Z[i])$ is given by
$$\chi_h(\Ga_1(2,{\mathfrak{a}}),\Q)=
\frac{1}{2}\varphi_{\Z[i]}({\mathfrak{a}}),$$
where $\varphi_{\Z[i]}({\mathfrak{a}})$ is the multiplicative function defined 
on the ideals of $\Z[i]$, generated by
$$\varphi_{\Z[i]}({\mathfrak{p}}^n)=N_{\Q(i)/\Q}({\mathfrak{p}})^n(1-
\frac{1}{N_{\Q(i)/\Q}({\mathfrak{p}})}).$$
\eth
\proof From theorem 6.11 and lemma 4.4 we have
$$
\begin{tabular}{rrl}
$\chi_h(\Ga_1(2,{\mathfrak{a}}),\Q)$&=
 &$\varphi_{\Z[i]}({\mathfrak{a}})(|N_{\Q(i)/\Q}(R([1,i]))|\chi(C([1,i]))+$\\
\\
    &&$+|N_{\Q(i)/\Q}(R([1,-i]))|\chi(C([1,-i]))+$\\
\\ 
   &&$+|N_{\Q(i)/\Q}(R([1,-1]))|\chi(C([1,-1])))=$\\
\\
 &=&$\varphi_{\Z[i]}({\mathfrak{a}})(\frac{1}{8}+\frac{1}{8}+\frac{1}{4})=$\\
\\
 &=&$\frac{1}{2}\varphi_{\Z[i]}({\mathfrak{a}}).$
\end{tabular}
$$

\subsection{Homological Euler characteristic
of $\Ga_1(m,{\mathfrak{a}})\subset GL_m(\Z[\xi_3]).$}
\ble{6.2}
Let $A$ be an $k$-torsion element of $\Ga_1(m,{\mathfrak{a}})$.
If $6k$ and ${\mathfrak{a}}$ are relatively prime then $A$ has an eigenvalue
$+1$.
\ele
\proof Let $L$ be the field obtained by adjoining all eigenvalues of $A$
to $\Q(\xi_3)$. Let ${\mathfrak{p}}$ be a prime ideal dividing 
${\mathfrak{a}}$. 
And let ${\mathfrak{P}}$ be
a prime ideal in $L$ sitting above ${\mathfrak{p}}$. Denote, also by $\mu(L)$ 
the roots of $1$ in $L$. And let 
$$\pi:L\rightarrow \F_q,$$
where
$$\F_q=L/{\mathfrak{P}}.$$ 
Let $p$ be the rational prime sitting below ${\mathfrak{p}}$.
Then $\pi$ maps $\mu(L)$ onto $\F_q$ because there are no $p$-roots of $1$,
in $L$ since ${\mathfrak{p}}$ does not divide $k$. 
So the eigenvalues on $A$ in $L$
are mapped onto $\F_q$. However, $1$ is an eigenvalue of 
$A\mbox{ }\mod \mbox{ }p$ in $\F_q.$ Therefore $1$ is an eigenvalue of
$A$ in $L$.


By proposition 2.8 part (b) we know that for a torsion element $A$, we have
$\chi(C(A))\neq 0$ only when the eigenvalues of $A$ are among
$$\{1,-1,\xi_3,-\xi_3, \xi_6, -\xi_6\}.$$

 \ble{6.3}
Let $A\e\Ga_1(m,{\mathfrak{a}})$ be a torsion element with $\chi(C(A))\neq 0$,
let ${\mathfrak{a}}$ be relatively prime with $(1-\xi_3)$, and let $f$ be the
characteristic polynomial of $A$. Then +1 is a root of $f$ and it
has multiplicity $1$. And the set 
$$\Ga_1(m,N)\backslash N^{\Ga_1(m,{\mathfrak{a}})}_{GL_m\Z[\xi_3]}(A)/
C_{GL_m \Z[\xi_3]}(A)$$ 
has $\frac{1}{6}\varphi_{\Z[|xi_3]}({\mathfrak{a}})$ elements,
where $\varphi_{\Z[|xi_3]}({\mathfrak{a}})$ is a multiplicative function
on the ideals of $\Z[\xi_3]$ generated by
$$\varphi_{\Z[\xi_3]}((\mathfrak{p})^n)=N_{\Q(\xi_3)/\Q}({\mathfrak{p}})^n(1
-\frac{1}{N_{\Q(\xi_3)/\Q}({\mathfrak{p}})}).$$
  \ele
\proof  In order to save some space set $$\Ga=\Ga_1(m,{\mathfrak{a}})$$ and
$$G=CL_m(\Z[\xi_3]).$$ Let $X\e  N^{\Ga}_G(A)$ and let $B=XAX^{-1}$. Then
$B\e \Ga$. We can write the matrices $A$, $B$ and $X$ in
$2\times 2$ block form with diagonal blocks $A_{11}$, $B_{11}$ and
$X_{11}$ of size $(m-1)\times (m-1)$ and $A_{22}$, $B_{22}$ and
$X_{22}$ of size $1\times1$. Then

 $$
   \left[\begin{tabular}{cc}
$X_{11}$ & $X_{12}$\\ $X_{21}$ & $X_{22}$\\
\end{tabular} \right]
\cdot
  \left[\begin{tabular}{cc}
 $A_{11}$ & $A_{12}$\\
   $0$      & $1$     \\
\end{tabular} \right]
\equiv
  \left[\begin{tabular}{cc}
  $B_{11}$ & $B_{12}$\\
   $0$      & $1$     \\
  \end{tabular} \right]
\cdot
 \left[\begin{tabular}{cc}
  $X_{11}$ & $X_{12}$\\
  $X_{21}$ & $X_{22}$\\
   \end{tabular} \right].
  \mod\: {\mathfrak{a}}
   $$
Consider the $(2,1)$-block of the product. The left hand side is
$$X_{21}A_{11}$$ and the right hand side is 
$$1\cdot X_{21}\: \mod \:{\mathfrak{a}}.$$ 
Examine the map $$P_{B_{22}A_{11}}:X_{21}\mapsto 
X_{21}A_{11}-B_{22}X_{21}$$ with $B_{22}=1$. The eigenvalues of
$A_{11}$ and $B_{22}$ are different. From lemma 1.8 follows that
the eigenvalues of $$P_{B_{22}A_{11}}$$ are $\l_i-1$ where $\l_i$
runs through the eigenvalues of $A_{11}$. From the previous lemma
we have that $\l_i - 1$ is invertible $\mod\:{\mathfrak{a}}$. Thus,
$$P_{B_{22}A_{11}}$$ is non-singular $\mod \: {\mathfrak{a}}$. Therefore,
$$P_{B_{22}A_{11}}(X_{21})= X_{21}A_{11}-B_{22}X_{21} \equiv 0
\mod\: {\mathfrak{a}}$$ implies that $$X_{21}\equiv 0\mod \:{\mathfrak{a}}.$$ 
Therefore,
$$N^{\Ga}_G(A) \subset \Ga_0(m,{\mathfrak{a}}),$$
where $\Ga_0(m,{\mathfrak{a}})$ is defined as the subgroup of $GL_m(\Z[\xi_3])$
that sends $[0,\dots,0,1]$ to $[0,\dots,0,\a]$ modulo ${\mathfrak{a}}$ for any 
$\a\neq 0 \mbox{ }\mod \mbox{ }{\mathfrak{a}}.$

On the other hand, $\Ga_0(m,{\mathfrak{a}})$
lies inside the normalizer of $\Ga_1(m,{\mathfrak{a}})$. 
Also,
$$\Ga_0(m,{\mathfrak{a}})\subset N^{\Ga}_G(A)$$ from the definition of
$N^{\Ga}_G(A)$. We obtain that $$N^{\Ga}_G(A)=\Ga_0(m,{\mathfrak{a}}).$$ 
Let $X$ be
an element in $\Ga_0(m,{\mathfrak{a}})$
of the same block form as
before. Inside the quotient $$\Ga_1(m,{\mathfrak{a}})\backslash 
\Ga_0(m,{\mathfrak{a}})$$ the
element $X$ is determined by $$X_{22} \:\mod \: {\mathfrak{a}}.$$ 
Note that there
are $\varphi_{\Z[\xi_3]} ({\mathfrak{a}})$ options for $X_{22}$ 
because it has to be
invertible modulo ${\mathfrak{a}}$. From the block-triangular theorem (theorem
2.3) we know that we can choose $A$ so that $A_{21}=0$. If $C$ is
in the centralizer of $A$ inside $GL_m(\Z[\xi_3])$, by lemma 3.2 we know
that $C_{21}=0$. Therefore the centralizer of $A$ modulo
$\Ga_1(m,{\mathfrak{a}})$ is determined by $C_{22}$ which could be $+1$,
$-1$, $i$ or $-i$. Therefore the double quotient consist of
$$\frac{1}{4}\varphi_{\Z[\xi_3]}({\mathfrak{a}})$$ elements.


\bth{6.4}
    Let $V$ be a representation of $GL_m(\Z[\xi_3])$, and let
${\mathfrak{a}}$ be an ideal relatively prime to $(1-\xi_3)$. Then
$$
\chi_h(\Ga_1(m,{\mathfrak{a}}),V)=
\varphi_{\Z[\xi_3]}({\mathfrak{a}})\sum_{A=[A_0,1]} 
R(A)\chi(C_{GL_{m-1}\Z[\xi_3]}(A_0))\tr(A^{-1}|V), 
$$
  where $$A_0=[A_{11},\dots A_{kk}],$$
$A_{ii}$ vary through $$\{-1,\xi_3,-\xi_3,\xi_6,-\xi_6\}$$ and $A_{ii}$ and
$A_{jj}$ are distinct,
$\varphi_{\Z[\xi_3]}({\mathfrak{a}})$ is the arithmetic function
defined on the ideals of $\Z[\xi_3]$
generated by 
$$\varphi_{\Z[\xi_3]}({\mathfrak{p}}^n)=N_{\Q(\xi_3)/\Q}({\mathfrak{p}})^{k}(1-
\frac{1}{N_{\Q(\xi_3)/\Q}({\mathfrak{p}})}).$$
       \eth
   \proof Let $A$ be a torsion element of $$\Ga=\Ga_1(m,{\mathfrak{a}})$$
such that $$\chi(C_{\Ga}(A))\neq 0.$$ Then $$\chi(C_{GL_m\Z[\xi_3]}(A))\neq
0.$$ By lemma 6.9 we have that $1$ is an eigenvalue of $A$. Also,
if  $$\chi(C_{\Ga}(A))\neq 0$$ then the multiplicity of the
eigenvalue $1$ is at most $1$. We are going to prove that 
  $$\chi(C_{\Ga}(A))=6\cdot \chi(C_{GL_m\Z[\xi_3]}(A)),$$ 
Assume that we have proven this formula. From the
generalization of Brown's formula, we have
 $$\chi_h(\Ga_1(m,{\mathfrak{a}}),V)=\sum_{A:\mbox{ torsion}}
\chi(C_{\Ga_1(m,{\mathfrak{a}})}(A))\tr(A^{-1}|V).$$ 
In the summation it is enough to
sum over torsion elements $A$ that have eigenvalue $1$ either with
multiplicity $1$. First take the sum over
elements $A$ that become conjugate to each other in $GL_m(\Z[\xi_3])$. 
Let
$A$ has eigenvalue $1$ with multiplicity $1$. By lemma 6.9 $$\sum
\chi(C_{\Ga_1(m,{\mathfrak{a}})} (A''))\tr(A''^{-1}|V)
=\varphi_{\Z[\xi_3]}({\mathfrak{a}})
\chi(C_{GL_n\Z[\xi_3]}(A'))\tr(A'^{-1}|V),$$ where the sum is taken
over all non-conjugate $A''$ in $\Ga_1(m,{\mathfrak{a}})$ that become conjugate
to $A'$ in $GL_n(\Z[\xi_3])$. 

We can assume that $A$ is in block-triangular 
$$
\left[\begin{tabular}{cc}
$A_{11}$ & $A_{22}$ \\
   0     & $1$\\
\end{tabular}\right]
$$
with $A_{22}=1$ and $A_{21}=0$. If $C$ is a matrix in $GL_m(\Z[\xi_3])$
commuting with $A$ then $C$ is if the same block form with
$C_{21}=0$ and $C_{22}=i^k$. Exactly one of the matrices $\xi_6^kC$,
$k=0,1,\dots 5$ belongs to $\Ga_1(m,{\mathfrak{a}})$ 
because of the $C_{22}$ entry. Thus
the centralizer of $A$ in $GL_m(\Z[\xi_3])$ 
contains 6 times as many element
as the centralizer inside $\Ga(m,{\mathfrak{a}})$. Therefore, 
$$\chi(C_{\Ga}(A))=6\cdot \chi(C_{GL_m\Z[\xi_3]}(A)).$$
By the previous lemma we have that the number of non-conjugate matrices in 
$\Ga$ that become conjugate in $G$ is 
$\varphi_{\Z[\xi_3]}({\mathfrak{a}})/6.$
This proves the theorem.
\bth{6.12}
The homological Euler characteristic of 
$\Ga_1(2,{\mathfrak{a}})\subset GL_2(\Z[\xi_3])$ is given by
$$\chi_h(\Ga_1(2,{\mathfrak{a}}),\Q)=
\frac{1}{3}\varphi_{\Z[\xi_3]}({\mathfrak{a}}),$$
where $\varphi_{\Z[\xi_3]}({\mathfrak{a}})$ 
is the multiplicative function defined 
on the ideals of $\Z[\xi_3]$, generated by
$$\varphi_{\Z[\xi_3]}({\mathfrak{p}}^n)=N_{\Q(\xi_3)/\Q}({\mathfrak{p}})^n(1-
\frac{1}{N_{\Q(\xi_3)/\Q}({\mathfrak{p}})}).$$
\eth
\proof From theorem 6.15 and lemma 4.5 we have
$$
\begin{tabular}{rrl}
$\chi_h(\Ga_1(2,{\mathfrak{a}}),\Q)$&=
 &$\varphi_{\Z[\xi_3]}({\mathfrak{a}})(
|N_{\Q(\xi_3)/\Q}(R([1,\xi_6]))|\chi(C([1,\xi_6]))+$\\
\\
&&$+|N_{\Q(\xi_3)/\Q}(R([1,\xi_6^2]))|\chi(C([1,\xi_6^2]))+$\\
\\
&&$+|N_{\Q(\xi_3)/\Q}(R([1,\xi_6^3]))|\chi(C([1,\xi_6^3]))+$\\
\\
&&$+|N_{\Q(\xi_3)/\Q}(R([1,\xi_6^4]))|\chi(C([1,\xi_6^4]))+$\\
\\
&&$+|N_{\Q(\xi_3)/\Q}(R([1,\xi_6^5]))|\chi(C([1,\xi_6^5])))+$\\
\\

 &=&$\varphi_{\Z[\xi_3]}({\mathfrak{a}})(
\frac{1}{36}+\frac{1}{12}+\frac{1}{9}+\frac{1}{12}+\frac{1}{36})=$\\
\\
 &=&$\frac{1}{3}\varphi_{\Z[\xi_3]}({\mathfrak{a}}).$
\end{tabular}
$$

\sectionnew{Dedekind zeta fuction at $-1$
and $\chi_h(SL_2({\cal{O}}_K))$ for totally real number fields}

In this section we examine the homological Euler characteristic of
$SL_2({\cal{O}}_K)$ for $K$ a totally real number field.
At the end of the section we consider two examples for the fields $\Q$ and
$\Q(\sqrt{3})$.
The relation
to the Dedekind zeta function is 
$$\chi(SL_2({\cal{O}}_K))=\zeta_K(-1).$$
We obtain the following theorem.
\bth{8.1}
 Let $K$ be a totally real number field. Then
$$\zeta_K(-1)=
-\frac{1}{4}\sum_{\xi}\sum_{I\in Cl({\cal{O}}_K[\xi]/{\cal{O}}_K)}
\frac{\#{\cal{O}}_K^{\times}/N_{K(\xi)/K}(R_I^{\times})}
{\#(R_I^{\times})_{tors}} +\frac{1}{2}N,$$
where the first summation is taken over all roots $\xi$ of $1$ 
such that $[K(\xi):K]=2,$ 
$Cl({\cal{O}}_K[\xi]/{\cal{O}}_K)$ is the set of ideal classes
 in ${\cal{O}}_K[\xi]$ that are free as ${\cal{O}}_K$-modules,
 the ring $R_I$ is an order inside 
${\cal{O}}_{K(\xi)}$ that contains ${\cal{O}}_K[\xi]$
which is isomorphic to the ring of matrices with coefficients in ${\cal{O}}_K$
that commute with  a matrix $A_I$ corresponding to the ideal class $I$
 (see theorem 1.3).
And $N=\chi_h(SL_2({\cal{O}}_K)).$
\eth
\proof From Brown's theorem we have 
$$N=\chi_h(SL_2({\cal{O}}_K))=\sum_A \chi(C(A)),$$
where the sum is taken over the torsion elements of $SL_2({\cal{O}}_K$
counted up to conjugation.
Note that besides $\pm I_2$ the rest of the torsion elements have
irreducible characteristic polynomial. Their eigenvalues are $\xi$ and
$\xi^{-1}$, where $\xi$ is a root of $1$ such that $[K(\xi):K]=2$.

The centralizers of $I_2$ and of $-I_2$ are the entire group. 
Therefore the Euler characteristic of their centralizer gives two copies
of $\zeta_K(-1).$ 

We need to examine closely the rest of the torsion elements $A$ in 
$SL_2({\cal{O}}_K.$ We have that the eigenvalues of $A$ are
$\xi$ and $\xi^{-1}$, where $$[K(\xi):K]=2.$$ For that reason, 
in the theorem we sum over such $\xi$. By proposition 1.7 
we have that the non-conjugate matrices  with the same characteristic 
polynomial as $A$ are parametrized by $$Cl({\cal{O}}_K[\xi]/{\cal{O}}_K),$$
where the above set denotes the set of ideal classes in ${\cal{O}}_K[\xi]$
that are free as ${\cal{O}}_K$-modules.
For that reason, in the theorem, 
we are summing over ideal classes from this set. 
Fix such a torsion matrix $A$. Let $I$ be the corresponding ideal class
from $Cl({\cal{O}}_K[\xi]/{\cal{O}}_K),$ and $\xi$ and $\xi^{-1}$
be its eigenvalues. We are going to define $R_I.$
Consider  the matrices with coefficients in $K$ that 
commute with $A$. Denote denote the collection by $C_{Mat_2 K}(A).$
We have that $C_{Mat_2 K}(A)\cong K(\xi),$ By sending $A$ to $\xi$,
this isomorphism $\psi$ is an isomorphism of $K$-algebras. Consider
the ring of matrices with coefficients in ${\cal{O}}_K$ that
commute with $A$. Denote it by   $C_{Mat_2 {\cal{O}}_K}(A)$. Define
$R_I=\psi(C_{Mat_2 {\cal{O}}_K}(A)).$ Then for the centralizers in
$GL_2{\cal{O}}_K$ and $SL_2{\cal{O}}_K$ we have
$$R^{\times}=\psi( C_{GL_2{\cal{O}}_K}(A)),$$ and
$$(R_I^{\times})_{tors}=\psi(C_{SL_2{\cal{O}}_K}(A)).$$
The second equality holds because the determinant of a matrix $B$ in 
the centralizer  $C_{GL_2{\cal{O}}_K}(A)$ corresponds to the norm 
$N_{K(\xi)/K}(\beta),$ where $\beta=\psi(B).$ Thus, $\det(B)=1$ 
is equivalent to $N_{K(\xi)/K}(\beta)=1$. The last equality holds 
only when $\beta$ is a root of $1$ because the ranks of the groups of units
in $K$ and in $K(\xi)$ coincide. In Brown's formula we need to compute
$\chi(C(A)),$ which we have done
$$\chi(C_{SL_2{\cal{O}}_K}(A))=\chi((R_I^{\times})_{tors})
=\frac{1}{\# (R_I^{\times})_{tors}}.$$
Proposition 1.3 classifies non-conjugate matrices in $GL_2({\cal{O}}_K).$
In order to pass to $SL_2({\cal{O}}_K),$ we use lemma 6.1.
It gives that the non-conjugate matrices in  $SL_2({\cal{O}}_K)$
that become conjugate to $A$ in $GL_2({\cal{O}}_K)$ are parametrized by
$$SL_2({\cal{O}}_K)\backslash N_{GL_2({\cal{O}}_K)}^{SL_2({\cal{O}}_K)}(A)/
C_{GL_2({\cal{O}}_K)}(A).$$ Note that 
$$ N_{GL_2({\cal{O}}_K)}^{SL_2({\cal{O}}_K)}(A)=GL_2({\cal{O}}_K).$$
Thus, the left quotient 
$SL_2({\cal{O}}_K)(A)\backslash GL_2({\cal{O}}_K)$
is parametrized by 
the determinant of the matrices, namely, by ${\cal{O}}_K^{\times}.$
The group by which we quotient from the right is 
$$C_{GL_2({\cal{O}}_K)}(A)\cong R_I^{\times}.$$
Its determinant leads to $N_{K(\xi)/K}( R_I^{\times})$. Thus, 
the double quotient is isomorphic to 
${\cal{O}}_K^{\times}/N_{K(\xi)/K}( R_I^{\times}),$
which gives the last ingredient in the theorem.
Thus, we have 
$$N=2\cdot \zeta_K(-1)+
\frac{1}{2}\sum_{\xi}\sum_{I\in Cl({\cal{O}}_K[\xi]/{\cal{O}}_K)}
\frac{\#{\cal{O}}_K^{\times}/N_{K(\xi)/K}(R_I^{\times})}
{\#(R_I^{\times})_{tors}},$$
where the $1/2$ occurs because we are summing over all roots of unity
$\xi$ in 
$$\sum_{\xi}\sum_{I\e Cl({\cal{O}}_K[\xi]/{\cal{O}}_K)},$$
while when we are summing over the torsion matrices $A$ we count 
the eigenvalues $\xi$ and $\xi^{-1}$ at the same time. 
\bco{8.2}
Let $K$ be a totally real number field. Suppose that for each root of $1$ $\xi$
such that $[K(\xi):K]=2$,
we have that ${cal{O}}_K[\xi]$ is integrally closed. Then
$$\chi_h(SL_2({\cal{O}}_K),\Q)=2\zeta_K(-1)+\frac{1}{2}\sum_{\xi}C_{\xi},$$
where the sum is taken over all roots of $1$ $\xi$ such that $[K(\xi):K]=2$,
and
$$C_{\xi}=\frac{\#\ker(K_0({\cal{O}}_{K(\xi)})\rightarrow K_0({\cal{O}}_K))
\#\coker(K_1({\cal{O}}_{K(\xi)})\rightarrow K_1({\cal{O}}_K))
}{\#\ker(K_1({\cal{O}}_{K(\xi)})\rightarrow K_1({\cal{O}}_K))},$$
where all the maps are norm maps.
\eco
\proof If $R$ is a number ring we have that $K_0(R)=Cl(R)\oplus \Z$ and
$K_1(R)=R^{\times}$. If ${\cal{O}}_K[\xi]$ is integrally closed 
then it coincides with ${\cal{O}}_{K(\xi)}$.
Then the ideal classes in  ${\cal{O}}_K[\xi]$ that are free 
${\cal{O}}_K$-modules are precisely 
$$\ker(K_0({\cal{O}}_{K(\xi)})\rightarrow K_0({\cal{O}}_K)).$$
Also, the rings $R_I$ from the previous theorem are precisely
${\cal{O}}_{K(\xi)}$. Then
$$\#{\cal{O}}_K^{\times}/N_{K(\xi)/K}({\cal{O}}_{K(\xi)}^{\times})=
\#\coker(K_1({\cal{O}}_{K(\xi)})\rightarrow K_1({\cal{O}}_K)),$$
and
$$\#({\cal{O}}_{K(\xi)}^{\times})_{tors}=
\#\ker(K_1({\cal{O}}_{K(\xi)})\rightarrow K_1({\cal{O}}_K)).$$
From the previous theorem the statement follows.

In the remaining portion of the section we give a simple method for
determining what $R_I$ is. And we end the section with two examples.

Let $A$ be a matrix in $SL_2({\cal{O}}_K)$ that has an irreducible over $K$
characteristic polynomial. Let $\xi$ and $\xi^{-1}$ be its eigenvalues.
By theorem 1.3 we have a correspondence between the ideal classes in 
${\cal{O}}_K[\xi]$ that are free as ${\cal{O}}_K$-modules. We give a way 
of determining what is the order $R_I$. Let $C_{Mat_{2,2}{\cal{O}}_K}(A)$
be the ring of matrices with coefficient in ${\cal{O}}_K$ that
commute with $A$. Let $C_{Mat_{2,2}K}(A)$ be the ring (actually field) 
of matrices that commute with $A$.
Consider that ${\cal{O}}_K$-isomorphism 
$$\psi:C_{Mat_{2,2}K}(A)\rightarrow K(\xi)$$
that send $A$ to $\xi$. Then $R_I$ is defined by
$$R_I=\psi(C_{Mat_{2,2}{\cal{O}}_K}(A)).$$
Every element $\la$ in ${\cal{O}}_{K(\xi)}$ can be written as a linear function 
$f(t)$ with coefficients in $K$, so that $\la=f(\xi)$.

The procedure of determining $R_I$ is the following:
Let $\la\in {\cal{O}}_{K(\xi)}$ and let $\la=f(\xi)$ where $f(t)$ is a 
linear polynomial with coefficients in $K$. Then $\la\e R_I$ if and only if
$f(A)=\psi^{-1}(\l)$ has coefficients in ${\cal{O}}_K$.
\bex{8.2} 
Let $K=\Q$. The root $\xi$ of $1$ that give a quadratic extensions are
$i^{\pm 1}$, $\xi_3^{\pm 1}$ and $\xi_6^{\pm 1}$. For each of the extensions
$Z[\xi]$ is integrally closed. Also ideal class in all the cases are trivial.
Thus, for each $\xi$ we have only one $R_I$ and $R_I=\Z[\xi]$. Then
$$\#\Z^{\times}/N_{\Q(\xi)/\Q}(R_I^{\times})=2.$$ The torsion elements in
$\Z[\xi]$ are $4$ if $\xi=i^{\pm 1}$ and $6$ if $\xi=\pm\xi_3^{\pm 1}$.
Then
$$\chi_h(SL_2\Z,\Q)=2\cdot \zeta(-1)
+\frac{1}{2}(2\cdot\frac{2}{4}+2\cdot\frac{2}{6}+2\cdot\frac{2}{6})=
2\cdot\zeta(-1)+\frac{7}{6}.$$
Since $$\zeta(-1)=-\frac{1}{12},$$ we obtain $$\chi_h(SL_2\Z)=1.$$
\eex
\bex{8.2}
Let $K=\Q(\sqrt{5})$. Then $\xi$ can be $3$-rd, $6$-th, $4$-th, $5$-th or
$10$-th root of unity. 

If $\xi=\xi_{10}^k$ for $k=1,\dots ,4, 6,\dots,9$ then
$${\cal{O}}_K[\xi]=\Z[\xi_{10}],$$
which is integrally closed with class number $1$ (see \cite{W} theorem 11.1). 
Let 
$N=N_{\Q(\xi_{10}/\Q(\sqrt{5}))}$ be the norm map. And let $\xi=\xi_{10}$.
We want to find the image of the norm map  $N$ on the units.
This leads to solving $$N(\a)=\frac{1+\sqrt{5}}{2}.$$
We have
$$N(a+b\xi+c\xi^2+d\xi^3)=$$
$$=(c-a/2)^2+(b-d/2)^2+\frac{3}{4}(a+2d/3)^2+\frac{5}{12}d^2
+(\mbox{ })\frac{1+\sqrt{5}}{2}.$$
If the norm is $(1+\sqrt{5})/2$ then 
all the constants $a$, $b$, $c$, and $d$ must be zero.
Therefore $(1+\sqrt{5})/2$ is not in the image of the norm map and
the cokernel of $N$ is 
$$Z[(1+\sqrt{5})/2]^{\times}/N(\Z[\xi_{10}]^{\times})=C_2\times C_2$$

The kernel of $N$ is $$\Z[\xi_{10}]^{\times}_{tors}=C_{10}.$$
Thus the contribution towards the homological Euler characteristic
is 
$$\frac{1}{2}\cdot 8 \cdot \frac{4}{10}=\frac{8}{5}.$$

If $\xi=\xi_6^k$ for $k=1,2,4,5$ then
$${\cal{O}}_K[\xi]=\Z[(1+\sqrt{5})/2,(1+\sqrt{-3})/2],$$
which is integrally closed. We are going to show that it has class number $1$.
First, we need to examine the units in the ring. Let $N$ be the norm map
from $\Q(\sqrt{5},\sqrt{-3})$ to $\Q(\sqrt{5})$. We try to solve 
$$N(\a)=\frac{1+\sqrt{5}}{2}.$$ We have
$$
N(a+b\frac{1+\sqrt{-3}}{2}+c\frac{1+\sqrt{5}}{2}
+d\frac{1+\sqrt{5}}{2}\cdot\frac{1+\sqrt{-3}}{2})=$$
$$=a^2+ab+b^2+c^2+cd+d^2+(\mbox{  })\frac{1+\sqrt{5}}{2}.$$
Thus $$a^2+ab+b^2+c^2+cd+d^2=0.$$
It can happen only when $a=b=c=d=0$.
Thus, the cokernel of $N$ is isomorphic to $C_2\times C_2.$
Also the regulator for the field $\Q(\sqrt{5},\sqrt{-3})$ is
$$R=\log((1+\sqrt{5})/2).$$ We are going to use the formula for the leading 
coefficient of the Dedekind zeta function at $1$. It is given by
$$\frac{2^{r_1}(2\pi)^{r_2}hR}{w\sqrt{|d|}},$$
where $r_1$ and $r_2$ are the number of the real and the complex imbeddings
if the field $h$ is the cless number, $R$ is the regulator, $w$ 
is the number of roots of unity in the field, and $d$ is the discriminant 
of the field. 

The Dedekind zeta function of $\Q(\sqrt{5},\sqrt{-3})$ decomposes into 
a product of
the Riemann zeta function and three $L$-functions corresponding to the
quadratic extensions $\Q(\sqrt{5})$, $\Q(\sqrt{-3})$ and $\Q(\sqrt{-15})$.
The first two of the quadratic field have class number $1$ and the last 
one class number $2$. Let $h$ be the class number of 
$\Q(\sqrt{5},\sqrt{-3})$. Then
$$\frac{(2\pi)^2h\log((1+\sqrt{5})/2)}{6\cdot 15}=
\frac{2^1\log((1+\sqrt{5})/2)}{2\cdot \sqrt{5}}\cdot
\frac{(2\pi)}{6\cdot \sqrt{3}}\cdot
\frac{(2\pi)\cdot 2}{2\cdot \sqrt{15}},
$$
which gives $h=1$.
The cokernel of $N$ acting on the units is isomporphic to
$C_2\times C_2$ and the 
kernel is isomorphic to $C_6$. Thus the contribution to the homological 
Euler characteristic is 
$$\frac{1}{2}\cdot 4\cdot \frac{4}{6}=\frac{4}{3}.$$

If $\xi=\pm i$ then $${\cal{O}}_K[\xi]=\Z[(1+\sqrt{5},i)],$$
which is integrally closed and of class number $1$. Let $N$ be the norm map 
from $\Q(\sqrt{5},i)$ to $\Q(\sqrt{5})$. We want to solve 
$$N(\a)=\frac{1+\sqrt{5}}{2}.$$ We have
$$N(a+ib+c\frac{1+\sqrt{5}}{2}+di\frac{1+\sqrt{5}}{2})=$$
$$=a^2+b^2+c^2+d^2+(\mbox{  })\frac{1+\sqrt{5}}{2}$$
Thus, $$a^2+b^2+c^2+d^2=0.$$
Therefore $(1+\sqrt{5})/2$ is not in the image of the norm map. 
Thus, the cokernel of $N$ acting on the units in $\Q(\sqrt{5},i)$ is 
isomorphic to $C_2\times C_2$. And the kernel is isomprphic to $C_4$.
Then the contribution to the homological Euler characteristic is
$$\frac{1}{2}\cdot 2\cdot \frac{4}{4}=1.$$
 
Finally,
$$\chi_h(SL_2(\Z[(1+\sqrt{5})/2]))=
2\zeta_{\Q(\sqrt{5})}(-1)+\frac{8}{5}+\frac{4}{3}+1=
2\zeta_{\Q(\sqrt{5})}(-1) +3+\frac{14}{15}.$$
We obtain that $$\zeta_{\Q(\sqrt{5})}(-1)=\frac{1}{30}+\frac{k}{2},$$
for some integer $k$. In fact $$\zeta_{\Q(\sqrt{5})}(-1)=\frac{1}{30},$$
It follows from  \cite{I} theorem 1, which expresses values of 
$L$-functions in terms of generalized Bernouli numbers.
 Therefore
$$\chi_h(SL_2(\Z[(1+\sqrt{5})/2]),\Q)=4.$$
\eex

\sectionnew{Generalization of Brown's formula}

The main result of this section is to prove a generalization of Brown's 
formula relating the homological Euler characteristic to
the Euler characteristics of certain centralizers.
Given an arithmetic group $\Ga$ and a finite dimensional representation
$V$ over a field $K$ of characteristic zero. We prove the following theorem.

\bth{9.1}
The homological Euler characteristic of $\Ga$ with coefficients in
$V$ is given by
$$\chi_h(\Ga,V)=\sum \chi(C(A))\cdot \tr(A^{-1}|V),$$
where the sum is taken over all torsion elements of $\Ga$ counted up 
to conjugation, $C(A))$ is the centralizer of $A$ inside $\Ga$ and
$\tr(A^{-1}|V)$ is the trace of the action of $A^{-1}$ on the finite 
dimensional vector space $V$.
\eth

Before we start with the proof we recall what is complete 
Euler characteristic, and what are some of its properties. For more 
detailed treatment of generalized Euler characteristic we refer to
Bass' article \cite{Ba} and Brown's book \cite{B1}.

Let $R$ be a unital non necessarily commutative ring. Let $F$ be a free 
$R$-module. And let $$f:F\rightarrow F$$ be an endomorphism of $R$-modules.
We can choose a basis for $F$ and then $f$ can be written as a matrix
$(a_{ij})$. We define the trace of $f$ to be 
$$\tr_R(f)=\sum_i \bar{a}_{ii},$$
where $\bar{a}$ is the projection of $a$ from $R$ to $T(R)=R/[R,R]$,
and $[R,R]$ is the additive group generated by $ab-ba$ for $a,b\e R.$
Note that $[R,R]$ is not an ideal so $T(R)$ is only an additive group.

The trace can be extended to endomorphisms of a finitely generated
projective moduleS $P$. Such a module
 can be imbedded in a free module of finite rank $F$. Denote by $i$ 
the imbedding, and by $\pi$ the projection from $F$ to $P$. The the trace 
of $$f:P\rightarrow P$$ if defined by $$\tr_R(f):=\tr_R(i\circ f\circ\pi).$$
Note that $i\circ f\circ\pi$ is an endomorphism of $F$. 

The trace can be also extended to modules $M$ that admit finite length
projective resolution by finitely generated projective modules.
Then an endomorphism  $$f:M \rightarrow M$$ extends to an endomorphism 
of its projective resolution $P_{\bullet}$
$$f_i:P_i\rightarrow P_i.$$
Then the trace of $f$ is defined by
$$\tr_R(f):=\sum_i (-1)^i \tr_R(f_i).$$

Now setting $R$ to be the group ring $\Q\Ga$ the complete Euler 
characteristic is defined by the trace of $id_\Q:\Q\rightarrow \Q,$
namely,
$$E(\Ga):=\tr_{\Q\Ga}(id_{\Q}).$$
Note that the complete Euler characteristic takes values in 
$$T(\Q\Ga)=\Q\Ga/[\Q\Ga,\Q\Ga]=\{\sum_{(A)}c_{(A)}\cdot (A)\},$$
where $(A)$ denotes the conjugacy class of $A\in \Ga$.

Let $$E(\Ga)= \sum_{(A)}c_{(A)}\cdot (A).$$
Let $\Ga$ be an arithmetic group.
The relation between $E(\Ga)$,
$\chi(\Ga)$ and $\chi_h(\Ga)$ is the following
$$\chi(\Ga)=c_{(I)},$$
$$\chi_h(\Ga)=\sum_{(A)}c_{(A)},$$
where the sum is taken over all conjugacy classes. In fact this sum is 
finite which follows from Brown's theorem
\bth{9.2}
If $\Ga$ is an arithmetic group and  $$E(\Ga)= \sum_{(A)}c_{(A)}\cdot (A).$$
Then $$a_{(A)}=\chi(C(A)).$$
\eth 

We generalize the complete Euler characteristic in the following way. 
For a finite dimensional representation $V$ of $\Ga$, let 
$$E(\Ga,V)=\tr_{\Q\Ga}(id_V).$$

In order to pass to homological Euler characteristic
we use a lemma due to Chiswell (see \cite{Ch}, lemma 12)

\ble{9.3}
If 
$$E(\Ga,V)=\sum_{(A)}c_{(A)}(A)$$
then
$$\chi_h(\Ga,V)=\sum_{(A)}c_{(A)}.$$
\ele
\proof
Let 
$$\dots \rightarrow P_1 \rightarrow P_0 \rightarrow V\rightarrow 0$$
be a projective resolution of $V$ of finite length by finitely generated 
projective $\Q\Ga$-modules. Then
$$E(\Ga,V)=\sum_i (-1)^i E(\Ga,P_i).$$
Tensor the projective resolution with $\Q\otimes_{\Q\Ga}.$
Then the augmentation map $\Q\Ga \rightarrow \Q$ induces a map between
$T(\Q\Ga)$ and $T(\Q)$. This map sends 
$$E(\Ga,V)=\tr_{\Q\Ga}(id_V)$$
to
$$\tr_\Q (id_{\Q\otimes_{\Q\Ga}V})=\sum_i (-1)^i 
\tr(id_{(\Q\otimes_{\Q\Ga}P_i)})
=\sum_i (-1)^i H^i(\Ga,V)$$
The augmentation map induced on $T(\Q\Ga)$ sends
$\sum_{(A)}c_{(A)}(A)$ to $\sum_{(A)}c_{(A)}$. And that proves the lemma. 

We are going to use another lemma by Chiswell (see \cite{Ch} lemma 2).
\ble{9.4}
If $f:R_1 \rightarrow R_2$ is a unital ring homomorphism. Then $f$ induces
a map $\bar{f}:T(R_1)\rightarrow T(R_2).$ If $V$ is a finitely generated
projective $R_1$-module then
$$\bar{f}(\tr(id_V))=\tr(id_{(R_2\otimes_{R_1}V)}).$$
\ele

Now we proceed with the proof of theorem 9.1.

\proof (of theorem 8.1) Let $\Ga$ be an arithmetic group, and let $V$
 be a finite dimensional representation of $\Ga$ over a field $K$ of 
characteristic $0$. In stead of the group ring $\Q\Ga$ we shall consider
a much bigger ring $R$. We are going to compute traces with respect to 
this ring $R$. And at the very end we shall compare the final result to traces
over the group ring $\Q\Ga$.

Let $$R_0=\{\sum_A c_A\cdot A|c_a\in K, A\in \Ga\},$$
where sum over infinitely many elements $A$ is allowed, and $K$ is
a field of characteristic zero. We define formal 
variables that correspond to infinite summation. Let
$$x=\sum_{A\in\Ga}1.$$
Denote by $(A)$ the conjugacy class of $A$. Denote by
$$x_{(A)}=\sum_{B\in(A)}1.$$
Let $x$ and $x_{(A)}$ be variables that commute with the elements of $R_0$.
Let $$R:=R_0[x,x^{-1},x_{(A)},\dots]$$ be the ring obtained from $R_0$ by 
adjoining central variables $x$, $x^{-1}$ and $x_{(A)}$ for all conjugacy 
classes $(A)$. Note that we allow to use only finitely many of these 
formulas at a time. Using infinitely many of them would lead to a 
contradiction. Another remark is that if any of the summation is finite , 
we can use the finite number.
 Let 
$$V_R=R\otimes_{K\Ga}V.$$
Let $$F=R\otimes_K V.$$
We are going to define an inclusion $$i:V_R\rightarrow F,$$
and a projection $$\pi:F\rightarrow V_R.$$
Let $v$ be an element of $V$. We define 
$$i(v)=\frac{1}{x}\sum_{A\e\Ga} A\cdot (A^{-1}v),$$ 
And extend it by $R$-linearity to $V_R$.
Let $$\pi(r\otimes v)=rv,$$ where $r\e R$ and $v\e V$
Then $$\pi\circ i=id_{V_R}.$$ Note that $F$ is a free $R$-module
of finite rank
 and $V_R$ is a projective module. By the definition of the trace, 
we have
$$\tr_R(id_{V_R})=\tr_R(i\circ\pi)=\sum_{(A)}
\frac{x_{(A)}}{x}\tr(A^{-1}|V)\cdot (A).$$
In order to find the right regularization of $x_{(A)}/x$, we compare
this result to the the Brown formula for the trivial representation. 
We have a natural inclusion $$\Q\Ga \rightarrow R.$$
Then this inclusion induces equality between 
$\tr_{\Q\Ga}(id_\Q)$ and $\tr_R(id_{\Q_R})$.
Using Brown's formula, we obtain that $$\frac{x_{(A)}}{x}=\chi(C(A)).$$
We also have that $\chi(C(A))=0$ for non-torsion elements $A$. Thus, 
the formula for the trace over $R$ leads to 
$$E(\Ga,V)=\sum_{A:\mbox{ torsion}}\chi(C(A))\tr(A^{-1}|V)\cdot (A).$$
Using Chiswell's lemma we obtain the formula for the homological Euler
characteristic.

\proof (second proof of theorem 8.1)
This proof follows more colosely the paper of Brown \cite{B2}.
And it turn out that the formula we need is a consequence of
a more general formula of Brown
For a finite length chain $C^{\bullet}$ of $\Q\Ga$-modules that admit 
generalized Euler characteristic, define
$$E(\Ga,C^{\bullet})=\sum_i(-1)^iE(\Ga,C^i).$$
Let $$E(\Ga,C^{\bullet})=\sum_{(A)}c_{(A)}\cdot(A).$$
Define $$e(\Ga,C^{\bullet})=c_{(I)},$$
where $I$ is the identity element in $\Ga$. 
We have that $$e(\Ga,\Q)=\chi(\Ga).$$
The general formula that Brown
 obtains is that
$c_{(A)}$ coincides with the coefficient next to $(I)$ in $\tr_{\Q\Ga}(M_A),$
where $M_A:(C^{\bullet})^A \rightarrow (C^{\bullet})^A$ 
is multiplication by $A$.
Now, let $X$ be a contractable CW complex on which $\Ga$ act 
properly and discontinuously. One can assume also that a cell is 
mapped to itself by an element of $\Ga$ then the entire cell is fixed 
by that element pointwise. One can achieve that by subdivision of the cells.
Then $$E(\Ga,V)=E(\Ga,C^{\bullet}(X,V)),$$
where $C^{\bullet}(X,V)$ is the cochain complex of $X$ with coefficients in 
$V$. Let $f\in C^i(X,V)$ and $A\in \Ga$. let also $\sigma$ be an $i$-th cell
and $v\in V$ so that $f(\sigma)=v$. Then $A$ acts on $f$ by 
$(A\cdot f)(A\sigma)=A\cdot v$. Since we want to find the trace we consider 
the fixed cells inder the action of $A$.
Consider a the following basis for $C^i(X,V)$. Let 
$v_1,\dots,v_n$ be a basis for $V$. Let
$f_{\sigma,v_i}$ be the function that sends the cell $\sigma$ to the vector
$v_i$, and sends all other cell to zero. Consider the contribution of 
 $f_{\sigma,v_i}$ to the trace. We have that 
$A\cdot f_{\sigma,v_i}$ maps $\sigma$ to
$a_{ii}\cdot v_i$, where $a_{ii}$ is a constant. But 
$A\cdot f_{\sigma,v_i}$ maps $A\sigma$ to $Av_i$. Thus, we have a 
contribution to the trace if $$A\sigma=\sigma.$$
Consider the space $X^A$ of cells which are fixed under the action of $A$.
On $X^A$ the centralizer $C(A)$ acts propersly and discontinuously.
Also $X^A$ is contractable. Thus 
$$c_{(A)}=e(C(A),C(X^A,V)).$$
We have 
$$(A\cdot f)(\sigma)=A^{-1}(f(A\sigma))=A^{-1}(f(\sigma)).$$
When we take the trace the last equation leads to $\tr(A^{-1}|V)$.
Thus, we obtain 
$$c_{(A)}=e(C(A),C^{\bullet}(X^A))\cdot\tr(A^{-1}|V)
=\chi(C(A))\cdot\tr(A^{-1}|V).$$
\renewcommand{\em}{\textrm}
\begin{small}
\renewcommand{\refname}{ {\flushleft\normalsize\bf{References}} }
    
\end{small}

\begin{thebibliography}{BHY1}
    \bibitem[Ba]{Ba}
Bass, H.: {\em{Euler Characteristic and Characters 
of Discrete Groups,} Inventiones math. 35, 155-196 (1976)}.
   
     \bibitem[B1]{B1}
Brown, K.: {\em{Cohomology of Groups,} Graduate Text in 
Mathematics, Springer-Verlag: New York, 1982}.

    \bibitem[B2]{B2}
Brown, K.: {\em{Complete Euler characteristics and 
fixed-point theory,} J. Pure Appl. Algebra 24 (1982), 103-121.}

   \bibitem[Ch]{Ch}
Chiswell, I.: {\em{Euler characteristics of groups,} 
Math. Z. 147, 1-11(1976)}.

   \bibitem[G1]{G1}
Goncharov, A.: {\em{The dihedral Lie algebras and the
Galois symmetries of $\pi_1(P^1-\{0,\infty\}\cup \mu_n)$,}
Duke Math. J. vol. 110, No. 3 (2001), 397-487.}

   \bibitem[H]{H}
Harder, G.: {\em{A Gauss-Bonnet formula for discrete 
arithmetically defined groups}, 
Ann. Sci. \'Ecole Norm. Sup.(4) 4 (1971), 409-455}.

    \bibitem[I]{I}
Iwasawa, K.: {\em{Lecutures on p-adic L-functions,} Annals 
of Math Studies no, 74. Princeton Univ. Press: Princeton, N.J., 1972.}

    \bibitem[M]{M}
Milnor, J.: {\em{Introduction to algebraic K-theory} Ann. of 
Math Studies 72, Princeton University Press, Princeton, 1971.}

    \bibitem[S]{S}
Serre, J.-P.: {\em{Cohomologie des groupes discretes,}
Ann. of Math. Studies 70 (1971), 77-169.}

\bibitem[W]{W}
Washinton, L.: {\em{Introduction to cyclotomic fields} 
Graduate Text 
in Mathematics, Springer-Verlag: New York, 1982.}
    \end{thebibliography}
\end{document}